\documentclass[11pt,a4paper]{article}
\usepackage{srcltx}
\usepackage{amsmath}
\usepackage{multirow}
\usepackage{slashbox}
\def\Size{{\hbox{\rm Size}}}
\def\mes{{\hbox{\rm mes}}}
\usepackage[bf,normalsize,tableposition=top]{caption}
\usepackage{subfig}
\usepackage{amsfonts}
\usepackage{color}
\usepackage{epsfig}
\usepackage{amssymb}
\usepackage{graphicx}
\usepackage{algorithm}
\usepackage[colorlinks=true]{hyperref}
\def\eexp{{\rm e}}

\newcommand{\cP}{{\cal P}}
\newcommand{\bZ}{{\mathbb{Z}}}

\def\H{{\cal H}}
\def\cI{{\cal I}}
\def\T{{\cal T}}

\def\cV{{\cal V}}

\def\Conv{\hbox{\rm Conv}}
\def\Col{\hbox{\rm Col}}

\def\Opt{{\mathop{\hbox{\rm Opt}}}}

\def\rint{\hbox{\rm rint}\,}
\newcommand{\bbr}{{\mathbb{R}}}
\newcommand{\bbz}{{\mathbb{Z}}}
\def\cD{{\cal D}}

\newcommand{\cN}{I\!\! N}
\def\ex{{\rm e}}

\def\M{{\cal M}}

\def\cP{{\cal P}}
\def\bE{{\mathbf{E}}}
\def\Argmax{{\mathop{\hbox{\rm Argmax}}}}

\newcommand{\F}{{\cal F}}
\newcommand{\E}{{\cal E}}
\newcommand{\C}{{\cal C}}
\newcommand{\D}{{\cal D}}

\newcommand{\I}{{\cal I}}

 \newcommand{\LL}{{\cal L}}
 
  \newcommand{\X}{{\cal X}}

\newcommand{\half}{ \mbox{\small$\frac{1}{2}$}}
\newcommand{\koverT}{ \mbox{\small$\frac{k}{T}$}}
\newcommand{\four}{ \mbox{\small$\frac{1}{4}$}}

\def\eps{\varepsilon}

\def\Ker {{\rm Ker}}

\def\prob{\mathop{\rm Prob}}

\def\inter{{\hbox{\rm int}\,}}

\def\cM{{\cal M}}
\def\Prob{{\hbox{\rm Prob}}}

\newcommand{\be}{\begin{eqnarray}}
\newcommand{\ee}[1]{\label{#1}\end{eqnarray}}
\newcommand{\nn}{\nonumber \\}
\newcommand{\ese}{\end{eqnarray*}}
\newcommand{\bse}{\begin{eqnarray*}}
\newcommand{\rf}[1]{~(\ref{#1})}

\newtheorem{lemma}{Lemma}[section]
\newtheorem{proposition}{Proposition}[section]

\newtheorem{theorem}{Theorem}[section]

\newtheorem{remark}{Remark}[section]

\def\argmin{\mathop{\hbox{\rm argmin$\,$}}}

\def\cU{{\cal U}}
\def\cA{{\cal A}}
\def\cS{{\cal S}}
\def\Card{{\hbox{\rm Card}}}
\def\ErfInv{{\hbox{\rm ErfInv}}}

\def\cN{{\cal N}}
\def\J{{\cal J}}
\def\qed{\hfill$\square$}

\def\mypar{\par}
\def\twohundred{220pt}
\def\onefifty{160pt}
\def\P{{\cal P}}
\def\O{{\cal O}}
\newcommand{\corr}[2]{{\color{blue}#1}}
\newcommand{\hide}[2]{{\color{blue}#1}}

\title{Hypotheses testing by convex optimization}
\author{Alexander Goldenshluger
\thanks{Department of Statistics,
University of Haifa,
31905 Haifa, Israel,
 {\tt goldensh@stat.haifa.ac.il}}
\and
Anatoli Juditsky
\thanks{LJK, Universit\'e Grenoble Alpes, B.P. 53, 38041 Grenoble Cedex 9, France,	
{\tt anatoli.juditsky@imag.fr}}
\and Arkadi Nemirovski
\thanks{Georgia Institute
 of Technology, Atlanta, Georgia
30332, USA, {\tt nemirovs@isye.gatech.edu}\newline
Research of the first author was supported by grants BSF 2010466, and ISF 104/11. The second author was supported by the CNRS-Mastodons project GARGANTUA,
and the LabEx PERSYVAL-Lab (ANR-11-LABX-0025). Research of
the third author was supported by NSF grants CMMI-1232623, CMMI-1262063, CCF-1415498, CCF-1523768.}
}
\begin{document}
\maketitle
\begin{abstract}
We discuss a general approach to hypothesis testing. The main ``building block'' of the proposed construction is a test for a pair of hypotheses in the situation where each particular hypothesis states that
the vector of parameters identifying the distribution of observations  belongs to a convex compact set associated with the hypothesis. This test, under appropriate assumptions,  is {\sl provably nearly optimal} and is yielded by a solution to a convex optimization problem, so that the construction admits computationally efficient implementation. We further demonstrate that our assumptions are satisfied in several important and interesting applications. Finally, we show how our approach can be applied to a rather general testing problems encompassing several classical statistical settings.
\end{abstract}

\section{Introduction}
In this paper we
promote a unified approach  to a class of decision problems, based on Convex Programming.
Our main building block (which we believe is important by its own right)  is a construction, based on Convex Programming (and thus computationally efficient) allowing, under appropriate assumptions,
to build  a {\sl provably nearly optimal} test for deciding between a pair of composite hypotheses on the distribution of observed random variable. Our approach is applicable in several  important situations, primarily, those when observation (a) comes from Gaussian distribution on $\bbr^m$ parameterized by its expectation, the covariance matrix being once for ever fixed, (b) is an $m$-dimensional vector with independent Poisson entries, parameterized by the collection of intensities of the entries, (c) is a randomly selected point from a given $m$-point set $\{1,...,m\}$, with the straightforward parametrization of the distribution by the vector of probabilities for the observation to take values $1$,..., $m$,
(d) comes from a ``direct product
of the outlined observation schemes,'' e.g., is a collection of $K$ independent realizations  of a random variable described by (a)-(c). In contrast to rather restrictive assumptions on the families of distributions we are able to handle, we are very flexible as far as the hypotheses are concerned: all we require  from a hypothesis is to correspond to a convex and compact set in the ``universe'' $\M$ of parameters of the family of distributions we are working with. \mypar
As a consequence, the spirit of the results to follow is quite different from that of a
``classical'' statistical inquiry, where one assumes that the signals underlying noisy observations belong to some ``regularity classes'' and the goal is to characterize analytically the minimax rates of detection for those classes. With our approach allowing for highly diverse hypotheses, an attempt to describe analytically the quality of a statistical routine seems to be pointless. For instance, in the two-hypotheses case, all we know in advance  is that the test yielded by our construction, assuming the latter applicable, is provably nearly optimal, with explicit specification of what ``nearly'' means presented in Theorem \ref{the1}.ii. By itself, this ``near optimality'' usually is not all we need --- we would like to know what actually are the performance guarantees (say, probability of wrong detection, or the number of observations sufficient to make an inference satisfying given accuracy and/or reliability specifications). The point is that with our approach, rather detailed information of this sort can be obtained by efficient situation-oriented computation.
In this respect our approach follows the one of \cite{Lecam1973,Lecam1975,Birge1980,Birge1982,Birge1983M,Lecam1986} where what we call below ``simple tests'' were used to test composite hypotheses represented by convex sets of distributions\footnote{These results essentially cover what in the sequel is called ``Discrete case,'' see section \ref{sect:appl0} for more detailed discussion.}; later this approach was successfully applied to nonparametric estimation of signals and functionals \cite{Birge1983,Donoho1987,Donoho1991,Birge2006}.
On the other hand, what follows can be seen as a continuation of another line of research focusing on testing \cite{Burnashev1979,Burnashev1982,Ingster2002} and on a closely related problem of estimating linear functionals \cite{ibragimov1985,ibragimov1988,Donoho1994} in white noise model. In the present paper we propose a general framework which mirrors that of \cite{JN2009}. Here the novelty (to the best of our understanding, essential) is in applying techniques of the latter paper to hypotheses testing  rather than to estimating linear forms, which allows to naturally encompass and extend the aforementioned approaches to get provably good tests for observations schemes mentioned in (a) -- (d).
We strongly believe that this approach allows to handle a diverse spectrum of applications, and in this paper our focus is on efficiently implementable testing routines\footnote{For precise definitions and details on efficient implementability, see, e.g., \cite{BN2001}. For the time being, it is sufficient to assume that the test statistics can be computed by a simple Linear Algebra routine with parameters which are optimal solutions to an optimization problem which can be solved using {\tt CVX} \cite{cvx2014}.} and related elements of the ``calculus of tests''.
\par
The contents and organization of the paper are as follows. We start with near-optimal testing of pairs of hypotheses, both in its general form and for particular cases of (a) -- (d) (section \ref{sec:mainres}).
We then demonstrate (section \ref{sec:multhyp}) that our  tests (same as other tests of similar structure) for deciding on {\sl pairs} of hypotheses are well suited for ``aggregation,'' via Convex Programming and simple Linear Algebra, into tests with efficiently computable performance guarantees deciding on $M
\geq2$ composite hypotheses. In the concluding section \ref{sec:cases} our focus is on applications. Here we illustrate the implementation of the approaches developed in the preceding sections by building  models and carrying out numerical experimentation for several statistical problems including Positron Emission Tomography, detection and identification of signals in a convolution model, Markov chain related inferences, and some others.

In all  experiments optimization was performed using {\tt Mosek} optimization software \cite{mosek}.
The proofs  missing in the main body of the paper can be found in the appendix.

\section{Situation and Main result}\label{sec:mainres}
In the sequel, given a parametric family $\cP=\{P_\mu,\mu\in\M\}$ of probability distributions on a space $\Omega$ and an observation $\omega\sim P_\mu$ with unknown $\mu\in\M$, we intend to test some composite hypotheses  about the parameter $\mu$. In the situation to be considered in this paper, provably near-optimal testing reduces to Convex Programming, and we start with describing this situation.
\subsection{Assumptions and  goal}\label{sect:Goal}
In what follows, we make the following assumptions on our ``observation environment:''
\begin{quote}
\begin{enumerate}
\item $\M\subset\bbr^m$ is a convex set which coincides with its relative interior;
\item $\Omega$ is a Polish (i.e., separable complete metric) space equipped with a Borel $\sigma$-additive $\sigma$-finite measure $P$, $\hbox{\rm supp}(P)=\Omega$, and distributions $P_\mu\in\cP$ possess densities $p_\mu(\omega)$ w.r.t. $P$. We assume that
    \begin{itemize}
    \item $p_\mu(\omega)$ is continuous in $\mu\in\M$, $\omega\in\Omega$ and is positive;
    \item the densities $p_\mu(\cdot)$ are ``locally uniformly summable:'' for every compact set $M\subset\M$,  there exists a Borel function $p^M(\cdot)$ on $\Omega$ such that
    $\int_\Omega p^M(\omega)P(d\omega)<\infty$ and $p_\mu(\omega)\leq p^M(\omega)$ for all $\mu\in M$, $\omega\in\Omega$;
    \end{itemize}
\item We are given a finite-dimensional linear space $\F$ of continuous functions on $\Omega$ containing constants such that $\ln(p_\mu(\cdot)/p_\nu(\cdot))\in\F$ whenever $\mu,\nu\in\M$.
    \par Note that the latter assumption implies that  distributions $P_\mu,\;\mu\in\M$, belong to an exponential family.
\item
For every $\phi\in \F$, the function
$
F_\phi(\mu)=\ln\left(\int_\Omega\exp\{\phi(\omega)\}p_\mu(\omega)P(d\omega)\right)
$
is well defined and concave in $\mu\in \M$.
 \end{enumerate}
 \end{quote}
In the just described situation, where assumptions 1-4 hold, we refer to the collection $\O=((\Omega,P),\{p_\mu(\cdot):\mu\in\M\},\F)$  as  {\em good observation scheme}.

\par
Now suppose that, on the top of a good observation scheme, we are given
two nonempty convex compact sets $X\subset \M$, $Y\subset \M$. Given an observation $\omega\sim P_\mu$
with some unknown $\mu\in\M$ {\sl known to belong either to $X$} (hypothesis $H_X$) {\sl or to $Y$} (hypothesis $H_Y$), our goal is to decide which of the two hypotheses takes place.
Let $T(\cdot)$ be a test, i.e. a Borel function on $\Omega$
taking values in $\{{-1},1\}$, which receives on input an observation $\omega$ (along with the data participating in the description of $H_X$ and $H_Y$).
Given observation $
\omega$, the test accepts $H_X$ and rejects $H_Y$ when $T(\omega)=1$, and accepts $H_Y$ and rejects $H_X$ when $T(\omega)=-1$. The quality of the test
is characterized by its error probabilities -- the probabilities of rejecting erroneously each of the hypotheses:
\[
 \epsilon_X=\sup_{x\in X}P_x\{\omega:T(\omega)= -1\},
\;\;\;\epsilon_Y=
\sup_{y\in Y} P_y\{\omega:T(\omega)= 1\},
\]
and we define the {\em risk of the test} as the maximal error probability:
$
\max\left\{
\epsilon_X, \epsilon_Y\right\}.
$
\par
In the sequel, we focus on {\sl simple} tests. By definition, a simple test is specified by a {\sl detector} $\phi(\cdot)\in\F$; it accepts $H_X$, the observation being $\omega$, if $\phi(\omega)\geq0$, and accepts $H_Y$ otherwise.
We define the {\sl risk}  of a detector $\phi$ on $(H_{X},H_{Y})$ as the smallest $\epsilon$ such that
\begin{equation}\label{equality}\begin{array}{rl}
\int_\Omega \exp\{-\phi(\omega)\}p_x(\omega) P(d\omega) \leq\epsilon\;\forall x\in X,\,\,\,
\int_\Omega \exp\{\phi(\omega)\}p_y(\omega) P(d\omega) \leq\epsilon\;\forall y\in Y.
\end{array}
\end{equation}
For a simple test with detector $\phi$ we have
\[
 \epsilon_X=\sup_{x\in X} P_x\{\omega:\phi(\omega)<0\},
\;\;\;\epsilon_Y=
\sup_{y\in Y} P_y\{\omega:\phi(\omega)\ge  0\},
\]
and the risk $\max\{\epsilon_X,\epsilon_Y\}$ of such test  clearly does not exceed the risk $\epsilon$ of the detector $\phi$.
\subsection{Main result}
We are about to show that in the situation in question, {\sl an efficiently computable via Convex Programming detector results in a nearly optimal test.} The precise statement is as follows:
\begin{theorem}\label{the1} In the just described situation and under the above assumptions,\\
\mypar {\rm (i)} The function
\begin{equation}\label{Phi}
\begin{array}{l}
\Phi(\phi,[x;y])=\ln\left(\int_\Omega\exp\{-\phi(\omega)\}p_{x}(\omega)P(d\omega)\right)+
\ln\left(\int_\Omega\exp\{\phi(\omega)\}p_{y}(\omega)P(d\omega)\right):\\
\multicolumn{1}{r}{ \F\times(X\times Y)\to\bbr.}\\
\end{array}
\end{equation}
is continuous on its domain, is convex in $\phi(\cdot)\in\F$, concave in  $[x;y]\in X\times Y$, and possesses a saddle point ($\min$ in $\phi$, $\max$ in $[x;y]$)
$(\phi_*(\cdot),[x_*;y_*])$ on $\F\times (X\times Y)$. $\phi_*$ w.l.o.g. can be assumed to satisfy the relation\footnote{Note that $\F$ contains constants, and shifting by a constant the $\phi$-component of a saddle point of $\Phi$ and keeping its $[x;y]$-component intact, we clearly get another saddle point of $\Phi$.}
\begin{equation}\label{balance}
\int_\Omega\exp\{-\phi_*(\omega)\}p_{x_*}(\omega)P(d\omega)= \int_\Omega\exp\{\phi_*(\omega)\}p_{y_*}(\omega)P(d\omega).
\end{equation}
Denoting the common value of the two quantities in {\rm (\ref{balance})} by $\varepsilon_\star$, the saddle point value
\[
\min_{\phi\in\F}\max_{[x;y]\in X\times Y}\Phi(\phi,[x;y])
 \]is $2\ln(\varepsilon_\star)$, and the risk of the simple test associated with the detector $\phi_*$ {on the composite hypotheses $H_X$, $H_Y$} is $\leq\varepsilon_\star$. Moreover, for every $a\in\bbr$, for the test with the detector $\phi_*^a(\cdot)\equiv \phi_*(\cdot)-a$, the probabilities $\epsilon_X$ to reject $H_X$ when the hypothesis is true and $\epsilon_Y$ to reject $H_Y$ when the hypothesis is true can be upper-bounded as
\begin{equation}\label{upperbounded}
\epsilon_X\leq\exp\{a\}\varepsilon_\star,\,\,\epsilon_Y\leq\exp\{-a\}\varepsilon_\star.
\end{equation}
{\rm (ii)}  Let  $\epsilon\geq0$ be such that there exists a (whatever) test for deciding between two simple hypotheses
\begin{equation}\label{AB}
\begin{array}{ll}
(A): \omega\sim p(\cdot):=p_{x_*}(\cdot),&
(B): \omega\sim q(\cdot):=p_{y_*}(\cdot)\\
\end{array}
\end{equation}
with the sum of error probabilities $\leq2\epsilon$. Then
\[
\varepsilon_\star\leq 2\sqrt{\epsilon(1-\epsilon)}.
 \]
 In other words, if the simple hypotheses $(A)$, $(B)$ can be decided, by a whatever test, with the sum of error probabilities $2\epsilon$, then the risk of the simple test with detector $\phi_*$ {on the composite hypotheses $H_X$, $H_Y$} does not exceed $2\sqrt{\epsilon(1-\epsilon)}$.
\mypar {\rm (iii)} The detector $\phi_*$ specified in (i) is readily given by the $[x;y]$-component $[x_*;y_*]$ of the associated saddle point of $\Phi$, specifically,
\begin{equation}\label{phistar}
\phi_*(\cdot)=\half\ln\left(p_{x_*}(\cdot)/p_{y_*}(\cdot)\right).
\end{equation}
\end{theorem}
\paragraph{Remark.} At this point let us make a small summary of the properties of simple tests in the problem setting and under assumptions of section \ref{sect:Goal}:
\begin{itemize}
\item[(i)] One has
\[
\varepsilon_*=\exp(\Opt/2)=\rho(x_*,y_*),
\]
where $[x_*;y_*]$ is the $[x;y]$-component of the saddle point solution of \rf{Phi}, and
\[
\rho(x,y)=\int_\Omega \sqrt{p_{x}(\omega)p_y(\omega)} P(d\omega),
\]
is the {\em Hellinger affinity} of distributions $p_x$ and $p_y$ \cite{Lecam1970,Lecam1986};
\item[(ii)] the optimal detector
$\phi_*$ as in \rf{phistar} satisfies \rf{equality} with $\epsilon=\varepsilon_*$;
 \item[(iii)] the simple test with detector  $\phi_*$ can be ``skewed'', by using instead of $\phi_*(\cdot)$ detector  $\phi^a_*(\cdot)=\phi_*(\cdot)-a$, to attain error probabilities of the test $\epsilon_X=e^a\varepsilon_*$ and $\epsilon_Y=e^{-a}\varepsilon_*$.
     \end{itemize}
As we will see in an instant, the properties (i) -- (iii) of simple tests allow to ``propagate'' the near-optimality property of the tests in the case of repeated observations and multiple testing, and underline all further developments.
\par
Of course, the proposed setting and construction of simple test are by no means unique. For instance, any test $\overline{T}$ in the problem of deciding between $H_X$ and $H_Y$,  with the risk bounded with $ \bar{\epsilon}\in(0,1/2)$, gives rise to the detector
\[
\bar{\phi}(\omega)=\half\ln\left({1-\bar{\epsilon}\over\bar{\epsilon}}\right)\overline{T}(\omega)
\]
(recall that $\overline{T}(\omega)=1$ when $\overline{T}$, as applied to observation $\omega$, accepts $H_X$, and $\overline{T}(\omega)=-1$ otherwise).
One can easily see that the risk of $\bar{\phi}(\cdot)$ satisfies the bounds of \rf{equality}
with
\[
\epsilon=2\sqrt{\bar{\epsilon}(1-\bar{\epsilon})}.
 \]
 In other words, in the problem of deciding upon $H_X$ and $H_Y$, any test $\overline{T}$ with the risk $\leq \bar{\epsilon}$   brings about a simple test with detector $\bar{\phi}$, albeit with a larger risk $\epsilon$.
\subsection{Basic examples}\label{sect:appl0}

We list here some situations where our assumptions are satisfied and thus Theorem \ref{the1} is applicable.
\subsubsection{Gaussian observation scheme}\label{sect:Gauss} In the {\em Gaussian observation scheme} we are given an observation
 $\omega
\in \bbr^m,\; \omega\sim \cN (\mu,\Sigma)$ with unknown parameter $\mu\in \bbr^m$ and known covariance matrix $\Sigma$.
Here the family $\cP$ is defined with $(\Omega,P)$ being $\bbr^m$ with the Lebesque measure, $p_\mu=\cN(\mu,\Sigma)$, $\M=\bbr^m$,
and $\F=\{\phi(\omega)=a^T\omega+b:\;a\in\bbr^m,\;b\in\bbr\}$ is the space of all affine functions on $\bbr^m$. Taking into account that
\[
\ln\left(\int_{\bbr^m}\ex^{a^T\omega+b}p_\mu(\omega)d\omega)\right)={b+a^T\mu+\half{a^T\Sigma a}},
\]
we conclude that Gaussian observation scheme is good.
The test yielded by Theorem \ref{the1} is particularly simple in this case: assuming that the nonempty convex compact sets $X\subset\bbr^m$, $Y\subset\bbr^m$ do not intersect\footnote{otherwise $\phi_*\equiv0$ and $\varepsilon_\star=1$, in full accordance with the fact that in the case in question no nontrivial (i.e., with both error probabilities $<1/2$) testing is possible.}, and that the covariance matrix $\Sigma$ of the distribution of observation is nondegenerate, we get
\be
&&\phi_*(\omega)=\xi^T\omega-\alpha,\;\xi=\half\Sigma ^{-1}[x_*-y_*],\;\alpha=\half\xi^T\Sigma^{-1}[x_*+y_*],\nn
&&\eps_\star=\exp\left(- \mbox{\small$\frac{1}{8}$}(x_*-y_*)^T\Sigma^{-1}(x_*-y_*)\right)\nn
&&\left[[x_*;y_*]\in\Argmax_{x\in X,y\in Y}\left[\psi(x,y)=-\four (x-y)^T\Sigma^{-1}(x-y)\right]\right].
\ee{gausscase}
One can easily verify that the error probabilities $\epsilon_X(\phi^*)$ and $\epsilon_Y(\phi^*)$ of the associated simple test do not exceed $\epsilon_*={\hbox{Erf}}\left(\half\|\Sigma^{-1/2}(x_*-y_*)\|_2\right)$, where ${\hbox{Erf}}(s)$ is the error function: \[\hbox{Erf}(t)=(2\pi)^{-1/2}\int_t^\infty \exp\{-s^2/2\}ds.\]
Moreover, in the case in question the sum of the error probabilities of our test is exactly the minimal, over all possible tests, sum of error probabilities when deciding between the simple hypotheses stating that $x=x_*$ and $y=y_*$.
\paragraph{Remarks.} Consider the simple situation where the covariance matrix $\Sigma$ is proportional to the identity matrix: $\Sigma=\sigma^2 I$ (the case of general $\Sigma$ reduces to this ``standard case'' by simple change of variables). In this case, in order to construct the optimal test, one should find the closest in the Euclidean distance points $x_*\in X$ and $y_*\in Y$, so that the affine form $\zeta(u)=[x_*-y_*]^Tu$ strongly separates $X$ and $Y$. On the other hand, testing in the white Gaussian noise between the closed half-spaces $\{u:\;\zeta(u)\leq \zeta(y_*)\}$ and $\{u:\;\zeta(u)\geq \zeta(x_*)\}$ (which contain $Y$ and $X$, respectively)  is exactly the same as deciding on two simple hypotheses stating that $y=y_*$, and $x=x_*$. Though this result is almost self-evident, it seems first been noticed in \cite{Burnashev1979} in the problem of testing in white noise model, and then exploited in \cite{Burnashev1982,Ingster2002} in the important to us context of hypothesis testing.
\par
As far as numerical implementation of the testing routines is concerned, numerical stability of the proposed test is an important issue. For instance, it may be useful to know the testing performance
when the optimization problem \rf{gausscase} is not solved to exact optimality, or when errors may be present in description of the sets $X$ and $Y$. Note that one can easily bound the error of the obtained test in terms of the magnitude of violation of first-order optimality conditions for \rf{gausscase}, which read:
\[
(y_*-x_*)^T\Sigma^{-1}(x-x_*)+(x_*-y_*)^T\Sigma^{-1}(y-y_*)\leq 0,\;\forall x\in X,\;y\in Y.
\]
Now assume that instead of the optimal test $\phi_*(\cdot)$  we have at our disposal an ``approximated'' simple test
associated with
\[
 \tilde{\phi}(\omega)=\tilde{\xi}^T\omega-\tilde{\alpha},\;\tilde{\xi}=
 \half\Sigma^{-1}[\tilde{x}-\tilde{y}],\;\tilde{\alpha}=\half \tilde{\xi}^T[\tilde{x}+\tilde{y}],
\]
where $\tilde{x}\in X,\;\tilde{y}\in Y$, $\tilde{x}\neq \tilde{y}$ satisfy
\be
(\tilde{y}-\tilde{x})^T\Sigma^{-1}(x-\tilde{x})+(\tilde{x}-\tilde{y})^T\Sigma^{-1}(y-\tilde{y})\leq \delta,\;\forall x\in X,\;y\in Y,
\ee{approx:gauss}
with some $\delta>0$. This implies the bound for the risk of the test with detector $\tilde{\phi}(\cdot)$:
\be
\max[\epsilon_X,\epsilon_Y]\leq \tilde{\epsilon}=
{\hbox{Erf}}\left(
\half\|\Sigma^{-1/2}(\tilde{x}-\tilde{y})\|_2-{\delta\over \|\Sigma^{-1/2}(\tilde{x}-\tilde{y})\|_2}\right).
\ee{bound:gauss}
{\small \begin{quotation}
Indeed, \rf{approx:gauss} implies that
$
\tilde{\xi}^T(x-\tilde{x})\geq -{\delta\over 2},\;\;\tilde{\xi}^T(y-\tilde{y})\leq {\delta\over 2},\;\;\forall x\in X, \;y\in Y.$
As a result,
\[
\tilde{\xi}^Tx-\tilde{\alpha}=\tilde{\xi}^T(x-\tilde{x})+\tilde{\xi}^T\Sigma\tilde{\xi}\geq -{\delta\over 2}+\tilde{\xi}^T\Sigma\tilde{\xi}\;\;\forall x\in X.
\]
and for all $x\in X$,
\[
\Prob_x \{\tilde{\phi}(\omega)<0\}=\Prob_x\{\tilde{\xi}^T(\omega-x)<-\tilde{\xi}^Tx+\tilde{\alpha}\}=
\Prob_x\left\{\|\Sigma^{1/2}\tilde{\xi}\|_2\eta<-\|\Sigma^{1/2}\tilde{\xi}\|^2_2+{\delta\over 2}\right\},
\]
where  $\eta\sim \cN(0,1)$. We conclude that
\[
\epsilon_X=\sup_{x\in X}\Prob_x\{\tilde{\phi}(\omega)<0\} \leq {\hbox{Erf}}\left(
\half\|\Sigma^{1/2}\tilde{\xi}\|_2-{\delta\over 2\|\Sigma^{1/2}\tilde{\xi}\|_2}\right)
\]
what implies the bound \rf{bound:gauss} for $\epsilon_X$.
The corresponding bound for $\epsilon_Y=\sup_{y\in Y}\Prob_y\{\tilde{\phi}(\omega)\ge 0\}$ is obtained in the same way.
\end{quotation}}

\subsubsection{Discrete observation scheme}  Assume that we observe a realization of a random variable $\omega$ taking values in $\{1,2,...,m\}$ with probabilities $\mu_i,\;i=1,...,m$:
\[
\mu_i=\prob\{\omega=i\},\;i=1,...,m.
\]
The just described {\em Discrete observation scheme} corresponds to $(\Omega,P)$ being $\{1,...,m\}$ with counting measure, $p_\mu(\omega)=\mu_\omega$,
$\mu\in \M=\{\mu\in \bbr^m:\;\mu_i>0,\;\sum_{i=1}^m\mu_i=1\}$, In this case $\F=\bbr(\Omega)=\bbr^m$, and for $\phi\in \bbr^m$,
\[
\ln\left(\sum_{\omega\in \Omega}e^{\phi(\omega)} p_\mu(\omega)\right)= \ln\left(\sum_{\omega=1}^m e^{\phi_\omega} \mu_\omega\right)
\]
is concave in $\mu\in \M$. We conclude that Discrete observation scheme is good.
Furthermore, when assuming the convex compact sets $X\subset\M$, $Y\subset\M$ (recall that in this case $\M$ is the relative interior of the standard simplex in $\bbr^m$) not intersecting, we get
\begin{equation}\label{discrcase}
\begin{array}{c}
\phi_*(\omega)=\ln\left(\sqrt{[x_*]_\omega/[y_*]_\omega}\right),\,\,
\eps_\star=\exp\{\Opt/2\}=\rho(x_*,y_*),\\
\left[{[x_*;y_*]}\in\Argmax_{x\in X,y\in Y}\left[\psi(x,y)=2\ln\rho(x,y),\;\Opt=\psi(x_*,y_*)\right],\right]
\end{array}
\end{equation}
where $\rho(x,y)=\sum_{\ell=1}^m\sqrt{x_\ell y_{\ell}}$ is the Hellinger affinity of  distributions $x$ and $y$. One has $\eps_\star=\rho(x_*,y_*)=1-h^2(x_*,y_*)$, the Hellinger affinity of the sets $X$ and $Y$, where
\[
h^2(x,y)=\half \sum_{\ell=1}^m \left(\sqrt{x_\ell}-\sqrt{y_\ell}\right)^{2}
 \]is the {\em Hellinger distance} between distributions $x$ and $y$. Thus the result of Theorem \ref{the1}, as applied to Discrete observation model, allows for the following simple interpretation: to construct the simple test $\phi_*$ one should find the closest in Hellinger distance points $x_*\in X$ and $y_*\in  Y$; then the risk of the likelihood ratio test $\phi_*$ for distinguishing $x_*$ from $y_*$, as applied to our testing problem, is bounded with $\rho(x_*,y_*)=1-h^2(x_*,y_*)$, the Hellinger affinity of  sets $X$ and $Y$.
\paragraph{Remarks.}
Discrete observation scheme considered in this section is a simple particular case -- that of finite $\Omega$ -- of the  result of \cite{birge1981,Birge1982} on distinguishing convex sets of distributions.
 Roughly, the situation considered in those papers is as follows: let $\Omega$ be a Polish space, $P$ be a $\sigma$-finite $\sigma$-additive Borel measure on $\Omega$, and  $p(\cdot)$ be a density  w.r.t. $P$ of probability distribution of observation $\omega$. Note that the corresponding observation scheme (with $\M$ being the set of densities with respect to $P$ on $ \Omega$) does not satisfy the premise of section \ref{sect:Goal} because the linear space $\F$ spanned by constants and functions of the form $
 \ln(p(\cdot)/q(\cdot))$, $p,q\in\M$ is not finite-dimensional.
Now assume that we are given two non-overlapping convex closed subsets $X$, $Y$ of the set of probability densities  with respect to $P$ on $\Omega$.
Observe that for every positive Borel function $\psi(\cdot):\Omega\to\bbr$, the detector $\phi$ given by
$
\phi(\omega)=\ln(\psi(\omega))
$
for evident reasons satisfies the relation
\[
\begin{array}{c}
\max\limits_{p\in X,q\in Y}\left[\int_\Omega \eexp^{-\phi(\omega)}p(\omega)P(d\omega),\int_\Omega \eexp^{\phi(\omega)}q(\omega)P(d\omega)\right]\leq\epsilon,\\
\epsilon=\max\left[\sup_{p\in X}\int \psi^{-1}(\omega)p(\omega)P(d\omega),\sup_{q\in Y}\int \psi(\omega)q(\omega)P(d\omega)\right]\\
\end{array}
\]
Let now
\begin{equation}\label{Tolya16n}
\Opt=\max_{p\in X,q\in Y}\left\{\rho(p,q)=\int_\Omega \sqrt{p(\omega)q(\omega)}P(d\omega)\right\},
\end{equation}
which is an infinite-dimensional convex program with respect to $p\in X$ and $q\in Y$. Assuming the program solvable with an optimal solution composed of distribution $p_*(\cdot)$, $q_*(\cdot)$ which are positive, and setting $\psi_*(\omega)=\sqrt{p_*(\omega)/q_*(\omega)}$, under some ``regularity assumptions'' (see, e.g., Proposition 4.2 of \cite{Birge1982}) the optimality conditions for \rf{Tolya16n} read:
\[
\min_{p\in X,q\in Y}\left[\int_\Omega\psi_*^{-1}(\omega)[p_*(\omega)-p(\omega)]P(d\omega)+\int_\Omega\psi_*(\omega)
[q_*(\omega)-q(\omega)]P(d\omega)\right]=0.
\]
In other words,
$$\max_{p\in X}\int_\Omega\psi_*^{-1}(\omega)p(\omega)dP(\omega)\leq  \int_\Omega\psi_*^{-1}(\omega)p_*(\omega)dP(\omega)=\Opt,$$ and similarly,
$$\max_{q\in Y}\int_\Omega\psi_*(\omega)q(\omega)dP(\omega)\leq  \int_\Omega\psi_*(\omega)q_*(\omega)dP(\omega)=\Opt,$$ so that for our $\psi_*$, we have $\epsilon=\Opt$.
\par
Note that, although this approach is not restricted to the Discrete case {\sl per se}, when $\Omega$ is not finite, the optimization problem in \rf{Tolya16n} is generally computationally intractable (the optimal detectors can be constructed explicitly for some special sets of distribution, see \cite{Birge1982,Birge1983M}).
\par
The bound  $\eps_\star$ for the risk of the simple test can be compared to the {\em testing affinity} $\pi(X,Y)$ between $X$ and $Y$,
\[
\pi(X,Y)=\max_{x\in X,y\in Y}\left\{\pi(x,y)=\sum_{\ell=1}^m \min[x_\ell,y_\ell]\right\},
\]
which is the least possible sum of error probabilities $\epsilon_X+\epsilon_Y$ when distinguishing between  $H_X$ and $H_Y$ (cf. \cite{Lecam1973,Lecam1986}).
 The corresponding  {\em minimax test} is a simple test with detector $\overline{\phi}(\cdot,\cdot)$, defined according to
\[\begin{array}{c}
\overline{\phi}(\omega)=\ln\left(\sqrt{[\overline{x}]_\omega/[\overline{y}]_\omega}\right),\\
\left[{[\overline{x};\overline{y}]}\in\Argmax_{x\in X,y\in Y}\left[\sum_{\ell=1}^m \min[x_\ell,y_\ell\right].\right].
\end{array}
\]
Unfortunately, this test cannot be easily extended to the case where repeated observations (e.g., independent realizations $\omega_k$, $k=1,...,K$, of $\omega$) are available. In \cite{huber1973} such an extension has been proposed in the case where $X$ and $Y$ are dominated by bi-alternating capacities (see, e.g., \cite{huber1974,bednarski1982,buja1986huber,augustin2010}, and references therein); explicit constructions of the test were proposed for some special sets of distributions \cite{Huber1965,rieder1977,osterreicher1978}. On the other hand, as we shall see in section \ref{secrepeated}, the simple test $\phi_*(\cdot,\cdot)$ allows for a straightforward generalization to the repeated observations case with the same (near-)optimality guaranties as those of Theorem \ref{the1}.ii.
\par
Finally, same as in the Gaussian observation scheme, the risk of a simple test with detector $\tilde{\phi}(\omega)=\half \ln\left({\tilde{x}_\omega/ \tilde{y}_\omega}\right),\;\omega\in \Omega$, defined by a pair of distributions
$[\tilde{x};\tilde{y}]\in X\times Y$, can be assessed through the magnitude of violation by $\tilde{x}$ and $\tilde{y}$ of the first-order optimality conditions for the optimization problem in \rf{discrcase}. Indeed, assume that
\[
\sum_{\ell=1}^m\sqrt{\tilde{y}_\ell\over \tilde{x}_\ell}(x_\ell-\tilde{x}_\ell)+
\sum_{\ell=1}^m\sqrt{\tilde{x}_\ell\over \tilde{y}_\ell}(y_\ell-\tilde{y}_\ell)\leq \delta\;\;
\forall x\in X,\;y\in Y.
\]
We conclude that
\bse
\epsilon_X&\leq& \max_{x\in X}\sum_{\ell=1}^m e^{-\tilde{\phi}_\ell} x_\ell=\max_{x\in X}\sum_{\ell=1}^m\sqrt{\tilde{y}_\ell\over \tilde{x}_\ell}x_\ell
\leq \sum_{\ell=1}^m\sqrt{\tilde{y}_\ell\tilde{x}_\ell}+\delta,\\
\epsilon_Y&\leq &\max_{y\in Y}\sum_{\ell=1}^m e^{\tilde{\phi}_\ell} y_\ell=\max_{y\in Y}\sum_{\ell=1}^m\sqrt{\tilde{x}_\ell\over \tilde{y}_\ell}y_\ell
\leq \sum_{\ell=1}^m\sqrt{\tilde{x}_\ell\tilde{y}_\ell}+ \delta,
\ese
 so that the risk of the test $\tilde{\phi}$ is bounded with $\rho(\tilde{x},\tilde{y})+\delta$.

\subsubsection{Poisson observation scheme}\label{sect:Poiss}  Suppose that we are given $m$ realizations of independent Poisson random variables
\[ \omega_i\sim \hbox{\rm Poisson}(\mu_i)\]
with parameters ${\mu_i},\;i=1,...,m$. The {\em Poisson observation scheme} is given by
$(\Omega,P)$ being  ${\bbz}_+^m$ with counting measure, $p_\mu(\omega)={\mu^\omega\over\omega!}e^{-\sum_i\mu_i}$ where $\mu\in\M=\inter\,\bbr^m_+$, and, similarly to the Gaussian case, $\F$ is comprised of the restrictions onto $\bZ_+^m$ of affine functions:
{$\F=\{\phi(\omega)=a^T\omega+b:\;a\in\bbr^m,\;b\in\bbr\}$}. Since
\[
\ln\left(\sum\limits_{\omega\in{\bbz}^m_+}
         \exp({a^T\omega+b})p_\mu(\omega)\right)={\sum_{i=1}^m(e^{a_i}-1)\mu_i}+b
\]
         is concave in $\mu$, we conclude that Poisson observation scheme is good.
         \par
Assume now that, same as above,  in the Poisson observation scheme, the convex compact sets $X\subset\bbr^m_{++}$, $Y\subset\bbr^m_{++}$ do not intersect. Then the data associated with the simple test yielded by Theorem \ref{the1} is as follows:
\begin{equation}\label{Poissoncase}
\begin{array}{c}
\phi_*(\omega)=\xi^T\omega-\alpha,\;\xi_\ell=\half \ln\left({[x_*]_\ell/[y_*]_\ell}\right),
\;\;\alpha=\half\sum_{\ell=1}^m[x_*-y_*]_\ell,\,\,\varepsilon_\star=\exp\{\Opt/2\}\\
\left[\begin{array}{rcl}{[x_*;y_*]}&\in&\Argmax_{x\in X,y\in Y} \left[\psi(x,y)=-2h^2(x,y)
\right],\;\;
\Opt=\psi(x_*,y_*),
\end{array} \right]
\end{array}
\end{equation}
where $h^2(x,y)=\half \sum_{\ell=1}^m\left[\sqrt{x_\ell}-\sqrt{y_\ell}\right]^2$ is the Hellinger distance between $x\in \bbr^m_+$ and $y\in \bbr^m_+$.
\paragraph{Remark.}
Let $\tilde{\phi}(\omega)=\tilde{\xi}^T\omega-\tilde{\alpha}$ be a detector, generated by
$[\tilde{x};\tilde{y}]\in X\times Y$, namely, such that
\[
\tilde{\xi}_\ell=\half \ln(\tilde{x}_\ell/\tilde{y}_\ell),\;\;\tilde{\alpha}=\half \sum_{\ell=1}^m (\tilde{x}_\ell-\tilde{y}_\ell).
\]
We assume that $[\tilde{x};\tilde{y}]$ is an approximate solution to \rf{Poissoncase} in the sense that the first-order optimality condition of \rf{Poissoncase} is `$\delta$-satisfied'':
\[
\sum_{\ell=1}^m\left[\left(\sqrt{\tilde{y}_\ell/ \tilde{x}_\ell}-1\right)(x_\ell-\tilde{x}_\ell)+
\left(\sqrt{\tilde{x}_\ell/ \tilde{y}_\ell}-1\right)(y_\ell-\tilde{y}_\ell)\right]\leq \delta\;\;\forall x\in X,\;y\in Y.
\]
One can easily verify that the risk of the test, associated with $\tilde{\phi}$, is bounded with
$
\exp(- h^2(\tilde{x},\tilde{y})+\delta)
$ (cf. the corresponding bounds for the Gaussian and Discrete observation schemes).

{\small\begin{quotation}
\end{quotation}}

\subsection{Repeated observations}\label{secrepeated} Good observation schemes admit naturally defined direct products. To simplify presentation, we start with explaining the corresponding construction in the case of {\sl stationary repeated observations} described as follows.
\subsubsection{$K$-repeated stationary observation scheme}
We are given a good observation scheme $((\Omega,P),\{p_\mu(\cdot):\mu\in\M\},\F)$ and a positive integer $K$, along with same as  above $X,Y$. Instead of a single realization $\omega\sim p_\mu(\cdot)$,  we now observe a sample of $K$ {\sl independent} realizations $\omega_k\sim p_\mu(\cdot)$, $k=1,...,K$. Formally, this corresponds to the observation scheme with the observation space $\Omega^{K}=\{\omega^K=(\omega_1,...,\omega_K):\omega_k\in\Omega\,\forall k\}$ equipped with the measure $P^{K}=P\times...\times P$, the family $\{p^{K}_\mu(\omega^K)=\prod_{k=1}^Kp_\mu(\omega_k),\mu\in\M\}$ of densities of repeated observations w.r.t. $P^{K}$, and $\F^{K}=\{\phi^{K}(\omega^K)=\sum_{k=1}^K\phi(\omega_k),\phi\in\F\}$. The components $X,Y$ of our setup are the same as for the original single-observation scheme, and the composite hypotheses we intend to decide upon state now that the $K$-element observation $\omega^K$ comes from a distribution $p_\mu^{K}(\cdot)$ with $\mu\in X$ (hypothesis $H_X$) or with $\mu\in Y$ (hypothesis $H_Y$).
\mypar
It is immediately seen that the just described {\sl $K$-repeated} observation scheme is good (i.e., satisfies all our assumptions), provided that the ``single observation'' scheme we start with is so. Moreover, the detectors $\phi_*$, $\phi_*^{K}$ and risk bounds $\varepsilon_\star$, $\varepsilon_\star^{(K)}$ given by Theorem \ref{the1} as applied to the original and the $K$-repeated observation schemes are linked by the relations
\begin{equation}\label{precrem}
\begin{array}{c}
\phi^{K}_*(\omega_1,...,\omega_K)={\sum}_{k=1}^K \phi_*(\omega_k),    \,\,\varepsilon_\star^{(K)}=(\varepsilon_\star)^K.
\end{array}
\end{equation}
As a result, the ``near-optimality claim'' Theorem \ref{the1}.ii can be reformulated as follows:
\begin{proposition}\label{newcol1}
Assume that for some integer $\bar{K}\geq1$ and some $\epsilon\in(0,1/4)$, the hypotheses $H_X$, $H_Y$ can be decided, {\sl by a whatever procedure utilising $\bar{K}$ observations}, with error probabilities $\leq \epsilon$. Then  with
\[
K^+=\left\rfloor {2\bar{K}\over 1-{2\ln[2]\over\ln[1/\epsilon]}}\right\lfloor
\]
observations, $\rfloor a\lfloor$ being the smallest integer $\geq a$, the simple test with the detector $\phi_*^{K^+}$ decides between $H_X$ and $H_Y$ with risk $\leq\epsilon$.
\end{proposition}
{\small \begin{quotation}{Indeed, applying (\ref{precrem}) with $K=\bar{K}$ and utilizing Theorem \ref{the1}.ii, we get $\varepsilon_\star\leq
(2\sqrt{\epsilon})^{1/\bar{K}}$ and therefore, by the same (\ref{precrem}),
$\varepsilon_\star^{(K)}=\varepsilon_\star^K\leq (2\sqrt{\epsilon})^{K/\bar{K}}$ for all $K$. Thus, $\varepsilon_\star{(K^+)}\leq\epsilon$, and therefore the conclusion of Proposition follows from Theorem \ref{the1}.i as applied to observations $\omega^{K^+}$.}
\end{quotation}
}
We see that for small $\epsilon$, the ``suboptimality ratio'' (i.e., the ratio $K^+/\bar{K}$) of the proposed test when $\epsilon$-reliable testing is sought is close to 2 for small $\epsilon$.
\subsubsection{Non-stationary repeated observations}\label{secrepeatedI}
We are about to define the notion of a general-type direct product of good observation schemes. The situation now is as follows: we are given $K$ good observation schemes
$$\O_k=\left((\Omega_k, P_k),\M_k\subset\bbr^{m_k},\{p_{k,\mu_k}(\cdot):\mu_k\in\M_k\},\F_k\right),\,k=1,...,K$$
and observe a sample $\omega^K=(\omega_1,...,\omega_K)$ of  realizations $\omega_k\in\Omega_k$ drawn independently of each other from the distributions with densities, w.r.t. $P_k$, being $p_{k,\mu_k}(\cdot)$, for a collection $\mu^K=(\mu_1,...,\mu_K)$ with $\mu_k\in \M_k$, $1\leq k\leq K$. Setting
\bse
&&\Omega^K=\Omega_1\times...\times\Omega_K=\{\omega^K=(\omega_1,...,\omega_K):\omega_k\in\Omega_k\,\forall k\leq K\},\\
&&P^K=P_1\times...\times P_K\\
&&\M^K=\M_1\times...\times \M_K=\{\mu^K=(\mu_1,...,\mu_K):\mu_k\in \M_k\,\forall k\leq K\},\\
&&p_{\mu^K}(\omega^K)=p_{1,\mu_1}(\omega_1)p_{2,\mu_2}(\omega_2)...p_{K,\mu_K}(\omega_K)\;\;\;\;[\mu^K\in\M^K,\omega^K\in\Omega^K],\\
&&\F^K=\{\phi^K(\omega^K)=\phi_1(\omega_1)+\phi_2(\omega_2)+...+\phi_K(\omega_K):\;\Omega^K\to\bbr:\phi_k(\cdot)\in \F_k\,\forall k\leq K\},\\
\ese
we get an observation scheme $((\Omega^K,P_K),\M^k,\{p_{\mu^K}(\cdot):\mu^k\in\M^k\},\F^K)$  which we call the {\em direct product of $\O_1,...,\O_K$} and denote $\O^K=\O_1\times...\times \O_K$. It is immediately seen that this scheme is good. Note that the already defined stationary repeated observation scheme deals with a special case of the direct product construction, the one where all factors in the product are identical to each other, and where, in addition, we replace $\M^K$ with its ``diagonal part'' $\{\mu^K=(\mu,\mu,...,\mu),\,\mu\in\M\}$.
 \par
 Let $\O^K=\O_1\times...\times \O_K$, where, for every $k\leq K$,
  {$$\O_k=((\Omega_k,P_k),\M_k,\{p_{\mu_k}(\cdot):\mu_k\in\M_k\},\F_k)$$}
  is a good observation scheme, specifically, either Gaussian, or Discrete, or Poisson (see section \ref{sect:appl0}).
 To simplify notation, we assume that all Poisson factors $\O_k$ are ``scalar,'' that is, $\omega_k$ is drawn from Poisson distribution with parameter $\mu_k$.\footnote{This assumption in fact does not restrict generality, since an $m$-dimensional Poisson observation scheme from section \ref{sect:Poiss} is nothing but the direct product  of $m$ scalar Poisson observation schemes. Since the direct product of observation schemes clearly is associative, we always can reduce the situation with multidimensional Poisson factors to the case where all these factors are scalar ones.}
For
$$
\phi^K(\omega^K)=\sum_{k=1}^K\phi_k(\omega_k)\in\F^K,\,\,\mu^K=(\mu_1,...,\mu_K)\in\M^K,
$$
let us set
\[
\Psi(\phi^K(\cdot),\mu^K)=\ln\left(\int_{\Omega^K}\exp\{-\phi^K(\omega^K)\}p_{\mu^K}(\omega^K)P^K(d\omega^K)\right)
=\sum_{k=1}^K\Psi_k(\phi_k(\cdot),\mu_k),
\]
with
\[
\Psi_k(\phi_k(\cdot),\mu_k)=\ln(\left(\int_{\Omega_k}\exp\{-\phi_k(\omega_k)\}p_{k,\mu_k}(\omega_k)P_k(d\omega_k)\right).
\]
The function $\Phi(\phi^K,[x,y])$, defined by (\ref{Phi}) as applied to the observation scheme $\O^K$, clearly is
\bse
&&\Phi(\phi^K,[x;y])=\sum_{k=1}^K \left[\Psi_k(\phi_k,x_k)+\Psi_k(-\phi_k,y_k)\right],\\
&&\left[\phi^K(\omega^K)=\sum_k\phi_k(\omega_k),\;x=[x_1;...;x_k]\in\M^K,\;\;\;y=[y_1;...;y_K]\in\M^K\right]
\ese
so that
\[
\min\limits_{\phi^K\in\F^K}\Phi(\phi^K,[x;y])=\sum_{k=1}^K \psi_k(x_k,y_k),
\]
where functions $\psi_k(\cdot,\cdot)$ are defined as follows (cf. \rf{gausscase}, \rf{discrcase} and \rf{Poissoncase}):
\begin{itemize}
\item $
\psi_k(\mu_k,\nu_k)=-\four{(\mu_k-\nu_k)^T\Sigma_k^{-1}(\mu_k-\nu_k)}
$
in the case of Gaussian $\O_k$ with $\omega_k\in \bbr^{m_k}$, $\omega_k\sim\cN(\mu_k,\Sigma_k)$, $\mu_k,\nu_k\in\bbr^{m_k}$;
\item $ \psi_k(\mu_k,\nu_k)=-(\sqrt{\mu_k}-\sqrt{\nu_k})^2$ for scalar Poisson $\O_k$,
with $\mu_k,\nu_k>0$;
\item
$ \psi_k(\mu_k,\nu_k)=2\ln\left(\sum_{i=1}^{m_k}\sqrt{[\mu_k]_i[\nu_k]_i}\right)$ for Discrete $\O_k$ with $\Omega_k=\{1,...,m_k\}$,\\ $\mu_k,\nu_k\in\M_k=\left\{\mu\in\bbr^{m_k}:\mu>0,\;\sum_i[\mu]_i=1\right\}$.
\end{itemize}
Let $X_k$ and $Y_k$ be compact convex subsets of $\M_k$, $k=1,...,K$; let $X=X_1\times ...\times X_K$ and $Y=Y_1\times ...\times Y_K$. Assume that
$[x_*;y_*]=\left[
[x_*]_1;...;[x_*]_K;[y_*]_1;...;[y_*]_K\right]$ is an optimal solution to the convex optimization problem
\be
\Opt=\max\limits_{x\in X,y\in Y}\left[
\sum_{k=1}^K
\psi_k(x_k,y_k)\right],
\ee{mult_opt}
and let
\be
\phi^k_*(\omega_k)=\left\{\begin{array}{l}\begin{array}{c}
\xi_k^T\omega_k-\alpha_k,\;
\xi_k=\half {\Sigma_k^{-1}[[x_*]_k-[y_*]_k]},
\\\alpha_k=\half {\xi_k^T[[x_*]_k+[y_*]_k]}
\end{array}
\;\;\mbox{for Gaussian $\O_k$},\\
\half \omega_k\ln\left([x_*]_k/[y_*]_k\right)-\half [[x_*]_k-[y_*]_k]\;\;\mbox{for scalar Poisson $\O_k$},\\
\half \ln\left([x_*]_{\omega_k}/[y_*]_{\omega_k}\right)\;\;\mbox{for Discrete $\O_k$}.
\end{array}\right.
\ee{phi*_k}
Theorem \ref{the1} in our current situation implies the following statement:
\begin{proposition}\label{propnonstI}
In the framework described in  section \ref{sect:Goal}, assume that the observation scheme $\O^K$ is the direct product of some Gaussian, Discrete and scalar Poisson factors. Let $[x_*;y_*]$ be an optimal solution to the convex optimization problem \rf{mult_opt} associated via the above construction with $\O^K$, and let
$$
\varepsilon_\star = \exp\{\Opt/2\}.
$$
Then the error probabilities of the simple test with detector $\phi^a_*(\omega^K)=\sum_{k=1}^K \phi^k_*(\omega_k)-a$, where $\phi^k_*(\cdot)$ are as in \rf{phi*_k}, and $a\in\bbr$,  satisfy
\[
\epsilon_X\leq\exp\{a\}\varepsilon_\star,\;\;\mbox{and}\;\;\epsilon_Y\leq\exp\{-a\}\varepsilon_\star.\]
Besides this,
no  test can distinguish between these hypotheses with the risk of test less than $\varepsilon_\star^2/4$.
\end{proposition}
\paragraph{Remarks.} Two important remarks are in order.
\par When $\O^K$ is a direct product of Gaussian, Poisson and Discrete factors,  finding the near-optimal simple test reduces to solving explicit well-structured convex optimization problem {\sl with sizes polynomial in $K$ and the maximal dimensions $m_k$ of the factors}, and thus can be done in reasonable time, whenever $K$ and $\max_k m_k$  are ``reasonable.'' This is so in spite of the fact that the ``formal sizes'' of the saddle point problem associated with $\Phi$ could be huge (e.g., when all the factors $\O_k$ are discrete, the cardinality of $\Omega^K$ can grow exponentially with $K$, rapidly making a straightforward  computation of $\Phi$ based on (\ref{Phi}) impossible).
\par We refer to the indexes $k$ and $k'$, $1\leq k,k'\leq K$, as equivalent  in the direct product setup, augmented by convex compact subsets $X,Y$ of $\M^K$, if  $\O_k=\O_{k'}$, $x_k=x_{k'}$ for all $x\in X$, and $y_k=y_{k'}$ for all $y\in Y$. Denoting by $K'$ the number of equivalence classes of indexes, it is clear that problem (\ref{mult_opt}) is equivalent to a problem of completely similar structure, but with $K'$ in the role of $K$.
     It follows that {\sl the complexity of solving (\ref{mult_opt}) is not affected by how large is the number $K$ of factors; what matters is the number $K'$ of equivalence classes of the indexes.} Similar phenomenon takes place when $X$ and $Y$ are direct products  of their projections, $X_k$ and $Y_k$, on the factors $\M_k$ of $\M^K$, and the equivalence of indexes $k$, $k'$ is defined as $\O_k=\O_{k'}$, $X_k=X_{k'}$, $Y_k=Y_{k'}$.

\section{Multiple hypotheses case}\label{sec:multhyp}
The examples outlined in section \ref{sect:appl0} demonstrate that the efficiently computable ``nearly optimal'' simple testing of composite hypotheses suggested by Theorem \ref{the1} and Proposition \ref{propnonstI}, while imposing strong restrictions on the underlying observation scheme, covers nevertheless some interesting and important applications. This testing ``as it is,'' however, deals only with ``dichotomies'' (pairs of hypotheses) of special structure. In this section, we intend to apply our results to the situation when we should decide on more than two hypotheses, or still on two hypotheses, but more complicated than those considered in Theorem \ref{the1}. Our general setup here is   as follows. We are given
a Polish observation space $\Omega$ along with a collection $X_1,...,X_M$ of (nonempty) families of Borel probability distributions on $\Omega$. Given an observation $\omega$ drawn from a distribution $p$ {\sl belonging to the union of these families} (pay attention to this default assumption!), we want to make some conclusions on the ``location'' of $p$. We will be interested in questions of two types:
 \begin{itemize}
 \item[A.] [testing multiple hypotheses] We want to identify the family (or families) in the collection to which $p$ belongs.
 \item[B.] [testing unions] Assume our families $X_1,...,X_M$ are split into two groups -- ``red'' and ``blue''  families. The question is, whether $p$ belongs to a red or a blue family.
 \end{itemize}
 When dealing with these questions, we will assume that for some pairs $(i,j)$, $i\neq j$,  of indexes from $1,...,M$ (let the set of these pairs be denoted $\I$) we are given ``pairwise tests'' $T_{ij}$ deciding on the pairs of hypotheses $H_i$, $H_j$ (where $H_k$ states that $p\in X_k$). To avoid ambiguities, we assume once for ever that the only possible outcomes of a test $T_{ij}$  are either to reject $H_i$ (and accept $H_j$), or to reject $H_j$ (and accept $H_j$).  For $(i,j)\in \I$, we are given the risks $\epsilon_{ij}$ (an upper bound on the probability for $T_{ij}$ to reject $H_i$ when $p\in X_i$) and $\bar{\epsilon}_{ij}$ (an upper bound on the probability for $T_{ij}$ to reject $H_j$ when $p\in X_j$).
 We suppose that whenever $(i,j)\in{\I}$, so is
  $(j,i)$, and the tests $T_{ij}$ and $T_{ji}$ are the same, meaning that when run on an observation $\omega$, $T_{ij}$ accepts $H_i$ if and only if $T_{ji}$ accepts $H_i$. In this case we lose nothing when assuming that $\epsilon_{ij}=\bar{\epsilon}_{ji}$.
  \par
  Our
  goal in this section is to ``assemble'' the pairwise tests $T_{ij}$ into a test for deciding on ``complex'' hypotheses mentioned in A and in B.
  For example, assuming that $T_{ij}$'s are given for all pairs $i,j$ with $i\neq j$, the simplest test for $A$ would be as follows: given observation $\omega$, we run on it  tests $T_{ij}$ for every pair $i,j$ with $i\neq j$, and accept $H_i$ when all tests $T_{ij}$ with $j\neq i$ accept $H_i$. As a result of this procedure, at most one of the hypotheses will be accepted. Applying the union bound, it is immediately seen that if $\omega$ is drawn from $p$ belonging to some $X_i$, $H_i$ will be rejected with probability at most
  $\sum_{j\neq i}\epsilon_{ij}$, so that the quantity $\max_i\sum_{j\neq i}\epsilon_{ij}$   can be considered as the risk of our aggregated test. \par
  The point in what follows is that when $T_{ij}$ are tests of the type yielded by Theorem \ref{the1}, {we have wider ``assembling options''}. Specifically, we will consider the case where
  \begin{itemize}
  \item $T_{ij}$ are ``simple tests induced by detectors $\phi_{ij}$,'' where $\phi_{ij}(\omega):\Omega\to\bbr$ are Borel functions; given $\omega$, $T_{ij}$ accepts $H_i$ when $\phi_{ij}(\omega)>0$, and accepts $H_j$ when $\phi_{ij}(\omega)<0$, with somehow resolved ``ties'' $\phi_{ij}(\omega)=0$. To make $T_{ij}$ and $T_{ji}$ ``the same,'' we will always assume that
 \begin{equation}\label{phiisasymmetric}
      \phi_{ij}(\omega)\equiv -\phi_{ji}(\omega),\;\omega\in\Omega,\;(i,j)\in{\I}.
 \end{equation}
  \item The risk bounds $\epsilon_{ij}$ ``have a specific origin'', namely, they are such that for all $(i,j)\in
  {\I}$,
\be
  \begin{array}{llcll}
  (a)&\int_\Omega\exp\{-\phi_{ij}(\omega)\}p(d\omega)\leq {\epsilon_{ij}}\,\,\forall p\in X_i;&
  (b)&\int_\Omega\exp\{\phi_{ij}(\omega)\}p(d\omega)\leq{\bar{\epsilon}}_{ij},\,\,\forall p\in X_j.\\
  \end{array}
  \ee{16ab}
  \end{itemize}
In the sequel, we refer to the quantities $\widehat{\epsilon}_{ij}:=\sqrt{\epsilon_{ij}\bar{\epsilon}_{ij}}$ as to the {\sl risks} of the detectors $\phi_{ij}$.
  Note that the simple tests provided by Theorem \ref{the1} meet the just outlined assumptions. Another example is the one where $X_i$ are singletons, and the distribution from $X_i$ has density $p_i(\cdot)>0$ with respect to a common for all $i$ measure $P$ on $\Omega$; setting $\phi_{ij}(\cdot)=\half\ln(p_i(\cdot)/p_j(\cdot))$ (so that $T_{ij}$ are the standard likelihood ratio tests) and specifying $\epsilon_{ij}=\bar{\epsilon}_{ij}$ as Hellinger affinities of $p_i$ and $p_j$, we meet our assumptions.
Furthermore,  {\sl every} collection of pairwise tests $\overline{T}_{ij}$, $(i,j)\in\I$, deciding, with risks $\delta_{ij}=\delta_{ji}\in(0,1/2)$, on the hypotheses $H_i$, $H_j$, $(i,j)\in \cI$, gives rise to pairwise detectors $\phi_{ij}$ meeting (\ref{phiisasymmetric}) and (\ref{16ab}) with $\epsilon_{ij}=\bar{\epsilon}_{ij}=2\sqrt{\delta_{ij}(1-\delta_{ij})}$ (cf. remark after Theorem \ref{the1}). Indeed, to this end it suffices to set
$\phi_{ij}(\omega)=\half \ln\left({1-\delta_{ij}\over\delta_{ij}}\right)\overline{T}_{ij}(\omega)$ where, clearly, $\overline{T}_{ij}(\omega)=-\overline{T}_{ji}(\omega)$.
  \par
  The importance of the above assumptions becomes clear from the following immediate observations:
  \begin{enumerate}
  \item By evident reasons, (\ref{16ab}.a) and (\ref{16ab}.b) indeed imply that when $(i,j)\in {\I}$ and $p\in X_i$, the probability for $T_{ij}$ to reject $H_i$ is $\leq\epsilon_{ij}$, while when $p\in X_j$, the probability for the test to reject $H_j$ is $\leq \bar{\epsilon}_{ij}$. Besides this, taking into account that $\phi_{ij}=-\phi_{ji}$, we indeed  ensure $\epsilon_{ij}=
      \bar{\epsilon}_{ji}$;
  \item Relations (\ref{16ab}.a) and (\ref{16ab}.b) are preserved by a shift of the detector -- by passing from $\phi_{ij}(\cdot)$ to $\phi_{ij}(\cdot)-a$ (accompanied with passing from $\phi_{ji}$ to $\phi_{ji}+a$) and simultaneous passing from  $\epsilon_{ij}$, $\bar{\epsilon}_{ij}$ to $\exp\{a\}\epsilon_{ij}$ and $\exp\{-a\}\bar{\epsilon}_{ij}$. In other words, all what matters is the product  $\epsilon_{ij}\bar{\epsilon}_{ij}$ (i.e., the squared risk $\widehat{\epsilon}_{ij}^2$ of the detector $\phi_{ij}$), and we can ``distribute'' this product between the factors as we wish, for example, making ${\epsilon}_{ij}=\bar{\epsilon}_{ij}=\widehat{\epsilon}_{ij}$;
  \item Our assumptions are ``ideally suited'' for passing from a single observation $\omega$ drawn from a distribution $p\in\bigcup\limits_{i=1}^M X_i$ to observing a $K$-tuple
  $\omega^K=(\omega_1,...,\omega_K)$ of observations drawn, independently of each other, from $p$. Indeed, setting $\phi^K_{ij}(\omega_1,...,\omega_K)=\sum_{k=1}^K
\phi_{ij}(\omega_k)$, relations (\ref{16ab}.a) and (\ref{16ab}.b) clearly imply similar relations for  $\phi^K_{ij}$ in the role of $\phi_{ij}$ and $[\epsilon_{ij}]^K$ and $[\bar{\epsilon}_{ij}]^K$ in the role of $\epsilon_{ij}$ and $\bar{\epsilon}_{ij}$. In particular, when $\max(\epsilon_{ij},\bar{\epsilon}_{ij})<1$, passing from a single observation to $K$ of them rapidly decreases the risks as $K$ grows.
\item The left hand sides in relations (\ref{16ab}.a) and (\ref{16ab}.b) are linear in $p$, so that (\ref{16ab}) remains valid when the families of probability distributions $X_i$ are extended to their convex hulls.
\end{enumerate}
 In the rest of this section, we derive ``nontrivial assemblings'' of pairwise tests, meeting the just outlined assumptions, in the context of problems A and B.
\subsection{Testing unions}\label{sec:2unions}
\subsubsection{Single observation case}\label{singleobservationcase}
Let us assume that we are given a family $\P$ of probability measures on a Polish space $\Omega$ equipped with a $\sigma$-additive $\sigma$-finite Borel measure $P$, and all distributions from $\P$ have densities w.r.t. $P$; we identify the distributions from $\P$ with these densities. Let $X_i\subset\P$, $i=1,...,m$ and $Y_j\subset \P$, $j=1,...,n$. Assume that pairwise detectors -- Borel functions $\phi_{ij}(\cdot):\Omega\to\bbr$,  with risk bounded with $\epsilon_{ij}>0$,  are available  for all pairs $(X_i, Y_j)$, $i=1,...,m,\,j=1,...,n$, namely,
\[
\begin{array}{lclclcl}
&\;&
\int_{\Omega} \exp\{-\phi_{ij}(\omega)\}p(\omega)P(d\omega)\leq\epsilon_{ij},\;\forall \,p\in X_{i},&
&\;&\int_{\Omega} \exp\{\phi_{ij}(\omega)\}q(\omega)P(d\omega)\leq\epsilon_{ij},\;\forall q\in Y_{j}.\\
\end{array}
\]
Consider now the problem of deciding between the hypotheses
\[
H_X:\;p\in X=\bigcup\limits_{i=1}^mX_i\;\;\mbox{and} \;\;H_Y:\;p\in Y=\bigcup\limits_{j=1}^n Y_j.\]
on the distribution $p$ of observation $\omega$.
\par Let $E=[\epsilon_{ij}]_{i,j}\in\bbr^{m\times n}$. Consider the matrix $H=\left[\begin{array}{cc}&E\cr E^T&\cr\end{array}\right]$. This is a symmetric entrywise nonzero nonnegative matrix. Invoking the Perron-Frobenius theorem, the leading eigenvalue of this matrix (which is nothing but the spectral norm $\|E\|_{2,2}$ of $E$) is positive, and the corresponding eigenvector can be selected to be nonnegative. Let us denote this vector $z=[g;h]$ with $g\in\bbr^m_{+}$ and $h\in\bbr^n_{+}$, so that
\begin{equation}\label{eigenvector}
Eh=\|E\|_{2,2} g,\,\,E^Tg=\|E\|_{2,2} h.
\end{equation}
We see that if one of the vectors $g$, $h$, is zero, both are so, which is impossible. Thus, both $g$ and $h$ are nonzero nonnegative vectors; since $E$ has all entries positive, (\ref{eigenvector}) says that in fact $g$ and $h$ are positive. Therefore we can
 set
\begin{equation}\label{olddetector}
\begin{array}{rcl}
a_{ij}&=&\ln(h_j/g_i),\,\,1\leq i\leq m,1\leq j\leq n,\\
\phi(\omega)&=&\max\limits_{i=1,...,m}\min\limits_{j=1,...,n}\left[\phi_{ij}(\omega)-a_{ij}\right]:
\Omega\to\bbr.
\end{array}
\end{equation}
Given observation $\omega$, we accept $H_X$ when $\phi(\omega^K)\geq0$, and accept $H_Y$ otherwise.
\begin{proposition}\label{propnonstation1pairnewnot} In the described situation, we have
\begin{equation}\label{standardreason1}
\begin{array}{ll}
(a)&\int_\Omega \exp\{-\phi(\omega)\}p(\omega)P(d\omega)\leq\varepsilon:=\|E\|_{2,2},\,\,p\in X,\\
(b)&\int_\Omega \exp\{\phi(\omega)\}p(\omega)P(d\omega)\leq\varepsilon,\,\,p\in Y.\\
\end{array}
\end{equation}
As a result, the risk of the just described test when testing $H_X$ versus $H_Y$ does not exceed $\varepsilon=\|E\|_{2,2}$.
\end{proposition}
\subsubsection{Case of repeated observations} The above construction and result admit immediate extension onto the case of non-stationary repeated observations. Specifically,
consider the following situation. For $1\leq t\leq K$, we are given
\begin{enumerate}
\item Polish space  $\Omega_t$ equipped with Borel $\sigma$-additive $\sigma$-finite measure $P_t$,
\item A family $\cP_t$ of Borel probability densities, taken w.r.t. $P_t$, on $\Omega_t$,
\item Nonempty sets $X_{it}\subset \cP_t$, $Y_{jt}\subset\cP_t$, $i\in\I_t=\{1,...,m_t\}$, $j\in\J_t=\{1,...,n_t\}$,
\item {\sl Detectors} -- Borel functions $\phi_{ijt}(\cdot):\Omega_t\to\bbr$, $i\in\I_t$, $j\in\J_t$, along with positive reals $\epsilon_{ijt}$, $i\in\I_t$, $j\in \J_t$, such that
\begin{equation}\label{suchthatnew177}
\begin{array}{ll}
(a)&\int_{\Omega_t} \exp\{-\phi_{ijt}(\omega)\}p(\omega)P_t(d\omega)\leq\epsilon_{ijt}\,\,\forall (i\in\I_t,j\in\J_t,p\in X_{it}),\\
 (b)&\int_{\Omega_t} \exp\{\phi_{ijt}(\omega)\}p(\omega)P_t(d\omega)\leq\epsilon_{ijt}\,\,\forall (i\in\I_t,j\in\J_t,p\in Y_{jt}),\\
 \end{array}
 \end{equation}
\end{enumerate}
{Given time horizon $K$, consider two hypotheses on observations $\omega^K=(\omega_1,...,\omega_K)$, $\omega_t\in\Omega_t$, $H_1:=H_X$ and $H_2:=H_Y$, as follows. According to hypothesis $H_\chi$, $\chi=1,2$, the observations $\omega_t$, $t=1,2,...,K$, are generated as follows:
 \begin{quote}
``In the nature'' there exists a sequence of  ``latent'' random variables $\zeta_{1,\chi},\zeta_{2,\chi},\zeta_{3,\chi},...$ such that $\omega_t$, $t\leq K$,
  is a deterministic function of $\zeta^t_\chi=(\zeta_{1,\chi},...,\zeta_{t,\chi})$, and the conditional, $\zeta^{t-1}_\chi$ being fixed, distribution of $\omega_t$ has density $p_t\in\cP_t$ w.r.t. $P_t$, the density $p_t$ being a deterministic function of $\zeta^{t-1}_\chi$.
  Moreover, {\sl when $\chi=1$,  $p_t$ belongs to $X_t:=\bigcup\limits_{i\in\I_t} X_{ti}$, and when $\chi=2$, it belongs to $Y_t:=\bigcup_{j\in\J_t} Y_{jt}$.  } \end{quote}
Our goal is to decide from observations $\omega^K=(\omega_1,...,\omega_K)$ on the hypotheses $H_X$ and $H_Y$.}
\paragraph{The test} we intend to consider is as follows. We set
\begin{equation}\label{Et}
E_t=[\epsilon_{ijt}]_{i,j}\in\bbr^{m_t\times n_t},\,\,H_t=\left[\begin{array}{cc}&E_t\cr
E_t^T&\cr\end{array}\right]\in\bbr^{(m_t+n_t)\times(m_t+n_t)},\,\,\varepsilon_t=\|E_t\|_{2,2}.
\end{equation}
As above, the leading eigenvalue of the symmetric matrix $H_t$ is $\varepsilon_t$, the corresponding eigenvector $[g^t;h^t]$, $g^t\in\bbr^{m_t},h^t\in\bbr^{n_t}$ can be selected to be positive, and we have
\begin{equation}\label{Perron}
E_th^t=\varepsilon_t g^t,\,\,E_t^Tg^t=\varepsilon_t h^t.
\end{equation}
We set
\begin{equation}\label{aijs}
\begin{array}{rcl}
a_{ijt}&=&\ln(h^t_j/g^t_i),\,\,1\leq i\leq m_t,1\leq j\leq n_t,\\
\phi_t(\omega_t)&=&\max\limits_{i=1,...,m_t}\min\limits_{j=1,...,n_t}\left[\phi_{ijt}(\omega_t)-a_{ijt}\right]:
\Omega\to\bbr,\\
\phi^K(\omega^K)&=&\sum_{t=1}^K\phi_t(\omega_t).\\
\end{array}
\end{equation}
Given observation $\omega^K=(\omega_1,...,\omega_K)$, we accept $H_X$ when $\phi^K(\omega^K)\geq0$, and accept $H_Y$ otherwise.
\par
We have the following analogue of Proposition \ref{propnonstI}
\begin{proposition}\label{propnonstationarynew} In the situation of this section, we have
\begin{equation}\label{standardreason2}
\begin{array}{ll}
(a)&\int_\Omega \exp\{-\phi_t(\omega)\}p(\omega)P(d\omega)\leq\varepsilon_t:=\|E_t\|_{2,2},\,\,p\in X_t,t=1,2,...\\
(b)&\int_\Omega \exp\{\phi_t(\omega)\}p(\omega)P(d\omega)\leq\varepsilon_t,\,\,p\in Y_t,t=1,2,...\\
\end{array}
\end{equation}
As a result, the risk of the just described test does not exceed $\prod_{t=1}^K\varepsilon_t$.
\end{proposition}
Some remarks are in order.
\paragraph{Symmeterizing the construction.} Inspecting the proof of Proposition \ref{propnonstationarynew},  we see that the validity of its risk-related conclusion is readily given by the validity
of (\ref{standardreason2}). The latter relation, in turn,  is ensured by the described in  (\ref{aijs}) scheme of ``assembling'' the detectors $\phi_{ijt}(\cdot)$ into $\phi_t(\cdot)$, but this is not the only assembling  ensuring  (\ref{standardreason2}). For example, swapping $X_t$ and $Y_t$, applying the assembling (\ref{aijs}) to these ``swapped'' data and ``translating'' the result back to the original data, we arrive at the detectors
$$
\overline{\phi}_t(\omega)=\min_{j=1,...,n_t}\max_{i=1,...,m_t} [\phi_{ijt}(\omega)-a_{ijt}],
$$
with $a_{ijt}$ given by (\ref{aijs}), and these new detectors, when used in the role of $\phi_t$,  still ensure (\ref{standardreason2}). Denoting by $\underline{\phi}_t$ the detector $\phi_t$ given by (\ref{aijs}), observe that $\underline{\phi}_t(\cdot)\leq\overline{\phi}_t(\cdot)$, and this inequality in general is strict. Inspecting the proof of Proposition \ref{propnonstationarynew}, it is immediately seen that Proposition remains true whenever $\phi^K(\omega^K)=\sum_{t=1}^K\phi_t(\omega_t)$ with $\phi_t(\cdot)$ satisfying the relations
$$
\underline{\phi}_t(\cdot)\leq\phi_t(\cdot)\leq \overline{\phi}_t(\cdot),
$$
for example, with the intrinsically symmetric ``saddle point'' detectors
$$
\phi_t(\cdot)=\max_{\lambda\in\Delta_{m_t}}\min_{\mu\in\Delta_{n_t}}\sum_{i,j}\lambda_i\mu_j[\phi_{ijt}(\cdot)-a_{ijt}]\eqno{[\Delta_k=\{x\in\bbr^k:x\geq0,\sum_{i=1}^kx_i=1\}]}
$$
Needless to say, similar remarks hold true in the context of Proposition \ref{propnonstation1pairnewnot}, which is nothing but the stationary (i.e., with $K=1$) case of Proposition \ref{propnonstationarynew}.
\paragraph{Testing convex hulls.} As it was already mentioned, the risk-related conclusions in Propositions \ref{propnonstation1pairnewnot}, \ref{propnonstationarynew} depend solely on the validity of relations (\ref{standardreason1}), (\ref{standardreason2}). Now,  density $p(\cdot)$ enters the left hand sides in (\ref{standardreason1}), (\ref{standardreason2}) linearly, implying that when, say, (\ref{standardreason2}) holds true for some $X_t$, $Y_t$, the same relation holds true when the families of probability densities $X_t$, $Y_t$ are extended to their convex hulls. Thus, in the context of Propositions \ref{propnonstation1pairnewnot}, \ref{propnonstationarynew} we, instead of speaking about testing {\sl unions}, could speak about testing {\sl convex hulls} of these unions.
\paragraph{Simple illustration.} Let  $p$ be a positive probability density on the real axis $\Omega=\bbr$ such that setting $\rho_i=\int \sqrt{p(\omega)p(\omega-i)}d\omega$, we have $\varepsilon:=2\sum_{i=1}^\infty\rho_i <\infty$. Let $p_i(\omega)=p(\omega-i)$, and let $I=\{\imath_1<...<\imath_m\}$ and $J=\{\jmath_1<...<\jmath_n\}$ be two non-overlapping finite subsets of $\bbz$.  Consider the case where $X_{it}=\{p_{\imath_i}(\cdot)\}$, $1\leq i\leq m=m_t$, $Y_{jt}=\{p_{\jmath_j}(\cdot)\}$, $1\leq j\leq n=n_t$, are singletons, and let us set
$$
\begin{array}{rcl}
\phi_{ijt}(\omega)&=&{1\over 2}\ln(p_{\imath_i}(\omega)/p_{\jmath_j}(\omega)),\,1\leq i\leq m,\,1\leq j\leq n,\\
\epsilon_{ijt}&=&\int\sqrt{p_{\imath_i}(\omega)/p_{\jmath_j}(\omega)}d\omega,\,1\leq i\leq m,\,1\leq j\leq n.\\
\end{array}
$$
This choice clearly ensures (\ref{suchthatnew177}), and for the associated matrix $E_t\equiv E$ we have $\|E\|_{2,2}\leq \varepsilon$.\footnote{We use the following elementary fact: {\sl Let $E$ be a matrix with sums of magnitudes of entries in every row and every column not exceeding $r$. Then $\|E\|_{2,2}\leq r$.} To be on the safe side, here is the proof: let $F=\left[\begin{array}{cc}&E\cr E^T&\cr\end{array}\right]$, so that $\|E\|_{2,2}=\|F\|_{2,2}$, and $\|F\|_{2,2}$ is just the spectral radius of $F$. We clearly have $\|Fx\|_\infty\leq r\|x\|_\infty$ for all $x$, whence the spectral radius of $F$ is at most $r$.} Thus, when $\varepsilon$ is small, we can decide with low risk on the hypotheses associated with $X_t:=\bigcup\limits_{i=1}^mX_{it}$, $Y_t:=\bigcup\limits_{j=1}^nY_{jt}$; note that $\varepsilon$ is independent of the magnitudes of $m,n$. Moreover, when $\varepsilon<1$, and repeated observations, of the structure considered in Proposition \ref{propnonstationarynew}, are allowed, $K=\rfloor\ln(1/\epsilon)/\ln(1/\varepsilon)\lfloor$ observations are sufficient to get a test with risk $\leq\epsilon$, and $K$ again is not affected by the magnitudes of $m,n$. Finally, invoking the above remark, we can replace in these conclusions the  finite  sets of probability densities $X_t$, $Y_t$ with their convex hulls.
\subsection{Testing multiple hypotheses}\label{sectionmultiple}
Let $X_1,...,X_m$ be nonempty sets in the space of Borel probability distributions on a Polish space $\Omega$, $E=[\epsilon_{ij}]$ be a symmetric $m\times m$ matrix with zero diagonal and positive off-diagonal entries, and let
$$
\phi_{ij}(\omega)=-\phi_{ji}(\omega):\Omega\to\bbr,\,1\leq i,j\leq m,\,i\neq j,
$$
be Borel {\sl detectors} such that
\begin{equation}\label{detectorsasusual}
\forall (i,j,1\leq i,j\leq m,i\neq j):\int_\Omega\exp\{-\phi_{ij}(\omega)\}p(d\omega)\leq\epsilon_{ij}\,\,\forall p\in X_i.
\end{equation}
Given a skew-symmetric matrix $[\alpha_{ij}]_{1\leq i,j\leq m}$ and setting $\bar{\phi}_{ij}(\cdot)=\phi_{ij}(\cdot)-\alpha_{ij}$, we get
\begin{equation}\label{detectorsasusual1}
\forall (i,j,1\leq i,j\leq m,i\neq j):\int_\Omega\exp\{-\bar{\phi}_{ij}(\omega)\}p(d\omega)\leq\exp\{\alpha_{ij}\}\epsilon_{ij}\,\,\forall p\in X_i.
\end{equation}
Consider the following test aimed to decide, given an observation $\omega$ drawn from a distribution $p$ known to belong to $X=\bigcup\limits_{i=1}^m X_i$, on $i$ such that $p\in X_i$ (we refer to the validity of the latter inclusion as to hypothesis $H_i$). The test is as follows: we compute $\bar{\phi}_{ij}(\omega)$ for all $i\neq j$, and accept all $H_i$'s such that all the quantities $\bar{\phi}_{ij}(\omega)$ with $j$ distinct from $i$ are positive. Note that since $\bar{\phi}_{ij}(\cdot)\equiv -\bar{\phi}_{ji}(\cdot)$, if some $H_i$ is accepted by our test, no $H_{i'}$ with $i'$ different from $i$ can be accepted; thus, our test, for every $\omega$, accepts at most one of the hypotheses $H_i$.  Let us denote by $\epsilon_i$ the maximal, over $p\in X_i$, probability for the test to reject $H_i$ when our observation $\omega$ is drawn from $p(\cdot)$. Note that since our test accepts at most one of $H_i$'s, for every $i$ the probability to accept $H_i$ when the observation $\omega$ is drawn from a distribution $p(\cdot)\in X\backslash X_i$ (i.e., when $H_i$ is false) does not exceed $\max_{j:j\neq i}\epsilon_j$.
\par
Now recall that the risks $\epsilon_i$ depend on the shifts $\alpha_{ij}$, and consider the problem as follows. Given ``importance weights'' $p_i>0$, $1\leq i\leq m$, we now aim to find the shifts $\alpha_{ij}$ resulting  in the smallest possible quantity
$$
\epsilon:=\max_{1\leq i\leq m} p_i\epsilon_i,
$$
or, more precisely, the smallest possible natural upper bound $\varepsilon$ on this quantity. We define this bound as follows.
\par
Let, for some $i$, an observation $\omega$ be drawn from a distribution $p\in X_i$. Given this observation, $H_i$ will be rejected if for some $j\neq i$ the quantity $\bar{\phi}_{ij}(\omega)$ is nonpositive. By (\ref{detectorsasusual}), for a given $j\neq i$, $p$-probability of the event in question is at most $\exp\{\alpha_{ij}\}\epsilon_{ij}$, which implies the upper bound on $\epsilon_i$, specifically, the bound
$$
\varepsilon_i=\sum_{j\neq i}\exp\{\alpha_{ij}\}\epsilon_{ij}=\sum_{j=1}^m\exp\{\alpha_{ij}\}\epsilon_{ij}
$$
(recall that $\epsilon_{ii}=0$ for all $i$). Thus, we arrive at the upper bound
\begin{equation}\label{varepsilonmultiple}
\varepsilon:=\max_ip_i\varepsilon_i=\max_i\sum_{j=1}^mp_i\epsilon_{ij}\exp\{\alpha_{ij}\}
\end{equation}
on $\epsilon$.
What we want is to select $\alpha_{ij}=-\alpha_{ji}$ minimizing this bound.
\par
Our goal is relatively easy to achieve: all we need is to solve the {\sl convex} optimization problem
\begin{equation}\label{varepsilonproblem}
\varepsilon_*=\min_{\alpha=[\alpha_{ij}]}\left\{f(\alpha):=\max_{1\leq i\leq m} \sum_jp_i\epsilon_{ij}\exp\{\alpha_{ij}\}:\alpha=-\alpha^T\right\}.
\end{equation}
The problem \rf{varepsilonproblem} allows for a ``closed form'' solution.
\begin{proposition}\label{prop_alpha1}
Let $\rho$ be the Perron-Frobenius eigenvalue of the entry-wise nonnegative matrix $\bar{E}=\left[p_i\epsilon_{ij}\right]_{1\leq i,j\leq m}$. The corresponding eigenvector $g\in \bbr^m$ can be selected to be {\em positive}, and for the choice
$
[\bar{\alpha}_{ij}:=\ln(g_j)-\ln(g_i)]_{i,j}, \;\;1\leq i,j\leq m,
$
one has
$
\varepsilon_*=f(\bar{\alpha})=\rho.
$
\end{proposition}

\paragraph{Remark.} The proof of Proposition \ref{prop_alpha1} demonstrates  that with the optimal assembling given by  $\alpha_{ij}=\bar{\alpha}_{ij}$ all the quantities $p_i\epsilon_i$ in (\ref{varepsilonmultiple}) become equal to $\varepsilon_*=\rho$. In particular, when $p_i=1$ for all $i$, for every $i$ the probabilities to reject $H_i$ when the hypothesis is true, and to accept $H_i$ when the hypothesis is false, are upper bounded by $\rho$.
\subsubsection{A modification}\label{multiple}
In this section we focus on multiple hypothesis testing in the case when all importance factors $p_i$ are equal to 1. Note that in this case the result we have just established can be void when the optimal value $\varepsilon_*$ in (\ref{varepsilonproblem}) is $\geq1$, as this is the case, e.g., when  some $X_i$ and $X_j$ with $i\neq j$ intersect. In the latter case, for every pair $i,j$ with $i\neq j$ and $X_i\cap X_j\neq\emptyset$, the best -- resulting in the smallest possible value of $\epsilon_{ij}$ -- selection of $\phi_{ij}$ is $\phi_{ij}\equiv 0$, resulting in $\epsilon_{ij}=1$. It follows that even with $K$-repeated observations (for which $\epsilon_{ij}$ should be replaced with $\epsilon_{ij}^K$) the optimal value in (\ref{varepsilonproblem}) is $\geq1$, so that our aggregated test allows for only trivial bound $\varepsilon\leq1$ on $\varepsilon$, see (\ref{varepsilonmultiple}).\footnote{Of course, the case in question is intrinsically difficult -- here no test whatsoever can make all the risks $\epsilon_i$ less than $1/2$.} Coming back to the general situation where $p_i\equiv 1$ and $\varepsilon_*$ is large, what can we do? A solution, applicable when $\epsilon_{ij}<1$ for all $i,j$, is to pass to $K$-repeated observations; as we have already mentioned, this is equivalent to passing from the original matrix $E=[\epsilon_{ij}]$ to its entrywise power $E^{(K)}=[\epsilon_{ij}^K]$; when $K$ is large,
the leading eigenvalue $\rho_K$ of $E^{(K)}$ becomes small. The question is what to do if some of $\epsilon_{ij}$ indeed are equal to 1, and a somewhat partial solution in this case may be obtained by substituting our original goal of highly reliable recovery of the true hypothesis with a less ambitious one. A natural course of action could be as follows. Let {$\I$} be the set of all ordered
 pairs $(i,j)$ with $1\leq i,j\leq m$, and let $\C $ be a given subset of this set containing all ``diagonal'' pairs $(i,i)$. We interpret the inclusion $(i,j)\in \C $ as the claim that $H_j$ is ``close'' to $H_i$. \footnote{Here the set of ordered pairs $\C $ is not assumed to be invariant w.r.t. swapping the components of a pair, so that in general ``$H_j$ is close to $H_i$'' is not the same as ``$H_i$ is close to $H_j$.''} Imagine that what we care about when deciding on the collection of hypotheses $H_1,...,H_m$ is not to miss a correct hypothesis and, at the same time, to reject all hypotheses which are ``far'' from the true one(s). This can be done by test as follows. Let us shift somehow the original detectors, that is, pass from $\phi_{ij}(\cdot)$ to the detectors $\phi^\prime_{ij}(\cdot)=\phi_{ij}(\cdot)-\alpha_{ij}$ with $\alpha_{ij}=-\alpha_{ij}$, thus ensuring that
\begin{equation}\label{thusthat}
 \phi^\prime_{ij}(\cdot):=-\phi^\prime_{ji}(\cdot)\ \&\ \int_\Omega \exp\{-\phi^\prime_{ij}(\omega)\}p(d\omega)\leq \epsilon^\prime_{ij}:=\exp\{\alpha_{ij}\}\epsilon_{ij} \,\,\forall p\in X_i.
 \end{equation}
Consider the test as follows:
\begin{quote}
{\bf Test $\T$:} Given observation $\omega$, we compute the matrix $[\phi^\prime_{ij}(\omega)]_{ij}$. Looking one by one at the rows $i=1,2,...m$ of this matrix, we accept $H_i$ if all the entries $\phi^\prime_{ij}(\omega)$ with $(i,j)\not\in \C $ are positive, otherwise we reject $H_i$.
\par
The outcome of the test is the collection of all accepted hypotheses (which now is not necessary either empty or  a singleton).
\end{quote}
What we can say about this test is the following. Let
\begin{equation}\label{setepsilon}
\epsilon=\max_i\sum_{j: (i,j)\not\in \C } \epsilon^\prime_{ij},
\end{equation}
and let the observation $\omega$ the test is applied to be drawn from distribution $p\in X_{i_*}$, for some $i_*$. Then
\begin{itemize}
\item if, for some $i\neq j$, $\T$ accepts both $H_i$ and $H_j$, then either $H_j$ is close to $H_i$, or $H_i$ is close to $H_j$, or both.\\
Indeed, if neither $H_i$ is close to  $H_j$, nor $H_j$ is close to $H_i$, both $H_i$, $H_j$ can be accepted only when $\phi^\prime_{ij}(\omega)>0$ and $\phi^\prime_{ji}(\omega)>0$, which is impossible due to $\phi^\prime_{ij}(\cdot)=-\phi^\prime_{ji}(\cdot)$.
\item $p$-probability for the true hypothesis $H_{i_*}$ not to be accepted is at most $\epsilon$.\\
Indeed, by (\ref{thusthat}), the $p$-probability for $\phi^\prime_{i_*j}$ to be nonpositive does not exceed $\epsilon^\prime_{i_*j}$. With this in mind, taking into account the description of our test and applying the union bound, $p$-probability
to reject $H_{i_*}$ does not exceed $\sum_{j: (i_*,j)\not\in \C }\epsilon^\prime_{i_*j}\leq\epsilon$.
\item $p$-probability of the event $\E$ which reads {\sl ``at least one of the accepted $H_i$'s is such that both $(i,i_*)\not\in \C $ and $(i_*,i)\not\in \C $''}  (that is, neither $i_*$ is close to $i$, nor $i$ is close to $i_*$)  does not exceed $\epsilon$.\\
Indeed, let $I$  be the set of all those $i$ for which $(i,i_*)\not\in \C$ and $(i_*,i)\not\in\C $. For a given $i\in I$, $H_i$ can be accepted by our test only when $\phi^\prime_{ii_*}(\omega)>0$ (since $(i,i_*)\not\in \C $), implying that $\phi^\prime_{i_*i}(\omega)<0$. By (\ref{thusthat}), the latter can happen with $p$-probability at most $\epsilon^\prime_{i_*i}$. Applying the union bound, the $p$-probability of the event $\E$ is at most
\[
\sum_{i\in I}\epsilon^\prime_{i_*i}\leq \sum_{i:(i_*,i)\not\in\C }\epsilon^\prime_{i_*i}\leq\epsilon
 \]
(we have taken into account that whenever $i\in I$, we have $(i_*,i)\not\in \C $, that is, $I\subset \{i: (i_*,i)\not\in\C\} $).
\end{itemize}
When $\epsilon$ is small (which, depending on how closeness is specified, can happen even when some of $\epsilon^\prime_{ij}$ are not small), the simple result we have just established is ``better than nothing:'' it says that up to an event of probability $2\epsilon$, the true hypotheses $H_{i_*}$ is accepted, and all accepted hypotheses $H_j$ are such that either $j$ is close to $i_*$, or $i_*$ is close to $j$, or both.
\par
Clearly, given $\C$, we would like to select $\alpha_{ij}$ to make $\epsilon$ as small as possible. The punch line is that this task is relatively easy: all we need is to solve the {\sl convex} optimization problem
\begin{equation}\label{selectingalphas}
\min_{[\alpha_{ij}]_{i,j}}\left\{\max_{1\leq i\leq m}\sum_{j: (i,j)\not\in \C}\epsilon_{ij}\exp\{\alpha_{ij}\}:\alpha_{ij}\equiv -\alpha_{ji}\right\}.
\end{equation}
\paragraph{Special case: testing multiple unions.} Consider the case when ``closeness of hypotheses'' is defined as follows: the set $\{1,...,M\}$ of hypotheses' indexes is split into $L\geq2$ nonempty non-overlapping subsets $\I_1,...,\I_L$, and $H_j$ is close to $H_i$ if and only if both $i,j$ belong to the same element of this partition. Setting $E=[\epsilon_{ij}]_{i,j}$, let $D=[\delta_{ij}]$ be the matrix obtained from $E$ by zeroing out all entries $ij$ with $i,j$ belonging  to $\I_\ell$ for some $1\leq \ell\leq L$.  Problem (\ref{selectingalphas}) now reads
$$
\min\limits_{[\alpha_{ij}]}\left\{\max\limits_{1\leq i\leq M}\sum\limits_{1\leq j\leq M}\delta_{ij}\exp\{\alpha_{ij}\}:\alpha=-\alpha^T\right\}.
$$
This problem, similarly to problem (\ref{varepsilonproblem}), admits a closed form solution: the Perron-Frobenius eigenvector $g$ of the entrywise nonnegative symmetric matrix $D$  can be selected to be positive, an optimal solution is given by
$
\alpha_{ij} =\ln(g_j)-\ln(g_i),
$
and  the optimal value is $\epsilon_*:=\|D\|_{2,2}$. Test $\T$ associated with the optimal solution can be converted into a test $\widehat{\T}$ deciding on $L$ hypotheses $\H_\ell=\bigcup\limits_{i\in \I_\ell}H_i$, $1\leq \ell\leq k$; specifically, when $\T$ accepts some hypothesis $H_i$, $\widehat{\T}$ accepts hypothesis $\H_\ell$ with $\ell$ uniquely defined by the requirement $i\in\I_\ell$. The above results on $\T$ translate in the following facts about $\widehat{\T}$:
\begin{itemize}
\item  $\widehat{\T}$ never accepts more than one hypothesis;
\item let the observation $\omega$ on which $\widehat{\T}$ is run be drawn from a distribution $p$ obeying, for some $1\leq i
\leq M$, the hypothesis $H_i$, and let $\ell$ be such that $i\in\I_\ell$. Then the $p$-probability for $\widehat{\T}$ to reject the hypothesis $\H_\ell$ is at most $\epsilon_*$.
\end{itemize}
When  $L=2$ we come back to the situation considered in section \ref{singleobservationcase}, and what has just been  said about $\widehat{\T}$ recovers the risk-related result of Proposition \ref{propnonstation1pairnewnot}; moreover, when $L=2$, the test  $\widehat{\T}$ is, essentially, the test based on the detector $\phi$ given by (\ref{olddetector}).\footnote{The only subtle difference, completely unimportant in our context,  is that the latter test accepts $\H_1$ whenever $\phi(\omega)\geq0$ and accepts $\H_2$ otherwise, while $\widehat{\T}$ accepts
 $\H_1$ when $\phi(\omega)>0$, accepts $\H_2$ when $\phi(\omega)<0$ and accepts nothing when $\phi(\omega)=0$.}
 Note that when $L>2$, one could use the detector-based tests, yielded by the construction in
 section \ref{singleobservationcase}, to build ``good'' detectors for the pairs of hypotheses $\H_\ell$, $\H_{\ell'}$ and then assemble these detectors, as explained in section  \ref{sectionmultiple}, into a test deciding on multiple hypotheses $\H_1,...,\H_L$, thus getting an ``alternative'' to $\widehat{\T}$ test $\widetilde{\T}$. Though both tests are obtained by aggregating detectors $\phi_{ij}$, $1\leq i,j\leq M$, in the test $\widehat{\T}$ we aggregate them ``directly'', while the aggregation in test $\widetilde{\T}$ is done in two stages where we first assemble $\phi_{ij}$ into pairwise detectors $\widetilde{\phi}_{\ell\ell'}$ for $\H_\ell$, $\H_{\ell'}$, and then assemble these new detectors into a test for multiple hypotheses $\H_1,...,\H_L$. However, the performance guarantees for the test $\widetilde{\T}$ can be only worse than those for the test $\widehat{\T}$ --
 informally, when assembling $\phi_{ij}$ into $\widetilde{\phi}_{\ell,\ell'}$, we take into account solely the ``atomic contents'' of the aggregated hypotheses $\H_\ell$ and $\H_{\ell'}$, that is, look only at the ``atoms'' $H_i$ with  $i\in\I_\ell\cup \I_{\ell'}$, while when assembling $\phi_{ij}$ into $\widehat{\T}$, we look at all $m$ atoms simultaneously.\footnote{The formal reasoning is as follows. On a close inspection, to get risk bound $\widetilde{\epsilon}$ for $\widetilde{\T}$,  we start with the $M\times M$ matrix $D$ partitioned into $L\times L$ blocks $D^{\ell\ell'}$ (this partitioning is induced by splitting the indexes of rows and columns into the groups $\I_1$,...,$\I_L$),  and form the $L\times L$ matrix $G$ with entries $\gamma_{\ell\ell'}=\|D^{\ell\ell'}\|_{2,2}$; $\widetilde{\epsilon}$ is nothing but $\|G\|_{2,2}$, while the risk bound $\epsilon_*$ for $\widehat{\T}$ is $\|D\|_{2,2}$. Thus, $\epsilon_*\leq\widetilde{\epsilon}$ by the construction of matrix $G$ from $D$.}
 \paragraph{Near-optimality.} Let the observation scheme underlying the just considered ``multiple unions'' situation be $K$-repeated version  $\O^K$ of a good observation scheme $\O=((\Omega,P),\{p_\mu(\cdot):\mu\in\M\},\F)$, meaning that our observation is $\omega=\omega^K:=(\omega_1,...,\omega_K)$ with $
 \omega_t$ drawn, independently of each other, from a distribution $p$, and $i$-th of our $M$ hypotheses, $H_i$, states that $p$ belongs to the set $X_i=\{p_\mu:\mu\in Q_i\}$, where $Q_i$ are convex compact subsets of $\M$.
Let $\phi_{ij}$ be the pairwise detectors for $H_i$ and $H_j$ yielded by Theorem \ref{the1}, and let  $\widehat{\T}^K$ be the test deciding on aggregated hypotheses $\H_\ell$'s from $K$-repeated observations $\omega^K$ and built by assembling detectors  $\phi_{ij}^K=\sum_{t=1}^K\phi_{ij}(\omega_t)$.
   We have the following near-optimality result (cf. Proposition \ref{newcol1}):
  \begin{proposition}\label{theverylatestopt} In the just described situation and given $\epsilon\in(0,1/4)$, assume that in the nature there exists a test $\overline{T}$, based on $\bar{K}$-repeated observations $\omega^{\bar{K}}$, deciding on $\H_1,...,\H_L$ and such that $\overline{T}$ never accepts more than one hypothesis and, for every $\ell\leq L$,  rejects $\H_\ell$ when the hypothesis is true with probability $\leq\epsilon$. Then the same performance guarantees are shared by the test $\widehat{\T}^K$, provided that
 \[K\geq{2\ln(M/\epsilon)\over\ln(1/\epsilon)-2\ln 2}\bar{K}.
 \]
 \end{proposition}
 \hide{}{
\subsection{Sequential Hypothesis Testing}
\subsubsection{Situation}
Let $\O=((\Omega,P),\{p_\mu(\cdot):\mu\in\M\},\F)$ be a good observation scheme; for the sake of simplicity (and basically without loss of generality) assume that the mapping $\mu\to p_\mu(\cdot)$ is an injection on $\M$, that is different values of the parameter $\mu$ correspond to different probability densities $p_\mu(\cdot)$.
   Let, further,  $X_j$, $1\leq j\leq J$, be nonempty convex compact subsets of $\M$; these sets define the sets of probability densities  $P_j=\{p_\mu(\cdot):\mu\in X_j\}$. \par
    Let
   $$
   \J:=\{1,2,...,J\}=\bigcup\limits_{i=1}^I\J_i
   $$
   be a partition of the set $\J$ of indexes of $X_j$'s into $I\geq 2$ non-overlapping nonempty groups $\J_1,...,\J_I$. We associate with the sets $X_j$ hypotheses $H_j$ on the distribution of stationary $K$-repeated observation $\omega^K=(\omega_1,...,\omega_K)$; according to $H_j$, $\omega_1,...,\omega_K$ are drawn, independently of each other, from a distribution $p_\mu\in P_j$. Our goal is to decide from observation $\omega^K$ on the hypotheses $\H_1,...,\H_I$, $\H_i$ stating that the probability density underlying observations belongs to $P^i=\bigcup\limits_{j\in\J_i}P_j$.
   \par
   From now on, we make the following assumption:
   \begin{quote}
   {\bf A:} {\sl When $j,j'\in\J$ do not belong to a common group $\J_i$, the sets of probability densities $P_j$, $P_{j'}$ are at positive Hellinger distance from each other.}
   \end{quote}
   \par
                     From the results of section \ref{multiple} it easily follows that under assumption {\bf A} and given  $\epsilon>0$, we can decide on the hypotheses $\H_1,...,\H_I$ with risk $
                     \leq\epsilon$, provided that the number $K$ of observations is large enough. This being said, the ``large enough'' $K$ could be indeed quite large, provided that for some pairs $j\in \J_i,j'\in\J_{i'}$ with $i\neq i'$, the Hellinger distances between $P_j$ and $P_{j'}$ are small. To overcome, to some extent, this difficulty, we can switch from decision rules based on $K$ observations to {\sl sequential} decision rules, where the decision is made on the basis of on-line adjustable number of observations, in hope that when we are lucky and the distribution $p_*$ underlying our observation is ``deeply inside'' of some $P^{i_*}$ and thus is ``far'' from all $P^i$, $i\neq i_*$, the true hypothesis $\H_{i_*}$ will be accepted much sooner than in the case when $p_*$ is close to some of ``wrong'' $P^i$'s. We are about to utilize our previous results to build sequential tests. What follows can be considered as a ``computationally friendly'' variation on the classical topic of Sequential Analysis originating from \cite{Wald1947}, see also \cite{Chernoff1972,Bakeman1997} and references therein.
\subsubsection{The construction}
\paragraph{Sequential test: setup.} The setup for our ``generic'' sequential test is given by
\begin{enumerate}
\item required {\sl risk} $\epsilon\in(0,1)$;
\item positive integer $S$ -- {\sl number of stages}, along with the following entities, defined for $1\leq s\leq S$ and forming {\sl $s$-th component} of the setup:
\begin{enumerate}
\item positive reals $\epsilon_s$, $1\leq s\leq S$, such that $2\sum_{s=1}^S\epsilon_s=\epsilon$;
\item representations $X_j=\bigcup\limits_{\iota=1}^{\iota_{js}} X_{j\iota s}$, $j\in \J$, where $X_{j\iota s}$ are nonempty closed convex subsets of $X_j$;\\
\item {\sl tolerances} $\delta_s\in(0,1)$.
\end{enumerate}
\end{enumerate}
\paragraph{Colors.} It is convenient to think about the (values of the) indexes $i=1,...,I$ of the groups $\J_i$ as of the {\sl colors} of these groups. Let us agree that these colors are inherited by different entities ``stemming'' from $\J_i$'s, like
 \begin{itemize}
\item indexes $j\in \J_i$,
\item the sets $X_j$, $X_{j\iota s}$, $P_j=\{p_\mu(\cdot):\mu\in X_j\}$, $P_{j\iota s}=\{p_\mu(\cdot):\mu\in X_{j\iota s}\}$ with $j\in \J_i$,
\item  the points $\mu\in X_j$ and the densities $p(\cdot)\in P_j$, $j\in \J_i$,
\item the hypotheses $H_j$ stating that the probability density underlying observations belongs to $P_j$, $j\in \J_i$,
\end{itemize}
etc. Note that under assumption {\bf A} assigning the above entities with colors is unambiguous, and that in terms of colors, our goal is to identify the color of the distribution underlying the observations. \par
\paragraph{Detectors.} To avoid messy notation, we enumerate, for every $s\leq S$, the sets $X_{j\iota s}$, $1\leq j\leq J, 1\leq \iota\leq\iota_{js}$, and call the resulting linearly ordered sets $Z_{1s},...,Z_{L_ss}$. Thus, $Z_{qs}$ is one of the sets $X_{j\iota s}$, and we assign $Z_{qs}$, same as its index $q$, with the same color as the one of $j$.\par
For every $s\leq S$ and every pair $(q,q')$, $1\leq q<q'\leq L_s$, let $\phi_{qq',s}$ be the detector, as given by Theorem \ref{the1}, associated with $X=Z_{qs}$ and $Y=Z_{q's}$; let also $\phi_{qq,s}\equiv 0$ and $\phi_{qq',s}(\omega)=-\phi_{q'q,s}(\omega)$, $1\leq q'<q\leq L_s$, so that
\begin{equation}\label{seq1}
\begin{array}{c}
\phi_{qq',s}(\cdot)\in\F,\,\phi_{qq',s}(\cdot)\equiv -\phi_{q'q,s}(\cdot),\\
 \int_\Omega\exp\{-\phi_{qq',s}(\omega)\}p_\mu(\omega)P(d\omega)\leq\epsilon_{qq',s}\,\,\forall \mu\in Z_{qs},\\
 \hbox{where\ } \epsilon_{qq',s}=\epsilon_{q'q,s}=\max\limits_{\mu\in Z_{qs},\nu\in Z_{q's}}\int_\Omega \sqrt{p_\mu(\omega)p_\nu(\omega)} P(d\omega)\in(0,1].\\
\end{array}
\end{equation}
\paragraph{$s$-closeness.}
Let $H_{qs}$ be the hypotheses on the probability density $p(\cdot)$ of an observation stating that this density is $p_\mu(\cdot)$ with some $\mu\in Z_{qs}$; the color of $H_{qs}$ is inherited from the set $Z_{qs}$, that is, it is the same as the color of all densities $p_\mu(\cdot)$ with $\mu\in Z_{qs}$. Let us say that hypothesis $H_{q's}$ is {\sl $s$-close} to hypothesis $H_{qs}$ (synonym: $q'$ is $s$-close to $q$), if either $H_{qs}$ and $H_{q's}$ are of the same color, or they are of different colors and
$\epsilon_{qq',s}>\delta_s$. The resulting set $\C=\C_s$ of pairs $(q,q')$ with $q'$ $s$-close to $q$ meets the only requirement imposed on $\C$ in section \ref{multiple}, specifically, it contains all diagonal pairs $(q,q)$ (since $\epsilon_{qq,s}=1$ by (\ref{seq1}) and $\delta_s<1$). In addition, $s$-closeness is symmetric:  $(q,q')\in\C_s$ if and only if $(q',q)\in C_s$ (this is due to $\epsilon_{qq',s}=\epsilon_{q'q,s}$).
\paragraph{Tests $T_s$.}
Now let us apply to the collection of hypotheses $\{H_{qs}:1\leq q\leq L_s\}$ and detectors $\{\phi_{qq',s}(\cdot):1\leq q,q'\leq L_s\}$, and to just defined $\C_s$ the construction from section \ref{multiple}, assuming that when deciding upon hypotheses $H_{1s},...,H_{L_ss}$, we have at our disposal $k$-repeated observation $\omega^k$, with a given $k$. Specifically, consider the optimization problem
\begin{equation}\label{seq9}
\Opt(k,s)=\min_{\alpha\in\bbr^{L_s\times L_s}}\left\{f_s(\alpha):=\max_p\sum_{q':(q,q')\not\in\C_s}\epsilon_{qq',s}^k\exp\{\alpha_{qq'}\}:\alpha=-\alpha^T\right\},\ \footnote{Essentially the same reasoning as in the proof of Proposition \ref{prop_alpha1} shows that $\Opt(k,s)$ is nothing but the spectral norm of the entry-wise nonnegative symmetric matrix
with the entries $\left\{\begin{array}{ll}\epsilon_{qq',s}^k,&(q,q')\not\in \C_s\\
0,&\hbox{otherwise}\\
\end{array}\right.$. Note, however, that now the optimal value in $(*)$ not necessarily is achieved.}
\end{equation}
Since $\epsilon_{qq',s}\leq \delta_s\in(0,1)$ when $(q,q')\not\in\C_s$,  $\Opt(k,s)$ goes to 0 as $k\to\infty$, so that the smallest $k=k(s)$ such that $\Opt(k,s)<\epsilon_s$ is well defined. And since $\Opt(k(s),s)<\epsilon_s$, problem (\ref{seq9}) with $k=k(s)$, solvable or not, admits a feasible solution $\bar{\alpha}^{(s)}$ such that $f_s(\bar{\alpha}^{(s)})\leq\epsilon_s$. Applying to detectors $\phi^\prime_{qq',s}(\cdot)=\phi_{qq',s}(\cdot)-\bar{\alpha}^{(s)}_{qq'}$ and to $\C=\C_s$ the construction from section \ref{multiple}, we get a test $T_{s}$ deciding on the hypotheses $H_{qs}$, $1\leq q\leq L_s$, via $k(s)$-repeated observation $\omega^{k(s)}$, with properties as follows:
\begin{quote}
Let $\omega_1,...,\omega_{k(s)}$ be drawn, independently of each other, from probability density $p_*(\cdot)$ obeying one of the hypotheses $H_{qs}$, say, the hypothesis $H_{q_*s}$. Then
\begin{enumerate}
\item the $p_*$-probability for $T_s$ not to accept $H_{q_*s}$ is at most $\epsilon_s$;
\item the $p_*$-probability of the event ``among the accepted hypotheses, there is a hypothesis $H_{qs}$ with $(q_*,q)\not\in\C_s$'' is $\leq\epsilon_s$.\\
(In general, in the latter claim ``$(q_*,q)\not\in \C_s$'' should be strengthened to  ``$(q_*,q)\not\in \C_s$ and $(q,q_*)\not\in \C_s$,'' but we are in the situation of symmetric closeness).
\end{enumerate}
\end{quote}
\paragraph{The last component of the setup.} We impose some restrictions on the last -- the $S$-th -- component of the setup, specifically, as follows. As we remember, $X_1,...,X_J$ are convex compact subsets of $\M$.  Let us agree that {\sl when $s=S$, the partitions of $X_j$ are trivial: $\iota_{jS}=1$ and $X_{j1S}=X_j$ for all $j\leq J$.} Moreover, we have assumed that the families  $P^1,...,P^I$ of probability densities are at positive Hellinger distances from each other, implying that {\sl with our $X_{j\iota S}$, all quantities $\epsilon_{qq',s}$, $1\leq q,q'\leq L_S=J$, with $q$, $q'$ of different colors are less than 1}. Our second agreement is that {\sl $\delta_S\geq\epsilon_{qq',S}$ whenever $q$, $q'$ are of different colors}, implying that {\sl $q$ and $q'$ are $S$-close if and only if $q$ and $q'$ are of the same color.}. With this specification of the  $S$-th component of our setup, we arrive at some $K=k(S)$. By reasons which will become clear in a moment, we are not interested in those components of our setup for which $k(s)>K$; if components with this property were present in our original setup, we can just eliminate them, reducing $S$ accordingly. Finally, we can reorder the components of our setup to make $k(s)$ nondecreasing in $s$. Thus, from now on we assume that
$$
k(1)\leq k(2)\leq...\leq k(S)=:K.
$$
\paragraph{Sequential test: construction.} We are ready to describe sequential test $\T$ corresponding to the outlined setup.  This test works by stages $s=1,2,...,S$ and formally (why ``formally,'', will be clear from the next sentence) makes decisions via observation $\omega^K$. Specifically, at stage $s$ we apply test $T_s$ to the initial fragment $\omega^{k(s)}$ of $\omega^K$. If the outcome of the latter test is acceptance of a nonempty set of hypotheses $H_{qs}$ {\sl and all these hypotheses are of the same color $i$}, test $\T$ accepts the hypothesis $\H_i$ and terminates, otherwise it proceeds to stage $s+1$ (when $s<S$) or terminates without accepting any hypothesis $(s=S)$.\par
Note that by construction $\T$ never accepts more than one of the hypotheses $\H_1,...,\H_I$.
\subsubsection{Analysis} For $\mu\in X=\bigcup\limits_{j=1}^J X_j$, let $s[\mu]\in\{1,2,...,S\}$ be defined as follows: for every $s$, $\mu$ belongs to (perhaps, several) of the sets $Z_{qs}$, $1\leq q \leq L_s$; let $Q_s[\mu]$ be the set of all $q$'s such that $\mu\in Z_{qs}$. Note that the color of every  $q\in Q_s[\mu]$ (recall that the sets $Z_{qs}$ and the corresponding values of $q$  have already been assigned colors) is the same as the color of $\mu$. Now, given $\mu\in X$, for some $s\leq S$ it may happen that
\begin{equation}\label{mayhappen}
\exists q\in Q_s[\mu]: \hbox{\ all $q'$ $s$-close to $q$ are of the same color as $q$.}
\end{equation}
In particular, the latter condition definitely takes place when $s=S$, since, as we have already seen,  $q$ and $q'$ are $S$-close if and only if $q$ and $q'$ are of the same color, whence for $s=S$ the conclusion in (\ref{mayhappen}) is satisfied whatever be the choice of $q\in Q_S[\mu]$. Now let $s[\mu]$ be the smallest $s$ such that (\ref{mayhappen}) takes place, and let the corresponding $q$ be denoted by $q[\mu]$; thus,
\begin{equation}\label{happens}
\begin{array}{l}
\forall \mu\in X:=\bigcup\limits_{j=1}^JX_j: \\
\multicolumn{1}{c}{\begin{array}{ll}s[\mu]\in\{1,2,...,S\},\,\, q[\mu]\in\{1,...,L_{s[\mu]}\},\,\, \mu\in Z_{q[\mu]s[\mu]};&(a)\\
\hbox{whenever $q'\in\{1,...,L_{s[\mu]}\}$ is $s[\mu]$-close to $q[\mu]$, $q'$ and $q[\mu]$ are of the same color.}&(b)\\
 \end{array}}\\
 \end{array}
\end{equation}

Our main result in this section is as follows.
\begin{proposition}\label{prop_sequential} Let $\omega^1,...,\omega^K$ be drawn, independently of each other, from some probability density $p_\mu(\cdot)$, $\mu\in \bigcup\limits_{j=1}^J X_j$, so that $\mu\in X_{j_*}$ for some $j_*\leq J$, and let $i_*$ be defined by the requirement $j_*\in \J_{i_*}$ (i.e., $i_*$ is the color of $\mu$). Then $p_\mu$-probability of the event
$$
\E=\left\{\omega^K: \hbox{\rm \begin{tabular}{l}$\T$, as applied to $\omega^K$, terminates not later than at the stage $s[\mu]$\\
 and accepts upon termintation the true hypothesis $\H_{i_*}$\\
 \end{tabular}} \right\}
$$
is at least $1-\epsilon$.
\end{proposition}
{\bf Proof.} Let $\mu$, $j_*$, $i_*$ be the entities from Proposition. For $1\leq s\leq S$, there exist $q_s\in\{1,...,L_s\}$ such that $\mu\in Z_{q_ss}$; without loss of generality, we can assume that $q_{s[\mu]}=q[\mu]$. For $1\leq s\leq S$ let
$$
\begin{array}{rcl}
\E_1^s&=&\{\omega^K: T_{s} \hbox{\ as applied to $\omega^{k(s)}$ does not accept the hypothesis $H_{q_ss}$}\},\\
\E_2^s&=&\{\omega^K: T_s \hbox{\ as applied to $\omega^{k(s)}$ accepts a hypothesis $H_{qs}$ with $q$ not $\C_s$-close to $q_s$}\}.\\
\end{array}
$$
From the just outlined properties of $T_s$ it follows that $p_\mu$-probabilities of $\E_1^s$ and $\E_2^s$ do not exceed $\epsilon_s$. Now let
$$
\E^*=\{\omega^K: \hbox{\ no one of the events $\E_1^s$, $\E_2^s$, $1\leq s\leq S$, takes place.}\},
$$
so that $p_\mu$-probability of $\E^*$ is at least $1-2\sum_{s=1}^S\epsilon_s = 1-\epsilon$. All we need is to verify that $\E\supset \E^*$, which is immediate. Indeed, let $s_*=s[\mu]$, $q_*=q[\mu]$ (so that $q_*=q_{s_*}$ due to $q_{s[\mu]}=q[\mu]$), and let $\omega^K\in \E^*$. By the latter inclusion, $\E_1^{s_*}$ does not take place, implying that \\
\indent (a) $T_{s_*}$ as applied to $\omega^{k(s_*)}$ accepts
$H_{q_*s_*}$.\\
Besides this, since $\omega^K\in \E^*$, we have $\omega^{k(s_*)}\not\in \E_2^{s_*}$, implying, by definition of $\E_2^{s_*}$, that all hypotheses $H_{qs_*}$ accepted by $T_{s_*}$ as applied to $\omega^{k(s_*)}$ are such that $q$ is $\C_{s_*}$-close to $q_*=q[\mu]$. By definition of $s_*=s[\mu]$ and $q_*=q[\mu]$, for all $q$ which are $\C_{s_*}$-close to $q_*$, the hypotheses $H_{qs_*}$ are of the same color as the accepted  by $T_{s_*}$, by (a), hypothesis $H_{q_*s_*}$. Thus, when $\omega^K\in\E^*$, we have
\\
\indent (b) the set of hypotheses $H_{qs_*}$ accepted by $T_{s_*}$ is nonempty, and all these hypotheses are of the same color, equal to the color of $\mu$.
\\
Invoking the description of $\T$, we conclude that when $\omega^K\in \E^*$, the test $\T$ terminates not later than at step $s_*=s[\mu]$ and the termination is productive -- some of the  hypotheses $\H_i$ indeed are accepted (in fact, in this case exactly one of the hypotheses $\H_i$ is accepted, since, as it was already mentioned, $\T$ never accepts more than one hypothesis).
\par
Now let $\omega^K\in \E^*$, and let $\bar{s}$ be the step at which $\T$ terminates. As we have already seen, $\bar{s}\leq s_*=s[\mu]$ and $\T$ terminates accepting exactly one hypothesis $\H_{\bar{i}}$, where $\bar{i}$ is a deterministic function of $\omega^{k(\bar{s})}$. Since $p_\mu$ obeys hypothesis $H_{q_{\bar{s}}\bar{s}}$ (by definition of $q_{\bar{s}}$) and $\omega^K\not\in \E_1^{\bar{s}}$, the test $T_{\bar{s}}$ as applied to $\omega^{k(\bar{s})}$ accepts the hypothesis $H_{q_{\bar{s}}\bar{s}}$ (which, by construction, has the same color $i_*$ as $\mu$). The latter observation combines with the termination rule for $\T$ to imply that the outcome of $\T$ in the case of $\omega^K\in\E^*$ is the hypothesis $\H_i$ with the same color $i$ as $H_{q_{\bar{s}}\bar{s}}$, that is, the hypothesis $\H_i$ with the same color as the color $i_*$ of $\mu$, that is, the true hypotheses $\H_{i_*}$, and that the outcome is obtained not later than at the stage $s[\mu]$. \qed
\subsubsection{An implementation} Observe that in order for the just described sequential test to  recover the ``non-sequential'' test built in section \ref{multiple} for the Special case of testing multiple unions, it suffices to utilize setup with $S=1$ (recall that our construction fully specifies the setup component with $s=S$). This being said, with our sequential test the running time (the number of observations used to make the inference) depends on the true distribution underlying the observations, and we can use the numerous degrees of freedom in our setup in order to save on the number of observations when the true distribution is ``deeply inside'' the true hypothesis. We are about to illustrate the corresponding options for our basic good observation schemes.
\paragraph{Situation.} Assume that our observation scheme is either Gaussian\footnote{We assume w.l.o.g. that the underlying covariance matrix is the unit one.}, or Poisson, or Discrete, see section \ref{sect:appl0}. We denote by $n$ the dimension of the associated parameter vector $\mu$: $\mu=[\mu_1;...;\mu_n]$.
\par
As in the already considered general case, the input to the sequential test is a risk level $\epsilon\in(0,1)$ and a collection of $J$ nonempty compact convex sets $X_j\subset \M$ colored in $I\geq 2$ colors (i.e., the set of indexes $\{1,...,J\}$ is split into $I\geq2$ non-overlapping nonempty sets $\J_1,...,\J_I$, $i$ being the common color of all sets $X_j$, $j\in\J_i$). Same as before, we assume that the sets $X_j$, $X_{j'}$ of different colors do not intersect. Given this input, we intend to specify the setup for sequential test deciding on the associated hypotheses $\H_i$, $1\leq i\leq I$.
\paragraph{Preliminaries.} Let $\psi(\mu,\nu):\M\times\M\to\bbr$ be the associated with the observation scheme in question function $\ln\left(\int_{\Omega}\sqrt{p_\mu(\omega)p_\nu(\omega)}P(d\omega)\right)$, so that
\begin{equation}\label{psiis}
\psi(\mu,\nu)=\left\{\begin{array}{ll}-{1\over 8}\|\mu-\nu\|_2^2,\,\mu,\nu\in\M=\bbr^n,&\hbox{Gaussian case}\\
-{1\over 2}\sum_{\ell=1}^n(\sqrt{\mu_\ell}-\sqrt{\nu_\ell})^2,\,\mu,\nu\in\M=\{x\in\bbr^n:x>0\},&\hbox{Poisson case}\\
\ln\left(\sum_{\omega=1}^n\sqrt{\mu_\omega\nu_\omega}\right),\,\mu,\nu\in\M=\{x\in\bbr^n:x>0,\sum_{\omega=1}^nx_\omega=1\},&\hbox{Discrete case}\\
\end{array}\right.
\end{equation}
Theorem \ref{the1} and its specifications presented in section \ref{sect:appl0} state that in the cases in question, given two nonempty convex compact subsets $X,Y$ of $\M$  and setting
$$
\Psi_{XY}=\max_{\mu\in X,\nu\in Y} \psi(\mu,\nu),
$$
the quantity $\exp\{\Psi_{XY}\}$ is exactly the risk of the detector yielded by Theorem \ref{the1} as applied to $X,Y$ and the observation scheme in question. Besides this, $\psi(\mu,\nu)$ clearly is smooth concave symmetric $(\psi(\mu,\nu)=\psi(\nu,\mu))$ function on its domain, and $\psi(\mu,\nu)<0$ whenever $\mu\neq\nu$.
\par
For $j\leq J$, we set
$$
\psi_j(\mu)=\max_{\nu\in X_j} \psi(\mu,\nu):\M\to\bbr.
$$
Since $\psi(\mu,\nu)$ is concave in $\mu,\nu\in\M$ and $X_j$ is a convex compact set, the functions $\psi_j(\cdot)$ are concave and continuous on $\M$.
  \par
We denote  by $\O$ the set of all ordered pairs $(j,j')$, $1\leq j,j'\leq J$, {\sl with $j,j'$ of different colors;} note that $(j,j')\in\O$ if and only if $(j',j)\in\O$. For $j\leq J$, we define  $\J^o_j$ as the set of all indexes $j'\leq J$ of color different from the one of $j$.
\par
 For a pair $j,j'$, $1\leq j,j'\leq J$, let
\begin{equation}\label{seq300}
\psi_{jj'}=\Psi_{X_jX_{j'}}:=\max_{\mu\in X_j,\nu\in X_{j'}}\psi(\mu,\nu)=\max_{\mu\in X_j}\psi_{j'}(\mu);
\end{equation}
note that $\psi_{jj'} <0$ when $(j,j')\in\O$, since in the latter case the convex compact sets $X_j$, $X_{j'}$ do not intersect, whence the objective in the right hand side optimization problem is negative at the compact feasible set of the problem. We set
\begin{equation}\label{seq330}
\chi=\min_{(j,j')\in\O}[-\psi_{jj'}],
\end{equation}
so that $\chi>0$.
\par
Finally, for  $(j,j')\in \O$ and a nonnegative $r$ we define {\sl $(jj'r)$-cut} as a linear inequality of the form $\ell(\mu)\leq 0$, $\mu\in \M$, such that
\begin{equation}\label{jjrineq}
\psi_{j'}(\mu)\leq -r\,\,\forall (\mu\in X_j: \ell(\mu)\leq 0).
\end{equation}
{\sl Example: default cuts.} For $(j,j')\in\O$, let $(a_{jj'},b_{jj'})$ form an optimal solution to the convex optimization problem (\ref{seq300}).
Setting
$$
e_{jj'} =\nabla_\mu\psi(a_{jj'},b_{jj'}), \,\,f_{jj'}=\nabla_\nu\psi(a_{jj'},b_{jj'}),\,\,(j,j')\in \O,
$$
and invoking optimality conditions for (\ref{seq300}) along with concavity of $\psi(\cdot,\cdot)$, we get for all $(j,j')\in\O$:
\begin{equation}\label{seq301}
\forall (\mu\in X_j,\nu\in X_{j'}):\left\{\begin{array}{ll}
e_{jj'}^T[\mu-a_{jj'}]\leq 0&(a)\\
f_{jj'}^T[\nu-b_{jj'}]\leq 0&(b)\\
\psi(\mu,\nu)\leq \psi_{jj'} +e_{jj'}^T[\mu-a_{jj'}]+f_{jj'}^T[\nu-b_{jj'}]&(c)\\
\end{array}\right.
\end{equation}
It follows that setting
$$
\ell_{jj'}(\mu)=\psi(a_{jj'},b_{jj'})+\langle e_{jj'},\mu-a_{jj'}\rangle-r,
$$
we get an affine function of $\mu\in\M$ which upper-bounds $\psi_{j'}(\cdot)-r$ on $X_j$.
\begin{quote} Indeed,
$$
\begin{array}{l}
(\nu\in X_{j'},\mu\in X_j)\Rightarrow
\psi(\mu,\nu)\leq \psi(a_{jj'},b_{jj'})+\langle e_{jj'},\mu-a_{jj'}\rangle +\langle f_{jj'},\nu-b_{jj'}\rangle\\
\leq \psi(a_{jj'},b_{jj'})+\langle e_{jj'},\mu-a_{jj'}\rangle
\end{array}
$$
(we have used (\ref{seq301}.$c$,$b$)). Taking in the resulting inequality the supremum over $\nu\in X_{j'}$, we arrive at $\ell_{jj'}(\mu)\geq\psi_{j'}(\mu)-r$, $\mu\in X_j$.
\end{quote}
The bottom line is that $\ell_{jj'}(\mu)\leq0$ is a $(jj'r)$-cut; we shall call this cut {\sl default}.
\paragraph{Specifying setup for the sequential test.}
The setup for our sequential test is as follows.
\begin{enumerate}
\item We select a sequence of positive integers $\{\bar{k}(s)\}_{s=1}^\infty$ satisfying
\begin{equation}\label{monotonicity}
\bar{k}(1)=1,\,\,\bar{k}(s)<\bar{k}(s+1)\leq 2\bar{k}(s),\,\,s=1,2,...
\end{equation}
and specify $S$ as the smallest positive integer such that
\begin{equation}\label{seq10}
\bar{k}(S)> {1\over \chi}\ln\left(2.1SJ^2/\epsilon\right);
\end{equation}
($S$ is well defined due to $\bar{k}(s)\geq s$).
\par
For $1\leq s\leq S$, we set
\begin{equation}\label{seq11}
\begin{array}{c}
\epsilon_s={\epsilon\over 2S},\,
r(s)=\ln(2.1SJ^2/\epsilon)/\bar{k}(s),\,\,\delta_s=\exp\{-r(s)\}.\\
\end{array}
\end{equation}
\item
For every $j\leq J$ and every $s$, $1\leq s\leq S$, we specify a number of closed convex subsets $X_{j\iota s}$, $1\leq \iota\leq \iota_{js}$, of $X_j$ as follows. For every pair $(j,j')\in\O$, we select somehow
a $(jj'r(s))$-cut $\ell_{jj's}(\cdot)\leq0$ and set
\begin{equation}\label{partition}
\begin{array}{rcl}
X^{j'}_{js}&=&\{\mu\in X_j: \ell_{jj's}(\mu) \geq 0\},\, j'\in\J^o_j,\\
X^j_{js}&=&\{\mu\in X_j: \ell_{jj's}(\mu)\leq0,\,\, j'\in \J^o_j\},\\
\end{array}
\end{equation}
 thus getting a number of convex compact sets; eliminating from their list all sets which are empty,  we end up with a number $\iota_{js}\leq J$ of {\sl nonempty} convex compact sets $X_{j\iota s}$, $1\leq\iota\leq\iota_{js}$, with $X_j$ being their union.
 \par
Observe that $r(S)<\chi$ by (\ref{seq10}), and for ($j,j')\in\O$ we clearly have
$$
\max_{\mu\in X_j}\psi_{j'}(\mu)=\psi_{jj'}\leq -\chi
$$
implying that $\ell_{jj'S}(\mu)\equiv -1$ are legitimate $(jj'S)$-cuts. These are exactly the cuts we use when $s=S$.
\end{enumerate}
We claim that the just defined entities form a legitimate setup for a sequential test. All we need in order to justify this claim is to verify that the $S$-component of our setup is as required, that is, that (a) for every $j\leq J$, $\iota_{jS}=1$, whence $X_{j1S}=X_j$ and $L_S=J$, and that (b) $q\leq L_S=J$ is $S$-close to $q'\leq L_S=J$ if and only if $q$ and $q'$ are of the same color.
\begin{quote}
To verify (a), note that $\ell_{jj'S}(\cdot)\equiv -1$ whenever $(j,j')\in\O$, implying that $X^{j'}_{jS}=\emptyset$ when $j'\in\J^o_j$ and $X^j_{jS}=X_j$, as claimed in (a).
To verify (b), note that as it was already mentioned, for $j$, $j'$ of different colors, the risk of the detector yielded by Theorem \ref{the1} as applied to the sets $X=X_j$, $Y=X_{j'}$, is $\exp\{\psi_{jj'}\}$, that is, this risk is $\leq \exp\{-\chi\}$. Invoking the already verified (a), we conclude that $\epsilon_{qq',S}\leq\exp\{-\chi\}$ whenever $1\leq q,q'\leq L_S=J$ and $q,q'$ are of different colors. As we have seen, $r(S)<\chi$, whence $\delta_S=\exp\{-r(S)\}>\exp\{-\chi\}$. The bottom line is that whenever $1\leq q,q'\leq J$ and $q,q'$ are of different colors, we have $\epsilon_{qq',S}<\delta_S$; this observation combines with the definition of $S$-closeness to imply that $q,q'$ are $S$-close if and only if $q,q'$ are of different colors, as claimed in (b).
\end{quote}
The legitimate setup we have presented induces a sequential test, let it be denoted by $\T$. We are about to analyse the properties of this test.
\paragraph{Analysis.}
Our  first observation is that for the sequential setup we have presented one has $k(s)\leq \bar{k}(s)$, $1\leq s\leq S$. To verify this claim, we need to check that when $k=\bar{k}(s)$, we have $\Opt(k,s)<\epsilon_s={\epsilon\over 2S}$. As we have already mentioned  (see (\ref{seq9})), $\Opt(k,s)$ is the spectral norm of the entrywise nonnegative symmetric matrix $D^{ks}$ of the size $L_s\times L_s$ with entries not exceeding $\delta_s^k$, and by construction, $L_s\leq J^2$, so that the spectral norm of $D^{ks}$ does not exceed $J^2\delta_s^k$. The latter quantity indeed is $< \epsilon_s$ when $k=\bar{k}(s)$, see (\ref{seq11}).
\paragraph{Worst-case performance.} In the analysis to follow, we assume that $0<\epsilon\leq 0.05$.
\par
By Proposition \ref{prop_sequential}, the sequential test $\T$ always accepts at most one of the hypotheses $\H_1,...,\H_I$, and the probability not to accept the true hypothesis is at most $\epsilon$; moreover, the number of observations used by $\T$ never exceeds $k(S)\leq \bar{K}:=\bar{k}(S)$.
 On the other hand, from the definition of $d$ and  Proposition \ref{newcol1} it follows that in order for a whatever test to decide on the hypotheses $\H_1,...,\H_I$ with risk $\epsilon$ via stationary repeated observations, the number of these observations should be {\sl at least}
\begin{equation}\label{Kstarequals}
K_*=\left[{{1\over 2}\ln(1/\epsilon)-\ln(2)\over\ln(1/\epsilon)}\right]{\log(1/\epsilon)\over \chi}\geq {\log(1/\epsilon)\over 4\chi}.
\end{equation}
It turns out that {\sl unless $\chi$ is ``astronomically small'', $\bar{K}$ is within logarithmic factor of $K_*$}, implying that as far as the worst-case performance  is concerned, our sequential test $\T$ is not that bad. The precise statement is as follows:
\begin{proposition}\label{proploglog} Let $\chi>0$, $J\geq2$, $\epsilon\leq 0.05$ satisfy, for some $\kappa\geq1$,  the relation
\begin{equation}\label{letdsatisfy}
{1\over \chi}\leq \Theta^\kappa,\,\,\Theta:=2.1J^2/\epsilon.
\end{equation}
Then
\begin{equation}\label{ifdsatisfies}
\bar{K}\leq \max\left[1,2\kappa{\ln(\Theta)\over
\chi}\right].
\end{equation}
\end{proposition}
For all practical purposes we can assume something like $\chi\geq 10^{-7}$, otherwise the {\sl lower} bound $K_*$ on the number of observations required by $(1-\epsilon)$-reliable test would be impractically large. Assuming $\chi\geq 10^{-7}$, (\ref{letdsatisfy}) is satisfied with $\kappa=3.15$ (recall that $\epsilon\leq 0.05$ and $J\geq2$). Thus, for all practical purposes we may treat the quantity $\kappa$ from the premise of Proposition \ref{proploglog} as a moderate absolute constant, implying that the upper bound $\bar{K}$ on the worst-case observation time of our $(1-\epsilon)$-reliable sequential test $\T$ indeed is within a logarithmic factor $O(1)\ln(J)/\ln(1/\epsilon)$ of the lower bound $K_*$ on the worst-case observation time of an ``ideal'' $(1-\epsilon)$-reliable test.
\begin{remark}\label{remloglog} It is easily seen that when $\bar{k}(s)$ grows with $s$ as rapidly as allowed by {\rm (\ref{monotonicity})}: $\bar{k}(s)=2^{s-1}$, the result completely similar to the one of Proposition \ref{proploglog} holds true in a much wider than {\rm (\ref{letdsatisfy})} range of values of $d$, specifically, in the range
$\ln(1/\chi)\leq C\Theta$, for a whatever constant $C\geq1$; in this range, one has $\bar{K}\leq \bar{C}\max[1,\ln(\Theta)/\chi]$, with $\bar{C}$ depending solely on $C$.
\end{remark}
\paragraph{Actual performance.} For $\mu\in X:=\bigcup\limits_{j=1}^J X_j$ let $s_*(\mu)$ be the smallest $s\leq S$ such that for some $j\leq J$ it holds $\mu\in X_{js}^j$ (see (\ref{partition})). Equivalently:
\begin{equation}\label{seq22}
s_*(\mu)=\min\left\{s: \exists j\leq J: \mu\in X_j \ \&\ \ell_{jj's}(\mu)\leq0\,\,\forall j'\in\J^o_j.
\right\}\end{equation}
Note that $s_*(\mu)$ is well defined -- we have already seen that $X^{j'}_{jS}=\emptyset$ whenever $j'\in\J^o_j$, so that $s=S$ is feasible for the right hand side problem in (\ref{seq22}).
\begin{proposition}\label{propactual}
One has $$s[\mu]\leq s_*(\mu).$$
\end{proposition}
Postponing for a moment the verification, let us look at the consequences. Assume that the observations are drawn from  density $p_\mu(\cdot)$, $\mu\in\cup_jX_j$. By Proposition \ref{prop_sequential}, with $p_\mu$-probability $\geq1-\epsilon$ our sequential test $\T$ terminates in no more than $s[\mu]\leq s_*(\mu)$ steps and upon termination recovers correctly the color $i[\mu]$ of $\mu$ (i.e., accepts the true hypothesis $\H_{i[\mu]}$, and only this hypothesis). Thus, if $\mu$ is ``deeply inside'' one of the sets $X_{j}$, meaning that $s_*(\mu)$ is much smaller than $S$, our sequential test will, with reliability $1-\epsilon$, identify correctly the true hypothesis $\H_{i[\mu]}$ much faster than in $S$ stages.
\begin{quote}
{\small
It remains to verify that $s[\mu]\leq s_*(\mu)$. Let $s_*=s_*(\mu)$, and let $j_*$ be such that $\mu\in X_{j_*s_*}^{j_*}$, see the definition of $s_*(\mu)$.
The set $X_{j_*s_*}^{j_*}\ni\mu$ is some $Z_{q_*s_*}$ with $q_*\in Q_{s_*}[\mu]$. We claim that\\[3pt]
 $\quad$ (!) {\sl all $q'$ which are $s_*$-close to $q_*$ are of the same color as $q_*$};\\[3pt]
note that the validity of (!) means that setting $s=s_*$, $q=q_*$, we meet the requirements in (\ref{mayhappen}), implying, by definition of $s[\mu]$, our target relation $s[\mu]\leq s_*$.\par
To verify (!), let $q'$ be $s_*$-close to $q_*$, so that by definition of $s$-closeness either\\
 $\quad$(a) $q'$ and $q_*$ are of the same color, or\\
  $\quad$(b) $q'$ and $q_*$ are of different colors and $\epsilon_{q'q_*,s_*}>\delta_{s_*}=\exp\{-r(s_*)\}$;\\
  all we need to prove is that (b) in fact is impossible. Assume, on the contrary, that (b) takes place, and let us lead this assumption to a contradiction. As it was already mentioned, $\epsilon_{q'q_*,s_*}=\exp\{\lambda\}$, where
  $$
  \lambda=\max_{\mu\in Z_{q_*s_*},\nu\in Z_{q's_*}}\psi(\mu,\nu).
  $$
  Thus, we are in the case when the color $i$ of $q'$
 differs from the color $i_*$ of $q_*$ and
 \begin{equation}\label{contradict}
\lambda:=\max_{\mu\in Z_{q_*s_*},\nu\in Z_{q's_*}}\psi(\mu,\nu)>-r(s_*).
 \end{equation}
 Since $q'$ is of color $i$, we have $Z_{q's_*}\subset X_{j'}$ for some $j'\in\J_i$, and since $i\neq i_*$, $j'$ and $j_*$ are of different colors, or, equivalently, $j'\in\J^o_{j_*}$. Now, by the definition of $q_*$ we have
 \begin{equation}\label{Pbyconstruction}
Z_{q_*s_*}=X^{j_*}_{j_*s_*}=\{\mu\in X_{j_*}: \ell_{j_*js_*}(\mu)\leq 0\,\,\forall j\in\J^o_{j_*}\},
\end{equation}
whence, taking into account that $j'\in \J^o_{j_*}$ and that $\ell_{j_*j's_*}(\cdot)\leq0$ is $(j_*j's_*)$-cut, we get
$$
\mu\in Z_{q_*s_*}\Rightarrow \psi_{j'}(\mu)\leq -r(s_*).
$$
The latter relation, due to $Z_{q's_*}\subset X_{j'}$ and that $\psi_{j'}(\cdot)=\max_{\nu\in X_{j'}}\psi(\mu,\nu)$, implies that
$$
\max_{\mu\in Z_{q_*s_*},\nu\in Z_{q's_*}}\psi(\mu,\nu)\leq-r(s_*),
$$
which contradicts (\ref{contradict}). We have arrived at the desired contradiction.}
\end{quote}
\paragraph{Selecting the cuts.} The major ``degrees of freedom'' in the just described construction are the cuts $\ell_{jj's}(\cdot)$. The simplest option here is to use the default cuts. A more advanced option is to select the cuts in order to improve the performance of the resulting test. Informally speaking, we would like these cuts to result in as small sets $X^{j'}_{js}$, $j'\in\J^o_j$ as possible,  thus increasing chances for the probability density underlying our observations to belong to $X^j_{js}$ for a small value of $s$ (such a value, as we remember, with probability $1-\epsilon$ upper-bounds the number of stages before termination). A natural formal goal here would be to minimize, given $s$ and $(j,j')\in \O$, the $m_j$-dimensional volume of the set $X^{j'}_{js}$, where $m_j$ is the dimension of $X_j$. Setting $m=m_j$, $X=X_j$ and $Y=\{x\in X: \psi_{j'}(x)\geq-r(s)\}$, this boils down to solving the following geometric problem:
\begin{quote}
(*): {\sl Given a convex compact set $X$ with a nonempty interior in some $\bbr^m$ and a nonempty convex compact subset $Y$ of $X$, find an affine function $\ell(x)$ such that the  linear inequality
$\ell(x)\geq0$ is valid on $Y$ and minimizes, under this requirement, the average linear size
$$
\Size(X_{\ell(\cdot)}):=\left[\mes_m(X_{\ell(\cdot)})\right]^{1/m}
$$
of the set $X_{\ell(\cdot)}=\{x\in X: \ell(x)\geq 0\}$.}
\end{quote}
Problem (*) seems to be heavily computationally intractable. We are about to present a crude {\sl sub}optimal solution to (*).
\par
We can assume w.l.o.g. that $Y$ intersects the interior of $X$ \footnote{Indeed, otherwise we can specify $\ell(\cdot)$ as an affine function separating $Y$ and $X$, so that $\ell(x)\geq0$ when $x\in Y$ and $\ell(x)\leq0$ on $X$ and  $\ell(\cdot)$ is nonconstant on $X$. Clearly, $\ell(\cdot)$ is a feasible solution to (*) with $\Size(X_{\ell(\cdot)})=0$.}. Let us equip $X$ with a $\vartheta$-self-concordant barrier $F(\cdot)$ \footnote{That is, $F$ is a three times continuously differentiable strictly convex function on $\inter X$ which is an interior penalty for $X$ (i.e., diverges  to $+\infty$ along every sequence of interior points of $X$ converging to a boundary point of $X$), and, in addition, satisfies specific differential inequalities; for precise definition and related facts to be used in the sequel, see \cite{NN1994}. A reader will not lose much when assuming that $X$ is a convex compact set given by a strictly feasible system $f_i(x)\leq b_i$, $1\leq i\leq m$,  of convex quadratic inequalities, $F(x)=-\sum_{i=1}^m\ln(b_i-f_i(x))$ and $\vartheta=m$.}. Then the minimizer $\bar{x}$ of $F$ on $Y$ is uniquely defined, belongs to $\inter X$  and can be found efficiently by solving the (solvable) convex optimization problem $\min_{x\in Y\cap X}F(x)$. By optimality conditions for the latter problem, the affine function
$\bar{\ell}(x)=\langle \nabla F(\bar{x}),x-\bar{x}\rangle$ is nonnegative on $Y$ and thus is a feasible solution to (*). Further, from the basic facts of the theory of self-concordant barriers it follows that
\begin{itemize}
\item [A.] The {\sl Dikin ellipsoid} of $F$ at $\bar{x}$ -- the set $D=\{x\in X: \langle x-\bar{x},\nabla^2F(\bar{x})[x-\bar{x}]\rangle \leq 1\}$ -- is contained in $X$;
\item [B.] The set $X^+=\{x\in X: \langle x-\bar{x},\nabla F(\bar{x})\rangle \geq0\}$ is contained in the set\footnote{What is $\vartheta+2\sqrt{\vartheta}$ below is $3\vartheta$ in \cite{NN1994}; refinement $3\vartheta\to\vartheta+2\sqrt{\vartheta}$ is due to F. Jarre, see \cite[Lemma 3.2.1]{NLN}.}
$$D^+=\{x\in\bbr^m:  \langle x-\bar{x},\nabla^2F(\bar{x})[x-\bar{x}]\rangle \leq \rho^2:=(\vartheta+2\sqrt{\vartheta})^2,\langle x-\bar{x},\nabla F(\bar{x})\rangle\geq0\}.$$
\end{itemize}
From B it follows that
\begin{equation}\label{volumeissmall}
\Size(X_{\bar{\ell}(\cdot)})\leq \Size(D^+)=\rho\Size(D'),\,\,D'=\{x\in D: \langle x-\bar{x},\nabla F(\bar{x})\rangle\geq0\}.
\end{equation}
On the other hand, if $\ell(\cdot)$ is a feasible solution to (*), then the set $X_{\ell(\cdot)}$ contains $Y$ and thus contains $\bar{x}$, implying that
$\ell(\bar{x})\geq0$. Consequently, by A we have
$$
\bar{D}:=\{x\in D:\ell(x)\geq0\}\subset \{x\in X:\ell(x)\geq0\}=X_{\ell(\cdot)},
$$
whence
$$
\Size(X_{\ell(\cdot)})\geq\Size(\bar{D}).$$
         Since $\ell(\bar{x})\geq0$ and $\ell(\cdot)$ is affine, we have $\bar{D}\supset D'':=\{x\in X:\ell(x)\geq\ell(\bar{x})\}$, and the $m$-dimensional volume of $D''$
         is at least half of the $m$-dimensional volume of $D'$ (since every one of the sets $D''$, $D'$ is either the entire ellipsoid $D$, or is the intersection of $D$ with half-space with the boundary hyperplane passing through the center of $D$).  It follows that $\Size(\bar{D})\geq\Size(D'')\geq 2^{-1/m}\Size(D')$. Thus, for every feasible solution $\ell(\cdot)$ to (*) it holds
$$
\Size(X_{\ell(\cdot)})\geq \Size(\bar{D})\geq 2^{-1/m}\Size(D')\geq 2^{-1/m}\rho^{-1} \Size(X_{\bar{\ell}(\cdot)}),
$$
where the concluding $\geq$ is due to (\ref{volumeissmall}). We see that the feasible solution $\bar{\ell}(\cdot)$ to (*) (which can be found efficiently) is optimal within the factor $2^{1/m}\rho=2^{1/m}[\vartheta+2\sqrt{\vartheta}]$. This factor is moderate when $X$ is an ellipsoid (or the intersection of $O(1)$ ellipsoids), and can be unpleasantly large when $\vartheta$ is large; this, however, is the best known to us computationally tractable approximation to the optimal solution of (*).

\paragraph{Gaussian case.} In Gaussian case {\sl with default cuts}, the quantity $s_*(\mu)$ which, as we have just seen, is a $(1-\epsilon)$-reliable upper bound on the number of stages in which $\T$ recognizes $(1-\epsilon)$-reliably  the true hypothesis $\H_{i[\mu]}$ provided the observations are drawn from  $p_\mu(\cdot)$, admits a transparent geometric upper bound. Observe, first, that in Gaussian case we have
$$
\chi={1\over 8}\min_{(j,j')\in \O}\min_{a\in X_j,b\in X_{j'}}\|a-b\|_2^2.
$$
Now let $\rho(\mu)$ be the largest $\rho$ such the $\|\cdot\|_2$-ball of radius $\rho$ centered at $\mu$ is contained in certain $X_j$. We claim that
\begin{equation}\label{gauss222}
s[\mu]\leq s_*[\mu]\leq \bar{s}(\mu):=\min\{s\leq S: r(s)\leq \chi+ \sqrt{\chi/2}\rho(\mu)\},
\end{equation}
meaning that the deeper $\mu$ is ``inside'' one of $X_j$ (the larger is $\rho(\mu)$), the smaller is the number of observations needed for $\T$ to identify correctly the color of $\mu$.
\begin{quote}
{\small
Justification of (\ref{gauss222}) is as follows. The first inequality in (\ref{gauss222}) has already been proved. To prove the second inequality, observe, first, that $\bar{s}:=\bar{s}(\mu)$ is well defined (indeed, as we have seen, $r(S)<\chi)$. Let $j_*$ be such that the $\|\cdot\|_2$-ball $B$ of radius $\rho(\mu)$ centered at $\mu$ is contained in $X_{j_*}$, and let $j\in\J^o_{j_*}$. By (\ref{seq301}.$a$), we have $e_{j_*j}^T[\mu'-a_{j_*j}]\leq0$
for all $\mu'\in X_{j_*}$ and thus for all $\mu'\in B$, and therefore $e_{j_*j}^T[\mu-a_{j_*j}]\leq-\rho(\mu)\|e_{j_*j}\|_2$, whence
\begin{equation}\label{seq401}
\forall (j\in \J^o_{j_*}): \psi_{j_*j}+e_{j_*j}^T[\mu-a_{j_*j}]\leq\psi_{j_*j}-\rho(\mu)\|e_{j_*j}\|_2.
\end{equation}
Denoting by $(a,b)$ an optimal solution to the problem $\min_{\mu\in X_{j_*},\nu\in X_{j}}\|\mu-\nu\|_2$, we clearly have $\psi_{j_*j}=-\|a-b\|_2^2/8$, $\|e_{j_*j}\|_2=\|a-b\|_2/4$, and $\|a-b\|^2/8\geq \chi$ whenever $j\in\J^o_{j_*}$ by origin of $\chi$. Thus, (\ref{seq401}) implies that
$$
\forall (j\in \J^o_{j_*}): \psi_{j_*j}+e_{j_*j}^T[\mu-a_{j_*j}]\leq -\chi-\rho(\mu)\sqrt{\chi/2},
$$
that is, $\mu\in X_{j_*\bar{s}}^{j_*}$ due to $r(\bar{s}) \leq \chi+\rho(\mu)\sqrt{\chi/2}$, implying that $\bar{s}\geq s_*(\mu)$, as claimed.}
\end{quote}
\paragraph{Numerical illustration.} The following numerical experiment provides some impression of the power of sequential testing in the Gaussian case with default cuts. In this experiment, we are given $J=4$ sets $X_j\subset\bbr^2$; $X_1$ is the square $\{0.01\leq x_1,x_2\leq 1\}$, $X_2,X_3,X_4$ are obtained from $X_1$ by reflections w.r.t. the coordinate axis and the origin. The partition of the index set into groups $\J_i$ is trivial -- these groups are just the points comprising $\J=\{1,2,3,4\}$, so that our goal is to recognize which of the sets $X_j$ contains the mean of our observation. Figure \ref{figseqprof} presents the graph of the {\sl logarithm} of the $0.99$-reliable upper bound on the number of observations used by the sequential test with $S=20$, $\bar{k}(s)=2^{s-1}$ and  and $\chi=5.0$e$-5$ as a function of the mean $\mu$ of the observation. We see that the savings from sequential testing are quite significant. The related numbers are as follows: when selecting $\mu$ in $\bigcup\limits_{j=1}^4 X_j$ at random according to the uniform distribution, the empirical average of the number of observations before termination is as large as $1.6\cdot10^5$, reflecting the fact that $X_j$ are pretty close to each other. At the same time, the median number of observations before termination is just 154, reflecting the fact that in our experiment $\mu$, with reasonably high probability, indeed is deeply inside the set $X_j$ to which $\mu$ happens to belong.
\begin{figure}
$$
\begin{array}{cc}
\begin{array}{c}
\\[-30pt]
\epsfxsize=140pt\epsfysize=100pt\epsffile{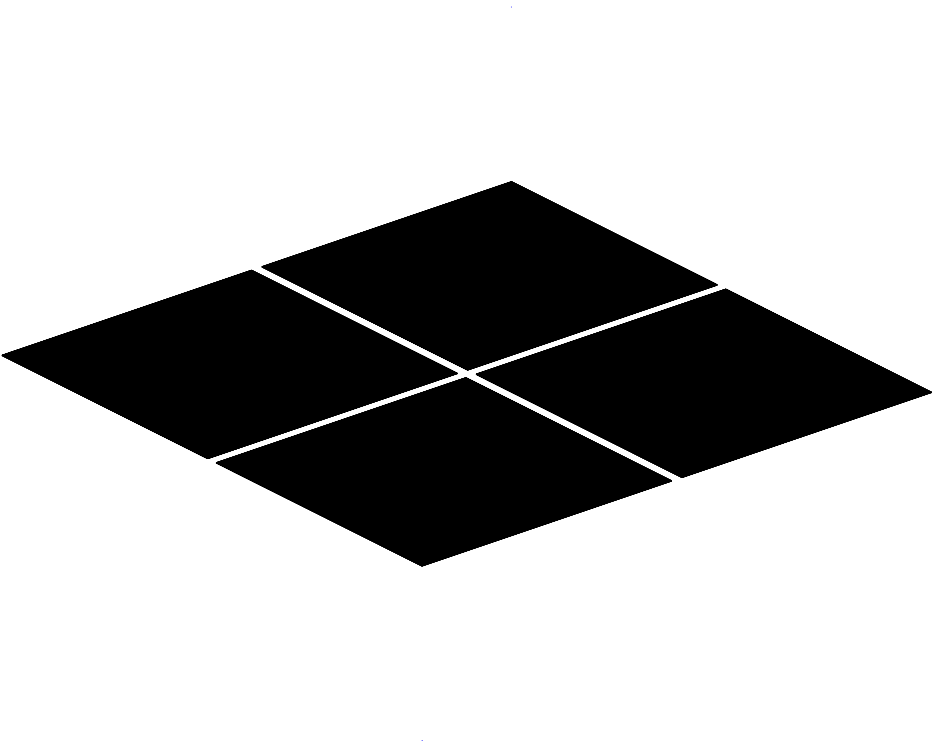}\\
\end{array}&\epsfxsize=230pt\epsfysize=140pt\epsffile{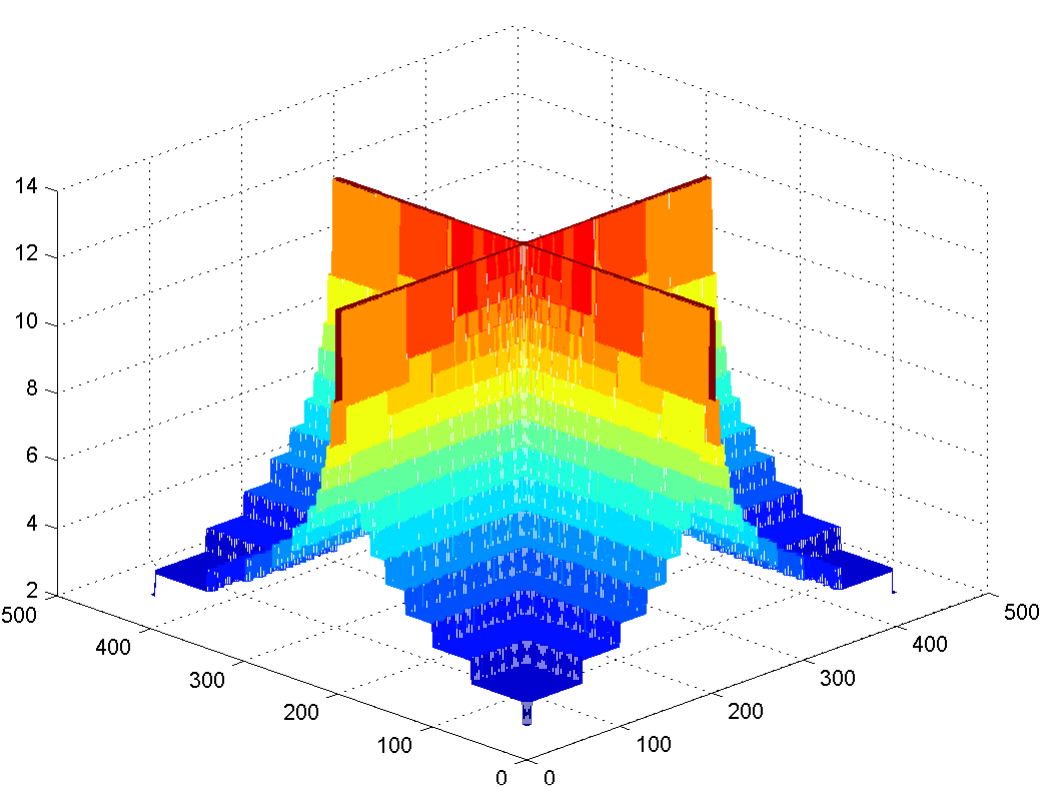}\\
\end{array}
$$
\caption{\label{figseqprof} ``Power'' of sequential test, Gaussian case. Left: $X_1,X_2,X_3,X_4$. Right: {\sl logarithm} of $k(s_*(\mu))$ as a function of $\mu\in\bigcup\limits_{i=1}^4X_i$, $\epsilon=0.01$.}
\end{figure}
\paragraph{Comparing cut policies.} To get impression on the effect of different cut policies (default cuts vs. cuts based on self-concordant barriers), consider the Gaussian case with $J=2$, $X_1=\{x\in\bbr^n:\delta\leq x_1\leq 1+\delta,0\leq x_i\leq1,\,2\leq i\leq n\}$, $X_2=\{x\in\bbr^n:-1\leq x_i\leq 0,\,1\leq i\leq n\}$, $I_1=\{1\}$, $J_2=\{2\}$ (i.e., color of $X_1$ is 1, and color of $X_2$ is 2), and let $\bar{k}(s)=2^{s-1}$, $s=1,2,...$.  We have
$$
\begin{array}{l}
\psi(\mu,\nu)=-{1\over 8}\|\mu-\nu\|_2^2,\,\,
\psi_1(\nu)=-{1\over 8}\left[(\nu_1-\delta)^2+\sum_{i=2}^n\nu_i^2\right],\,\,\psi_2(\mu)=-{1\over 8}\|\mu\|_2^2;\\
\psi_{1,2}=\psi_{2,1}=-{1\over 8}\delta^2,\,\,d:=\min[-\psi_{1,2},-\psi_{2,1}]={\delta^2\over 8};\\
a_{1,2}=b_{2,1}=e:=[1;0;...;0],\,\,b_{1,2}=a_{2,1}=0\in\bbr^n,\,\,
e_{1,2}=f_{2,1}=e,\,\,f_{1,2}=e_{2,1}=0.\\
\end{array}
$$
Further,
$$
r(s)={2\ln(8.4S/\epsilon)\over 2^s},\,1\leq s\leq S,
$$ where $\epsilon$ is the target risk, and $S$ is the smallest positive integer such that
$r(S)$ as given by the above formula is $<\chi={\delta^2\over 8}$
\par
Now, the default cuts are
$$
\begin{array}{ll}
(1,2,r):&\ell_{1,2,r}(\mu):=-{1\over 8}\delta^2-{\delta\over 4}(\mu_1-\delta)+r\leq0\\
(2,1,r):&\ell_{2,1,r}(\nu):=-{1\over 8}\delta^2+{\delta\over 4}\nu_1+r\leq0\\
\end{array}
$$
and the sets $Z_{qs}$ associated with the default cuts are:
$$
\begin{array}{ll}
\hbox{color \# 1:}&\left\{\begin{array}{rcl}Z_{1s}&=&X_{1,s}^1=\{\mu:\max[\delta,4{r(s)\over\delta}+{\delta\over 2}]\leq\mu_1\leq 1+\delta\ \&\ 0\leq \mu_i\leq1,\,2\leq i\leq n\},\\
Z_{2s}&=&X_{1,s}^2=\{\mu:\delta\leq \mu_1\leq \min[\delta,4{r(s)\over\delta}+{\delta\over 2}]\ \& \ 0\leq mu_i\leq 1,\,2\leq i\leq n\};\\
\end{array}\right.\\
\hbox{color \# 2:}&\left\{\begin{array}{rcl}Z_{3s}&=&X_{2,s}^2=\{\nu:-1\leq \nu_1\leq \min[0,-{4r(s)\over \delta}+{\delta\over 2}]\ \&\ -1\leq \nu_i\leq 0,\,2\leq i\leq n\},\\
Z_{4s}&=&X_{2,s}^1=\{\nu:\min[0,-{4r(s)\over\delta}+{\delta\over 2}]\leq \nu_1\leq0\ \&\ -1\leq \nu_i\leq 0,\,2\leq i\leq n\}.\\
\end{array}\right..\\
\end{array}
$$
The ``good'' sets here are $Z_{1s}$ and $Z_{3s}$, meaning that if a single observation is distributed  according to $\cN(\mu,I_d)$ with $\mu$ belonging to $Z_{1s}$ or $Z_{3s}$, our sequential test $\T$ with probability at least $1-\epsilon$ terminates  in course of the first $s$ stages with correct conclusion on the color of $\mu$. We see also that when $\delta\ll 1$ and $s$ is ``moderate,'' meaning that $\delta\gg r(s)\gg\delta^2$, the bad sets $Z_{2s}$, $Z_{4s}$, in terms of their $n$-dimensional volume, form $O(1)r(s)/\delta$-fractions of the respective boxes $X_1$, $X_2$, provided default cuts are used. When using ``smart'' cuts -- those induced by self-concordant barriers for $X_1$, $X_2$ -- it is immediately seen that in the range $1\gg\delta\gg r(s)\gg\delta^2$ (what exactly $\gg$ means, depends on $n$), the bad sets $Z_{2s}$ and $Z_{4s}$ become simplexes:
$$
\begin{array}{rcl}
Z_{2s}&=&\{\mu=[\delta;0;...;0]+[\lambda_1;...;\lambda_n]:\lambda\geq0,\sum_{i=1}^nc_i\lambda_i\leq \sqrt{r(s)n}\},\\
Z_{4s}&=&\{\nu=-[\lambda_1;...;\lambda_n]:\lambda\geq0,\sum_{i=1}^nc_i\lambda_i\leq \sqrt{r(s)n}\},\\
\end{array}
$$
where $c_i=c_i(n,s)$ are of order of 1; the good sets $Z_{1s}$, $Z_{3s}$ are the closures of the complements of these simplexes to the respective boxes $X_1$, $X_2$. The new bad sets are much smaller than the old ones; e.g., their volumes are of order of $r^{n/2}(s)$ -- much smaller than volumes $O(r(s)/\delta)$ of the old bad sets, see Figure \ref{newfig} and Table \ref{newtable}.
\begin{figure}
$$
\epsfxsize=180pt\epsfysize=170pt\epsffile{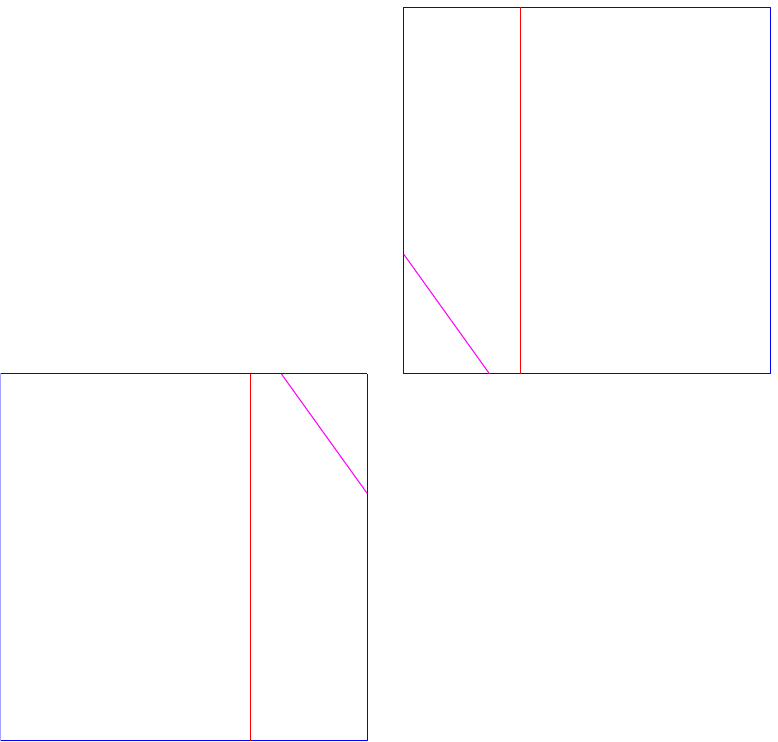}
$$
\caption{\label{newfig} \small $X_1$, $X_2$ (blue squares) and cuts (red -- default, magenta -- smart) for $s=11$ \\
{[$n=2$, $\epsilon=0.01$, $\delta=0.1$, $S=14$, $r(11)=0.0092$]}.}
\end{figure}
\begin{table}
\begin{center}
\begin{tabular}{|c|c|c|c|c|c|c|}
\hline
cuts&$n=2$&$n=3$&$n=4$&$n=5$&$n=6$&$n=20$\\
\hline
\hline
default&0.32&0.32&0.32&0.32&0.32&0.32\\
\hline
smart&3.8e-2&7.6e-3&1.4e-3&2.2e-4&3.4e-5&$<$1.3e-18\\
\hline
\end{tabular}
\end{center}
\caption{\label{newtable}\small $n$-dimensional volumes of ``bad'' sets $Z_{qs}$ for $s=11$
{[$\epsilon=0.01$, $\delta=0.1$, $S=14$, $r(11)=0.0092$]}.}
\end{table}

\subsection{Dynamical hypothesis testing}
We are about to present a particular application of the approach developed in section \ref{multiple}, aimed at {\sl dynamical hypotheses testing.}
\subsubsection{Situation} Let $\O_t=((\Omega_t,P_t),\{p_{\mu_t,t}(\cdot):\mu_t\in\M_t\},\F_t)$, $t=1,...,K$, be good observation
schemes, and consider their (nonstationary) $K$-factor direct product
{\small $$
\begin{array}{rcl}
\O^K&=&\big\{(\Omega^K,P^K)=(\Omega_1\times...\times\Omega_K,P_1\times..\times P_K),\\
&&\{p_{\mu^K}(\omega_1,...,\omega_K)=p_{\mu_1,1}(\omega_1)...p_{\mu_K,K}(\omega_K):\mu^K=[\mu_1;...;\mu_K]\in\M^K=\M_1\times...\times \M_K\},\\
&&\multicolumn{1}{r}{\F^K=\{f(\omega_1,...,\omega_K)=\sum_{t=1}^Kf_t(\omega_t):f_t\in\F_t,1\leq t\leq K\}\big\},}\\
\end{array}
$$}
so that our observation is
$$
\omega^K=(\omega_1,...,\omega_K)
$$
with independent of each other $\omega_t$, the densities of $\omega_t$ w.r.t. $P_t$ being $p_{\mu_t,t}(\omega_t)$, with $\mu_t\in \M_t$,  $t=1,...,K$. We treat $t$ as time and assume that the observations $\omega_t$ are consecutive portions of information we get. Now, we are given $N$ composite hypotheses on the distribution of $\omega^K$, $j$-th of them, $H_j$, stating that the density of $\omega^K$ w.r.t. $P^K$ is $p_{\mu^K}(\omega^K)$, with $\mu^K$ belonging to a given nonempty convex compact set $X_j\subset\M^K$; we refer to $\mu^K$ as to the {\sl signal} underlying the observation $\omega^K\sim p_{\mu^K}(\cdot)$.  The first of the hypotheses, $H_1$,  is treated as the null hypothesis (``no actual signal, just nuisance''), and the remaining hypotheses, let us call them {\sl signal} ones, are interpreted as stating the presence of a non-nuisance signal. We address here the simplest problem of deciding between the null hypothesis and the union of the signal ones, and focus on the situation when we are allowed for on-line decisions, meaning that at every time instant $t$, $1\leq t\leq K$, we decide, via the already observed part $\omega^t=[\omega_1;...;\omega_t]$ of observation $\omega^K$, whether or not the signal underlying our observations is a non-nuisance one; if we decide at time $t$ that this signal is not a nuisance (``non-nuisance conclusion at time $t$''), we terminate and ignore all subsequent observations $\omega_s$, otherwise we proceed to time instant $t+1$ (if $t<K$) or terminate (when $t=K$) with the conclusion that {\sl so far} the null hypothesis was valid. Our goal, informally, is to design a sequence of tests $T_t$ based on observations $\omega^t$, $1\leq t\leq K$, in such a way that the probability of {\sl false alarm} -- the worst, over $\mu^K\in X_1$, $p_{\mu^K}(\cdot)$-probability to make at some time $t\leq K$ the non-nuisance conclusion and thus to reject $H_1$ -- does not exceed a given risk level $\epsilon_{\rm f}$; under this restriction, we are interested to detect the presence of a non-nuisance signal, if any, as early as possible. Thus, what follows should be considered as a one more variation on {\sl change point detection}, the topic which is the subject of vast literature, see, e.g., \cite{Bass1988,Bass1993,Poor2009} and references therein.

\subsubsection{Dynamical test}
For a vector $z^K=[z_1;...;z_K]\in\M^K$, we denote by $z^t$ the vector $[z_1;...;z_t]$ \footnote{A less ambiguous notation would be something like $z^t[z^K]$; we prefer to use simplified notation, with the convention that {\sl whenever we are simultaneously speaking about vectors $\mu^t\in\M^t$ and $\mu^\tau\in\M^\tau$ with $\tau\geq t$} (or about $\omega^t$ and $\omega^\tau$, etc.), {\sl the first of these vectors is the initial fragment of the second one.}}. For $t\leq K$, let us set
$$
X^t_j=\{\mu^t:\mu^K\in X_j\},
$$
so that $X^t_j$ are nonempty closed convex subsets in $\M^t=\M_1\times...\times\M_t$; we denote by $H_j^t$ the hypothesis stating that the observation $\omega^t$ stems from a $t$-long signal $\mu^t\in X^t_j$:
$$\omega^t\sim p_{\mu^t}(\omega^t)=\prod_{s=1}^tp_{\mu_s,s}(\omega_s).$$
Note that $X^t_1$ can intersect with some $X^t_j$, $j>1$ even when $X_1=X^K_1$ is pretty far from $X_j=X^K_j$, implying that while the null hypothesis $H_1$ can be easily distinguished from  a particular signal hypothesis $H_j$ via $K$ observations $\omega^K$, it is impossible to distinguish reliably between these hypotheses via a given number $t<K$ observations.\par
 The approach we are about to develop is aimed at designing tests $T_t$, $t=1,2,...,K$, which respect a given upper bound $\epsilon_{\rm f}$ on the probability of false alarm and at the same time try to
recognize the presence of a non-nuisance signal $\mu^t\in X^t_j$ whenever $X^t_j$ is ``far enough'' from $X^t_1$. The seemingly simplest implementation of this desire results in the following construction:
\begin{quote}Given risk levels $\epsilon_{\rm f}\in(0,1)$, $\epsilon_{\rm m}\in(0,1)$ (upper bounds on probabilities of false alarm and miss, respectively), we act as follows.\par
{\bf A.} We consider $N$ hypotheses $H^t_j$, $1\leq j\leq N$, along with the pairwise detectors $\phi_{jj'}^t(\omega^t)=-\phi_{j'j}^t(\omega^t)$, $1\leq j,j'\leq N$, yielded by Theorem \ref{the1}. Thus,
\begin{equation}\label{asusual}
\begin{array}{l}
\forall (1\leq j,j'\leq N,\,\mu^t\in X_j^t):\\
\int_{\Omega^t} \exp\{-\phi_{jj'}^t(\omega^t)\}p_{\mu^t}(\omega^t)P^t(d\omega^t)\leq\epsilon_{jj't}:=\max\limits_{\mu^t\in X^t_j,\nu^t\in X^t_{j'}}
\int_{\Omega^t}\sqrt{p_{\mu^t}(\omega^t)p_{\nu^t}(\omega^t)}P^t(d\omega^t),\\
\multicolumn{1}{c}{\phi_{jj'}^t(\cdot)\equiv -\phi_{j'j}(\cdot),\,\,\epsilon_{jj't}=\epsilon_{j'jt}\in(0,1]\
\left[\Omega^t=\Omega_1,...,\Omega_t,P^t=P_1\times...\times P_t\right]}\\
\end{array}
\end{equation}
\par{\bf B.} Let us eliminate from the vector  $\varepsilon=[\epsilon_{11t};\epsilon_{12t};...;\epsilon_{1Nt}]$ a number $L$ of the largest entries in such a way that the $\|\cdot\|_2$-norm of the resulting vector $\widehat{\varepsilon}=[\widehat{\epsilon}_1;...;\widehat{\epsilon}_N]$ is $\leq\widehat{\epsilon}:=\sqrt{\epsilon_{\rm f}\epsilon_{\rm m}/K}$; under the latter requirement, we select $L$ to be as small as possible.  Let $J_t$ be the set of indexes $j$ of the removed entries; note that $1\in J_t$ due to $\epsilon_{11t}=1>\widehat{\epsilon}$. The test $T_t$, as applied to observation $\omega^t$, is as follows\footnote{Essentially, the test below is nothing but the test from section \ref{multiple} as applied to the hypotheses $H^t_j$, $1\leq j\leq N$, and properly defined closeness $\C_t$.}:
\begin{quote} Given $\omega^t$, we compute the quantities
$$
\begin{array}{rcl}
\overline{\phi}_{j1}^t(\omega^t)&=&\phi_{j1}^t(\omega^t)+\ln\left({\epsilon_{j1t}\over\epsilon_{\rm m}}\right),\,1\leq j\leq N.
\end{array}
$$
If there exists $j\not\in J_t$ such that $\overline{\phi}_{j1}^t(\omega^t)\geq0$, we make the no-nuisance conclusion and terminate the inference process, otherwise claim that so far the null hypothesis is accepted, and pass to the next time instant (when $t<K$), or terminate (when $t=K$).
\end{quote}
 \end{quote}
By reasons which will become clear in a while, in the sequel we refer to the signal hypotheses $H_j$ (and $H_j^t$) with $j\not\in J_t$ as to signal hypotheses {\sl recognizable} at time $t$.\par
What we can say about the resulting dynamical inference procedure $\T$ is as follows.
\begin{proposition}\label{propdynamic} Let the signal $\mu^K$ underlying observation $\omega^K$ belong, for some $j_*$, to $X_{j_*}$, and let $p(\cdot)=p_{\mu^K}(\cdot)$ be the corresponding probability density of $\omega^K$. Then
\begin{itemize}
\item[{\rm (i)}]  When $j_*=1$ {\rm (i.e., the null hypothesis holds true)}, $p(\cdot)$-probability of false alarm, i.e. the probability  for $\T$ to terminate somewhere on the time horizon $1,...,K$ with the non-nuisance conclusion, is at most $\epsilon_{\rm f}$;
\item[{\rm (ii)}] Whenever $t$, $1\leq t\leq K$, is such that  $j_*\not\in J_t$ {\rm (i.e., $H_{j_*}$ is a signal hypothesis recognizable at time $t$)}, $p(\cdot)$-probability for $\T$ to terminate with the non-nuisance conclusion, and to do so  not later than at time $t$, is at least $1-\epsilon_{\rm m}$.
\end{itemize}
\end{proposition}
{\bf Proof.} (i): Let $j_*=1$, and let us fix $t\leq K$. If $\T$ terminates with a non-nuisance conclusion at time instant $t$, then
$\overline{\phi}_{j1}^t(\omega^t)\geq0$ for some $j=j(\omega^t)\not\in J_t$, implying that $\phi_{j1}^t(\omega^t)\geq- \ln(\epsilon_{j1t}/\epsilon_{\rm m})$, whence $\phi_{1j}^t(\omega^t)\leq \ln(\epsilon_{ijt}/\epsilon_{\rm m})$. By the union bound, $p(\cdot)$-probability $\pi^t$ of this event is at most $\sum_{j\not\in J_t}\pi_j$, where $\pi_j$ is the $p(\cdot)$-probability of the event $\phi_{1j}^t(\omega^t)\leq \ln(\epsilon_{j1t}/\epsilon_{\rm m})$. Since the null hypothesis is true, the density, induced by $p(\cdot)$, of $\omega^t$ obeys $H^t_1$, so that (\ref{asusual}) implies that $\pi_j\leq {\epsilon_{j1t}^2\over\epsilon_{\rm m}}$ and thus
$$\pi^t\leq \sum_{j\not\in J_t} {\epsilon_{j1t}^2\over\epsilon_{\rm m}}=\sum_{j\not\in J_t} {\epsilon_{1jt}^2\over\epsilon_{\rm m}} \leq \epsilon_{\rm f}/K,$$
where the concluding inequality is due to the definition of $J_t$. Thus, the $p(\cdot)$-probability for $\T$ to terminate at time $t$ with a non-nuisance conclusion is at most $\epsilon_{\rm f}/K$, which combines with the union bound to imply (i).
 \par
 (ii) Now let $p(\cdot)\in X_{j_*}$ with $j_*>1$ such that $H_{j_*}$ becomes recognizable at some time instant $t$, so that  $j_*\not\in J_t$. The density of $\omega^t$ induced by $p(\cdot)$ obeys $H^t_{j_*}$, whence by (\ref{asusual}) the $p(\cdot)$-probability of the event $\overline{\phi}_{j_*1}^t(\omega^t)<0$, or, which is the same, the event $\phi_{j_*1}^t(\omega^t) <-\ln(\epsilon_{j_*1t}/\epsilon_{\rm m})$, is at most $\epsilon_{\rm m}$. When the just defined event does {\sl not} take place (``good case''), the test $T_t$, by construction and due to $j_*\not\in J_t$, terminates with the no-nuisance conclusion, implying that in the good case, $\T$ terminates at time $t$ or earlier  with the no-nuisance conclusion. As we have seen, the $p(\cdot)$-probability of the complement to the good case is $\leq \epsilon_{\rm m}$. (ii) is proved. \qed
\par
Note that (i-ii) say nothing on what may happen at time $t$ when the true hypothesis is a signal hypothesis $H^t_j$ with $j\in J_t$ (i.e., a signal hypothesis which is not recognizable at time $t$). In principle, this is unavoidable, since on time horizon $1,...,t$ hypotheses of this type can be close to, or even undistinguishable from, the null hypothesis.
\subsubsection{Numerical illustration}
It seems to be impossible to say once for ever how close our dynamical inference scheme is to ``limits of performance'' of schemes of this type. To get an impression of what can we expect here, we are about to report on a simple numerical experiment which we find instructive.
\paragraph{The setup} of our experiment is as follows. The good observation schemes $\O_t$, $t=1,...,K$, are identical to each other standard scalar Gaussian schemes, so that
$$
\omega_t=\mu_t+\xi_t
$$
with $\mu_t\in\M_t\equiv \bbr$ and independent across $t=1,2,...,K$ observation noises $\xi_t\sim\cN(0,1)$.  As we have seen in section \ref{sect:Gauss}, this corresponds to
$$
\begin{array}{c}
\phi_{jj'}^t(\omega^t)={1\over 2}[a_{jj'}^t-b_{jj'}^t]^T\omega^t + {1\over 4}[\|b_{jj'}^t\|_2^2-\|a_{jj'}^t\|_2^2],\,\,
\epsilon_{j'jt}=\exp\{-\|a_{jj'}^t-b_{jj'}^t\|_2^2/8\},\\
\end{array}
$$
where $(a_{jj'}^t\in X^t_j,b_{jj'}^T\in X^t_{j'})$ form an optimal solution to the optimization problem
$$
\min_{a\in X^t_j,b\in X^t_{j'}}\|a-b\|_2.
$$
We assume w.l.o.g. that $a_{jj'}^t=b_{j'j}^t$, $1\leq j,j'\leq N$. \par
Further, the sets $X_j$ in our experiment are specified as follows. We are given a discrete time linear time invariant system with finite, of duration $T$, scalar impulse response and observe on time horizon $1,...,K$ the output $\mu^K=Ax$ of this system, corrupted by the standard Gaussian noise; here $x\in\bbr^{T-1+K}$ is the vector comprised of inputs to the system on the time horizon $D=\{-T+2,-T+3,...,K\}$, and $A$ is $K\times (T-1+K)$ matrix readily given by the impulse response. In our experiment, $X_1$ is just the singleton $\{0\}$, that is, the null hypothesis states that the input to the system on the time horizon $D$ is identically zero. The remaining signal sets have two indexes $\tau\in D$ and $\nu$, and correspond to inputs starting at time $\tau$ and known up to their amplitude, specifically,
 $$
 X_{\tau\nu}=\{\mu^K=rAx^\tau:r\geq r_\nu\},
 $$
 where $x^\tau$ is the ``step at time $\tau$'':
 $$
 x^\tau_s=\left\{\begin{array}{ll} 0,&-T+1\leq s<\tau\\
 1,&\tau\leq s\leq K\\
 \end{array}\right.
 $$
 and
 $$
r_\nu=0.1(1.2)^{\nu-1},\,\nu=1,2,...,44
$$
 (44 is the largest value of $\nu$ resulting in $r_\nu<256$). Finally, we used the same impulse response as in sections \ref{sec:numerics}, \ref{signalident}, of duration $T=32$, and our time horizon was $K=32$. With this setup, the number of signal hypotheses was $N=(2T-1)\cdot44=2772$. The tolerances $\epsilon_{\rm f}$, $\epsilon_{\rm m}$ underlying our dynamical test were set to $0.01$.
 \paragraph{The primary goal} of the experiment was to understand how conservative under circumstances is our dynamical inference procedure $\T$. Here is how we measured the conservatism. Let us fix a time instant $t$, $1\leq t\leq K$.
 \par
 All (output) signals $\mu^K$ from the signal set $X_{\tau\nu}$ are zero prior to a common time instant $t(\tau)$ depending solely on $\tau$. When
 $t(\tau)\leq t$, observations $\omega^t$ do not allow to distinguish between the null hypothesis and the signal hypothesis $H_{\tau\nu}$ stating that $\mu^K\in X_{\tau\nu}$. Besides this, depending on  $\tau$, $\nu$, some signals $\mu^t\in X_{\tau\nu}^t$, while not identically zero on the time horizon $1,...,t$, are too small to be distinguishable from $\mu^t=0$ postulated by the null hypothesis $H_1^t$. Let us fix a reasonable reliability level, say, 0.9. In order for the smallest amplitude signal in $X_{\tau\nu}$, that is, $r_\nu Ax^\tau$, to be distinguished 0.9-reliably from the zero signal via noisy observations from time 1 to time $r$ we should have
 $$
\sqrt{{\sum}_{s=1}^r[r_\nu Ax^\tau]_s^2}\geq 2 d,
$$
where $d\approx 1.28$ is the 0.9-quantile of the standard Gaussian distribution. Given $\tau$ and $\nu$, let $r=r(\tau,\nu)$ be the smallest $r\leq K$ for which the  above relation holds true (if no such $r$ exists, we set $r(\tau,\nu)=+\infty$). We refer to $r(\tau,\nu)$ as to the time instant when the signal set $X_{\tau,\nu}$ and the hypothesis $H_{\tau\nu}$ become {\sl visible}. When making decisions on the basis of $t$ observations $\omega^t$, we have no hope to distinguish somehow reliably a signal hypothesis $H^t_{\tau\nu}$ with $r(\tau,\nu)>t$ (i.e., a signal hypothesis still invisible at time $t$) from the null hypothesis, even if instead of multiple signal hypotheses we were dealing with the single hypothesis $H_{\tau\nu}$. Let $N(t)$ be the total number of hypotheses $H_{\tau\nu}$ which become visible at time $t$ or earlier. One way to quantify the conservatism of our dynamical decision making $\T$ is to compare with $N(t)$, for each $t=1,2,...,K$, the number $N^+(t)$ of signal hypotheses which are recognizable at time $t$ or earlier\footnote{recall that by Proposition
\ref{propdynamic}.ii, if a signal hypothesis is recognizable at some time $t'$, the probability, the hypothesis being valid,  for $\T$ to terminate with the ``no nuisance'' conclusion at time $t'$ or earlier is at least $1-\epsilon$.}.  The ratios $N^+(t)/N(t)$, $t=1,2,...,K$ can be thought of as some measure of conservatism of $\T$ -- the smaller these ratios, the larger the conservatism. \par Another, closely related way to quantify the conservatism of $\T$ is as follows. Along with $r(\tau,\nu)$ -- the first time instant where $H_{\tau\nu}$ becomes visible or, which is the same, ``in the nature'' there are chances to distinguish somehow reliably between the signal hypothesis $H_{\tau\nu}$ and the null hypothesis -- let us define $r^+(\tau,\nu)$ as the first time instant $t'\leq K$ where $H_{\tau\nu}$ becomes recognizable (and thus $\T$ terminates, with probability $1-\epsilon$, with the ``no nuisance'' conclusion, at time $t'$ or earlier, provided $H_{\tau\nu}$ holds true); here again we set $r^+(\tau,\nu)=+\infty$ when no required $t'$ exists. Given $t$, we can look at all signal hypotheses
$H_{\tau\nu}$ which are both visible and recognizable at time $t$ (i.e., those with $r(\tau,\nu)\leq t$ and $r^+(\tau,\nu)\leq t$). For every one of these  hypotheses $H_{\tau\nu}$, the quantity $r^+(\tau,\nu)/r(\tau,\nu)-1$ can be thought of as the (upper bound on the) {\sl relative delay} in our acceptance of the ``no nuisance'' conclusion, the true hypothesis being $H_{\tau\nu}$, as compared to arriving at the same conclusion by ``the best existing in the nature'' decision procedure. We can now think about the empirical characteristics, like the mean $\hbox{\rm mean}(t)$ and the median $\hbox{median}(t)$, of the array of relative delays as seen at time $t$ as of additional indicators of how conservative $\T$ is (the less the indicator, the less the conservatism).
\paragraph{The results} of our experiment are presented in  Table \ref{tablecompare}. To the best of our understanding, these results suggest that our simple dynamical hypotheses testing $\T$ is not too conservative.\par
 Along with quantifying the conservatism of $\T$, we carried out intensive numerical simulation aimed at quantifying empirical reliability of this procedure. Recall that when building $\T$, the tolerances $\epsilon_{\rm f}=\epsilon_{\rm m}=0.01$ were used, meaning that theoretically speaking, we could expect about 1\% of false alarms (non-nuisance conclusions under the null hypothesis) and about 1\% of ``bad misses'' (a signal obeying signal hypothesis recognizable at time $t$ was qualified by $T_t$ as nuisance). In fact, in about 50,000 simulations the fraction of false alarms  was well below $0.001$, and  not a  {\sl single} bad miss was observed. Thus, empirically speaking, $\T$ is essentially more reliable than our theoretical risk bounds say\footnote{See section \ref{sect:Gauss}  for  the comparison of the actual risk of the simple Gaussian case test yielded by Theorem \ref{the1}
and the upper bound on this risk as given by this Theorem.}, and we can use this phenomenon, at a heuristic level, to reduce the conservatism of the dynamical test. Specifically,
when passing to dynamical test $\T'$ associated with the tolerances as large (and by themselves meaningless) as $\epsilon_{\rm f}=\epsilon_{\rm m}=0.8$, the empirical probabilities of false alarm and of bad miss, as observed in an extensive simulation study, remain reasonably small (less than 1.1\% for false alarm and less than 0.5\% for bad miss), but the conservatism reduces significantly: the means mean$(t)$ and the medians median$(t)$ of the relative delays, which were about 43\% (means) and $7\%$ (medians) for $\T$, drop to just 10\% (means) and 0\% (medians) for $\T'$, and the ratios $N^+(t)/N(t)$, which for $\T$ varied, depending on $t$, from 0.85 to 0.95, are ``lifted'' to the range 0.94 -- 0.98.
\begin{table}
{\tiny
$$
\begin{array}{||c||c|c|c|c|c|c||}
\cline{2-7}
\multicolumn{1}{c||}{}
&t=1&t=2&t=3&t=4&t=5&t=6\\
\hline
N_{\hbox{\tiny tot}}(t)&1408&1452&1496&1540&1584&1628\\
\hline
N^+(t)/N(t)&850/1016\approx0.84&932/2203\approx0.85&992/1165\approx0.85&1047/1216\approx0.86&1098/1264\approx0.87&1145/1309\approx0.88\\
\hline
\hbox{mean}(t)&0.00&0.09&0.17&0.26&0.32&0.35\\
\hline
\hbox{median}(t)&0.00&0.00&0.00&0.00&0.00&0.00\\
\hline
\multicolumn{7}{c}{}\\
\cline{2-7}
\multicolumn{1}{c||}{}&t=7&t=8&t=9&t=10&t=11&t=12\\
\hline
N_{\hbox{\tiny tot}}(t)&1672&1716&1760&1804&1848&1892\\
\hline
N^+(t)/N(t)&1193/1355\approx0.88&1240/1400\approx0.89&1284/1445\approx0.89&1331/1490\approx0.89&1378/1534\approx0.90&1423/1578\approx0.90\\
\hline
\hbox{mean}(t)&0.39&0.41&0.42&0.43&0.44&0.45\\
\hline
\hbox{median}(t)&0.00&0.00&0.00&0.00&0.00&0.00\\
\hline
\multicolumn{7}{c}{}\\
\cline{2-7}
\multicolumn{1}{c||}{}&t=13&t=14&t=15&t=16&t=17&t=18\\
\hline
N_{\hbox{\tiny tot}}&1936&1980&2024&2068&2112&2156\\
\hline
N^+(t)/N(t)&1467/1622\approx0.90&1513/1666\approx0.91&1560/1710\approx0.91&1605/1754\approx0.92&1650/1798\approx0.92&1695/1842\approx0.92\\
\hline
\hbox{mean}(t)&0.45&0.46&0.46&0.46&0.46&0.46\\
\hline
\hbox{median}(t)&0.00&0.00&0.00&0.00&0.00&0.00\\
\hline
\multicolumn{7}{c}{}\\
\cline{2-7}
\multicolumn{1}{c||}{}&t=19&t=20&t=21&t=22&t=23&t=24\\
\hline
N_{\hbox{\tiny tot}}(t)&2200&2244&2288&2332&2376&2420\\
\hline
N^+(t)/N(t)&1740/1886\approx0.92&1784/1930\approx0.92&1828/1974\approx0.93&1873/2018\approx0.93&1917/2063\approx0.93&1961/2106\approx0.93\\
\hline
\hbox{mean}(t)&0.45&0.45&0.44&0.44&0.43&0.43\\
\hline
\hbox{median}(t)&0.00&0.05&0.05&0.06&0.06&0.06\\
\hline
\multicolumn{7}{c}{}\\
\cline{2-7}
\multicolumn{1}{c||}{}&t=25&t=26&t=27&t=28&t=29&t=30\\
\hline
N_{\hbox{\tiny tot}}(t)&2464&2508&2552&2596&2640&2684\\
\hline
N^+(t)/N(t)&2005/2150\approx0.93&2049/2194\approx0.93&2093/2238\approx0.94&2137/2282\approx0.94&2181/2326\approx0.94&2225/2370\approx0.94\\
\hline
\hbox{mean}(t)&0.43&0.42&0.42&0.41&0.41&0.40\\
\hline
\hbox{median}(t)&0.06&0.06&0.07&0.07&0.07&0.07\\
\hline
\multicolumn{7}{c}{}\\
\cline{2-3}
\multicolumn{1}{c||}{}&t=31&t=32&\multicolumn{4}{c}{}\\
\cline{1-3}
N_{\hbox{\tiny tot}}(t)&2728&2772&\multicolumn{4}{c}{}\\
\cline{1-3}
N^+(t)/N(t)&2269/2414\approx0.94&2313/2458\approx0.94&\multicolumn{4}{c}{}\\
\cline{1-3}
\hbox{mean}(t)&0.39&0.39&\multicolumn{4}{c}{}\\
\cline{1-3}
\hbox{median}(t)&0.07&0.07&\multicolumn{4}{c}{}\\
\cline{1-3}
\end{array}
$$
}
\caption{\label{tablecompare} Dynamical hypotheses testing.
\footnotesize
$N_{\hbox{\tiny tot}}(t)$: $\#$ of signal hypotheses $H_{\tau\nu}$ with $\tau\leq t$;
$N(t)$: $\#$ of signal hypotheses $H_{\tau\nu}$ visible at time $t$;
$N^+(t)$: $\#$ of signal hypotheses $H_{\tau\nu}$ recognizable at time $t$;
mean$(t)$, median$(t)$: mean and median of relative delay.}
\end{table}}

\section{Case studies}\label{sec:cases}
\subsection{Hypotheses testing in PET model}
To illustrate applications of the simple test developed in section \ref{sect:Poiss} we discuss here a toy testing problem in the {\em Positron Emission Tomography (PET) model}.
\par
A  model of PET which is  accurate enough for medical purposes is as follows. The patient is injected a radioactive tracer and is placed inside a cylinder with the inner surface split into detector cells. Every tracer disintegration act gives rise to two $\gamma$-quants flying in opposite directions along a  randomly oriented line (Line of Response, LOR) passing through the disintegration point. Unless the LOR makes too small angle with the cylinder's axis, the $\gamma$-quants activate (nearly) simultaneously a pair of detector cells; this event (``coincidence'') is registered, and the data acquired in a PET study is the list of the detector pairs in which the coincidences occurred. The goal of the study is to infer about the density of the tracer on the basis of these observations.
\par After appropriate discretization of the field of view into small cells, disintegration acts in a particular cell form a Poisson processes with intensity proportional to the density of the tracer in the cell. The entries of the observations vector $\omega$ are indexed by {\sl bins} $i$ -- pairs of detectors, $\omega_i$ being the number of coincidences registered during the study by bin $i$.  Mathematically, $\omega_i$, $i=1,...,m$, are the realizations of independent across $i$'s Poisson random variables with parameters $\mu_i=(tP\lambda)_i$, where $t$ is the observation time, $\lambda$ is the vector of intensities of disintegration in the cells of the field of view, and the entries $P_{ij}$ in the matrix $P$ are the probabilities for a LOR originating in cell $j$ to be registered by bin $i$; this matrix is readily given by the geometry of PET's device. We observe that PET model meets the specifications of what we call Poisson observation scheme.
\par
Let $\M$ be the image, under the linear mapping $\lambda\mapsto {t}P\lambda$, of the set $\Lambda=\Lambda_{L,R}$ of non-vanishing on $\bbr^n$ densities $\lambda$ satisfying some regularity restrictions, specifically, such that the uniform norm of discrete Laplacian of $\lambda$ is upper-bounded by $L$,  and the average of $\lambda$, over all pixels, is upper-bounded by $R$, i.e.
\[
\Lambda_{L,R}=\left\{
\begin{array}{c}
\lambda\in \bbr^n:\;\lambda\ge 0,\;n^{-1}\sum_{j=1}^n \lambda_{j}\leq R,\\
 \four|4\lambda_{j(k,\ell)}-\lambda_{j(k-1,\ell)}-\lambda_{j(k,\ell-1)}-\lambda_{j(k+1,\ell)}-\lambda_{j(k,\ell+1)}|\leq L,\;\;1\leq j\leq n
 \end{array}
 \right\},
\]
$(k,\ell)$ being the coordinates of the cell $j$ in the field of view  (by convention, $\lambda_{j(k,\ell)}=0$ when the cell $(k,\ell)$ is not in the field of view).
Our goal is to distinguish two hypotheses, $H_1$ and $H_2$, about $\lambda$:
$$
H_1:\;\lambda\in \Lambda_1=\{\lambda\in \Lambda:\;g(\lambda)\leq \alpha\},\;\;H_1:\;\lambda\in \Lambda_2=\{\lambda\in \Lambda:\;g(\lambda)\geq \alpha+\rho\},\eqno{(\P_{g,\alpha}[\rho])}
$$
$g(\lambda)=g^T\lambda$ being a given linear functional of $\lambda$. From now on we assume that $g\notin \Ker(P)$ and $\rho>0$, thus the described setting corresponds to the Poisson case of the hypotheses testing problem of section \ref{sect:Poiss}, $X=tP\Lambda_1$ and $Y=tP\Lambda_2$ being two nonintersecting convex sets of observation intensities.
Let us fix the value $\epsilon\in (0,1)$, and consider the optimization problem
\be
t_*=\min_t\max_{\lambda,\lambda'}\left\{
t:\;\begin{array}{l}
-{t\over 2} \sum_{i=1}^{m}\left[\sqrt{[P\lambda]_i}-\sqrt{[P \lambda']_i}\right]^2
\geq \ln \epsilon,\\
\lambda,\lambda'\in \Lambda,\;g(\lambda)\leq \alpha,\;g(\lambda')\geq \alpha+\rho.
\end{array}\right\}
\ee{eq:P2}
Suppose that the problem parameters are such that both hypotheses in $(\P_{g,\alpha}[\rho])$ are not empty.  It can be easily seen that in this case problem \rf{eq:P2} is solvable and its optimal value $t_*$ is positive $0$. Let $[\lambda_*;\lambda'_*]$ be the $[\lambda;\lambda']$-component of an optimal solution to \rf{eq:P2}, consider the test $T_*$ associated with the detector
 \be
 \phi_*(\omega)=\half \sum_{i=1}^m\ln\left[{[P\lambda_*]_i\over [P\lambda'_*]_i}\right] \omega_i-\half\sum_{i=1}^m[P\lambda_*-P\lambda'_*]_i.
 \ee{eq:phiP}
By applying  Theorem \ref{the1} in the Poisson case (cf. \rf{Poissoncase}) we conclude that
the risk of the test $T_*$ associated with detector $\phi_*$, when applied to the problem testing problem $(\P_{g,\alpha}[\rho])$ is bounded with $\epsilon$, as soon as the observation time $t\geq t_*$.
\par
\begin{figure}
$$
\begin{array}{cc}
\epsfxsize=220pt\epsfysize=220pt\epsffile{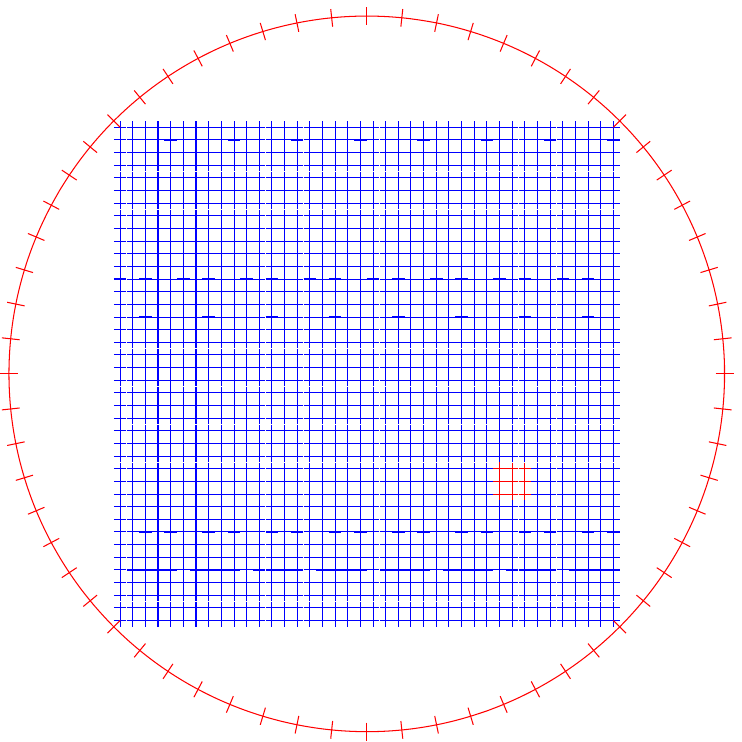}&\epsfxsize=220pt\epsfysize=220pt\epsffile{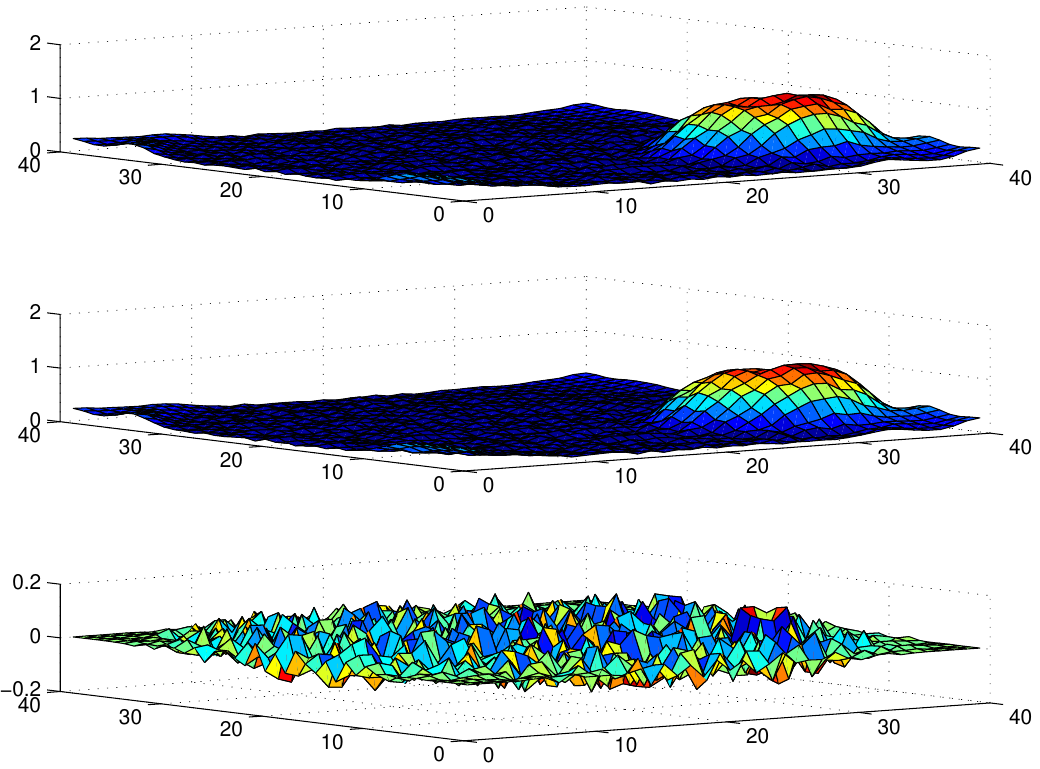}
\end{array}
$$
\caption{\label{figPET}
Toy PET experiment. Left: $40\times 40$ field of view with $3\times 3$ ``suspicious spot'' (in red) and the ring of 64 detector cells. Right: the hardest to distinguish tracer densities $\lambda_*$ (top) and $\lambda'_*$ (middle), and the difference of these densities (bottom).}
\end{figure}

In {the} numerical experiment we are about to describe we simulate a 2D PET device with square field of view split into $40\times40 $ pixels (i.e., dimension of $\lambda$ was $n=1600$). The detector cells are represented by $k=64$ equal arcs of the circle circumscribing the field of view, resulting in the observation space (pairs of detectors which may be activated during the experiment) of dimension $m=1536$.
We choose $g(\cdot)$ to be the density average over a specific $3\times 3$ ``suspicious spot'' (see the left plot on figure \ref{figPET}), and values of $\alpha=1.0$ and $\rho=0.1$, so that under $H_1$ the average of the density $\lambda$ of the tracer on the spot is upper-bounded by 1, while under $H_2$ this average is at least 1.1. The regularity parameters of the density class $\Lambda_{L,R}$ were set to $L=0.05$ and $R=1$, the observation time $t^*$ and parameters of the detector $\phi_*$ were selected according to \rf{eq:P2} and \rf{eq:phiP} with $\epsilon=0.01$.
 \par
 On the right plot on figure \ref{figPET} we present the result of computation of the hardest to distinguish densities $\lambda_*\in \Lambda_1$ and $\lambda'_*\in \Lambda_2$.
We have also measured the actual performance of our test by simulating $2000$ PET studies with varying from study to study density of the tracer. In the first 1000 of our simulations the true density was selected to obey $H_1$, and in the remaining 1000 simulations -- to obey $H_2$, and we did our best to select the  densities which make decision difficult.
 In the reported experiment the empirical probabilities to reject the true hypothesis were $0.005$ when the true hypothesis was $H_1$, and $0.008$  when the true hypothesis was $H_2$.

\subsection{Event detection in sensor networks}\label{sec:applex}
\subsubsection{Problem description} Suppose that $m$ sensors are deployed on the domain $G\subseteq \bbr^d$. The signals are real-valued functions  $x:\,\Gamma\to \bbr^n$ on a grid $\Gamma=(\gamma_i)_{i=1,...,n}\subset G$, and  the observation $\omega_j$ delivered by  $j$th sensor, $j=1,...,m$, is a linear form of the signal, contaminated with random noise. So we have at our disposal {an} observation $
\omega\sim P_\mu$ --
a random vector in $\bbr^m$  with the distribution parameterized by $\mu\in \bbr^m$, where $\mu=Ax$ {and} $A\in \bbr^{m\times n}$ is a known matrix of sensor responses ($j$th row of $A$ is the response of the $j$th sensor). Further, we assume that the signal $x$ can be decomposed into
$x=s+v$, where $v\in {\cal V}$ is a background (nuisance) signal, $\cV$ is a known convex and compact set in $\bbr^n$.
We assume that {at most} one event can take place during the observation period, and an event {occurring} at a node $\gamma_i$ of the grid produces the signal $s=re[i]\in \bbr^n$ on the grid of known signature $e[i]$ {with unknown real factor $r$}.
\mypar
We want to decide whether an event occurred during the observation period, i.e. to test the null hypothesis that no event happened against the alternative  that exactly one event took place. To make
{a} consistent decision possible we need the alternative to be separated from the {null} hypothesis,
so we require, first, that $Ae[i]\neq0$ for all $i$, and, second, that under the alternative, when an event occurs at a node $\gamma_i\in \Gamma$, we have $s=re[i]$ with $|r|\geq\rho_i$ with some given $\rho_i>0$. Thus we come to the testing problem as follows:
\be
({\cal D}_\rho)\;\;&\hbox{\begin{tabular}{l}Given $\rho=[\rho_1;...;\rho_n]>0$, test the hypothesis $H_0:\;s=0$ against the\\
alternative $H_1(\rho):$ $s=re[i]$ for some $i\in\{1,...,n\}$ and $r$ with $|r|\geq\rho_i$.\\
 \end{tabular}}
\ee{ed1}
Our goal is,
 given an $\epsilon\in(0,1)$, to construct a test with risk $\leq\epsilon$  for  as wide as possible (i.e., with as small $\rho$ as possible) alternative $H_1(\rho)$.
 \mypar
The problem of multi-sensor detection have recently received much attention in the signal processing and statistical literature
(see e.g., \cite{TartVeer2004,TartVeer2008} and references therein). Furthermore, a number of classical detection problems, extensively studied in statistical literature, such as detecting
jumps in derivatives of a function and cusp detection \cite{AntoniaGuib2002,GJTZ2008,goldenshluger2008change,kuznetsov1976,Mull1999,neumann1997optimal,Wang1995,yin1988detection},  detecting a nontrivial signal on input of a dynamical system \cite{Gust2000}, or parameter change detection \cite{Bass1988} can be posed as $(\cD_\rho)$.
\par
Our current objective is to apply the general approach described in section \ref{singleobservationcase}
 to the problem
$({\cal D}_\rho)$.
Note that, in terms of the parameter $\mu$ underlying the distribution of the observation $\omega$, the hypothesis $H_0$ corresponds to $\mu \in X:=A\cV$, a convex compact set, while the alternative $H_1$ is represented by the union $Y=\bigcup\limits_{i=1}^n Y_{i}$ of the sets $Y_{i}=\{Are[i]+\nu,\;\nu\in \cV, \,|r|\ge \rho_i\}$.
To comply with assumptions of section \ref{sec:mainres} we bound the sets $Y_i$ by imposing an upper bound on the amplitude $r$ of the useful signal: from now on we assume that $\rho_i\leq |r|\leq R$ in the definition of $({\cal D}_\rho)$.\footnote{Imposing a finite upper bound $R$ on $|r|$ is a minor  (and non-restrictive, as far as applications are concerned) modification of the problem stated in the introduction; the purely technical reason for this modification is our desire to work with compact sets of parameters. It should be stressed that $R$ does not affect the performance bounds to follow.}
\mypar
Given a test $\phi(\cdot)$ and $\epsilon>0$, we call a collection $\rho=[\rho_1;...;\rho_n]$ of positive reals {\sl an $\epsilon$-rate profile of the test $\phi$} if whenever the signal $s$ underlying our observation is $re[i]$ for some $i$ and $r$ with $\rho_i\leq |r|\leq R$, the hypothesis $H_0$ will be rejected by the test with probability $\geq1-\epsilon$, whatever be the nuisance $v \in\cV$, and whenever $s=0$, the probability for the test to reject $H_0$ is $\leq\epsilon$, whatever be the nuisance $v \in\cV$. Our goal is to design a test with $\epsilon$-rate profile ``nearly best possible'' in the sense of the following definition:
\begin{quotation}
\noindent Let $\kappa\geq 1$. A test $T$ with risk $\epsilon$ in the problem $({\cal D}_\rho)$
is said to be {\em $\kappa$--rate optimal}, if there is no test with the risk $\epsilon$ in the problem $({\cal D}_{\underline{\rho}})$ with $\underline{\rho}<\kappa^{-1} \rho$ (inequalities between vectors
are understood componentwise).
\end{quotation}

\subsubsection{Poisson case} Let the sensing matrix $A$ be nonnegative and without zero rows, let the signal $x$ be nonnegative, and let the entries $\omega_i$ in our observation be independent and obeying Poisson distribution with the intensities $\mu:=[\mu_1;...;\mu_m]=Ax$.
In this case the null hypothesis is that the signal is a pure nuisance:\\
$$
H_0:\; \mu\in X=\{\mu=Av ,\;v \in \cV\},
$$
where $\cV$ is the nuisance set assumed to be a nonempty compact convex set belonging to the interior of the nonnegative orthant.
The alternative $H_1(\rho)$ is the union over $i=1,...,n$ of the hypotheses
$$
\begin{array}{l}
H^{i}(\rho_i):\;\mu\in Y(\rho_i)=\{rAe[i]+Av ,\;v \in \cV,\;\rho_i\leq r\leq R\},
\end{array}
$$
where $e[i]\geq0$, $1\leq i\leq n$, satisfy $Ae[i]\neq0$.
For $1\leq i\leq n$, let us set (cf. section \ref{sect:Poiss})
$$
\rho^P_i(\epsilon)=\max_{\rho,r,u,v}\left\{\rho:
\begin{array}{l}
{1\over 2}{\sum}_{\ell=1}^m\left[\sqrt{[Au]_\ell}-\sqrt{[A(re[i]+v)]_\ell}\right]^2\leq\ln(\sqrt{n}/\epsilon)\\
u\in\cV,v\in\cV,r\geq\rho\\
\end{array}\right\},\eqno{(P^i_\epsilon)}
$$
\begin{equation}\label{phiP}
\phi_i(\omega)={\sum}_{\ell=1}^m\ln(\sqrt{[Au^i]_\ell/[A(r^ie[i]+v^i)]_\ell})\omega_\ell -{1\over 2}{\sum}_{\ell=1}^m[A(u^i-r^ie[i]-v^i)]_\ell,
\end{equation}
where $r^i$, $u^i$, $v^i$ are the $r,u,v$-components of an  optimal solution to  $(P^i_\epsilon)$ (of course, in fact $r^i=\rho^P_i(\epsilon)$). Finally, let
\[
\rho^P[\epsilon]=[\rho^P_1(\epsilon);...;\rho^P_n(\epsilon)],\quad
\widehat{\phi}_{P} (\omega)=
\min_{i=1,...,n} \phi_i(\omega)+\half \ln(n).
\]
Detector $\widehat{\phi}_P(\cdot)$ specifies a test
which accepts $H_0$, the observation being $\omega$, when $\widehat{\phi}_P(\omega)\geq0$ (i.e., with observation $\omega$, all pairwise tests with detectors $\phi_{i}$, $1\leq i\leq n$, $\chi=\pm1$, when deciding on $H_0$ vs. $H^{i}$, accept $H_0$), and accepts $H_1(\rho)$ otherwise.
\begin{proposition}\label{propPoisson} Whenever $\rho\geq \rho^P[\epsilon]$ and $\max_i\rho_i\leq R$,  the risk of the detector $\widehat{\phi}_P$ in the Poisson case of problem $({\cal D}_\rho)$ is $\leq\epsilon$. When
$\rho=\rho^P[\epsilon]$ and $\epsilon<1/4$, the test associated with $\widehat{\phi}_P$ is $\kappa_n$-rate optimal with
$\kappa_n=\kappa_n(\epsilon):={\ln(n/\epsilon^2)\over \ln(1/(4\epsilon))}.$ Note that $\kappa_n(\epsilon)\to 2$ as $\epsilon\to+0$.
\end{proposition}

\subsubsection{Gaussian case}\label{sec:gauss_multi}
Now let the distribution $P_\mu$ of  $\omega$ be normal with the mean $\mu$ and known variance $\sigma^2>0$, i.e. $\omega\sim \cN(\mu,\sigma^2 I)$. For the sake of simplicity, assume also that the (convex and compact) nuisance set $\cV$ is symmetric w.r.t. the origin.
In such a case,  the null hypothesis is\\
\be
H_0:\; \mu\in X:=\{\mu=Av ,\;v \in \cV\},
\ee{norm_null}
while the alternative $H_1(\rho)$ can be represented as the union, over $i=1,...,n$ and $\chi\in\{-1,1\}$, of $2n$ hypotheses
\be
\begin{array}{ll}
H^{\chi,i}(\rho_i):&\;\mu\in \chi Y_i(\rho_i)=\chi \left\{rA e[i]+Av:v\in\cV, \rho_i\leq r\leq R\right\}\\
\end{array}
\ee{gausgep}
(note that $\{x=re[i]+v :v \in\cV,-R\leq r\leq-\rho_i\}=-\{x=re[i]+v :v \in\cV,R\geq r\geq \rho_i\}$ due to $\cV=-\cV$).
Let $\hbox{ErfInv}(\cdot)$ be the inverse error function: $\hbox{Erf}(\hbox{ErfInv}(s))=s$, $0<s<1$. For $1\leq i\leq n$ and $\chi\in\{-1,1\}$,  let us set (cf. section \ref{sect:Gauss})
$$
\rho^G_{i}(\epsilon)=\max_{\rho,r,u,v}\left\{\rho: \begin{array}{l}\|A(u-re[i]-v)\|_2\leq \sigma \,\left[\ErfInv\left({\epsilon\over 4n}\right)+\ErfInv\left({\epsilon\over 2}\right)\right]\\
\chi r\geq\rho,\,u,v\in\cV\\
\end{array}\right\}
\eqno{(G^{i,\chi}_\epsilon)}
$$
(the left hand side quantity clearly is independent of $\chi$ due to $\cV=-\cV$),
and let
\begin{equation}\label{phiG}
\begin{array}{rcl}
\phi_{i,\chi}(\omega)&=&[A(u^{i,\chi}-r^{i,\chi}e[i]-v^{i,\chi})]^T\omega-\alpha_i,\\
\alpha_i&=&\lambda\, [A(u^{i,\chi}-r^{i,\chi}e[i]-v^{i,\chi})]^T[A(u^{i,\chi}+r^{i,\chi}e[i]+v^{i,\chi})],\\
\lambda&=&{\ErfInv\left({\epsilon\over 2}\right)\over \ErfInv\left({\epsilon\over 4n}\right)+\ErfInv\left({\epsilon\over 2}\right)},
\end{array}
\end{equation}
where $u^{i,\chi},v^{i,\chi},r^{i,\chi}$ are the $u,v,r$-components of an  optimal solution to $(G^{i,\chi}_\epsilon)$ (of course, in fact $r^{i,1}=-r^{i,-1}=\rho^G_i(\epsilon)$, and, besides, we can assume w.l.o.g. that $u^{i,-1}=-u^{i,1}$, $v^{i,-1}=-v^{i,1}$). Finally, let
\begin{equation}\label{eq:detection-proc}
\rho^G[\epsilon]=[\rho^G_1(\epsilon);...;\rho^G_n(\epsilon)],\quad
\widehat{\phi}_{G} (\omega)=
\min_{1\leq i\leq n,\chi=\pm1} \phi_{i,\chi}(\omega).
\end{equation}
Properties of the test associated with detector $\widehat{\phi}_G$ can be described as follows:
\begin{proposition}\label{pro:normal}
Whenever $\rho\geq \rho^G[\epsilon]$ and $\max_i\rho_i\leq R$,  the risk of the test $\widehat{\phi}_G$ in the Gaussian case of problem $({\cal D}_\rho)$ is $\leq\epsilon$. When
$\rho=\rho^G[\epsilon]$, the test is $\kappa_n$-rate optimal with
\[
\kappa_n=\kappa_n(\epsilon):={\ErfInv({\epsilon\over 4n})\over 2\ErfInv({\epsilon\over2})}+\half.
\]
Note that $\kappa_n(\epsilon)\to 1$ as $\epsilon\to +0$.
\end{proposition}

\paragraph{Remarks.} The results of Propositions \ref{propPoisson}, \ref{pro:normal}  imply that
testing procedures  $\widehat{\phi}_{G}$ and $\widehat{\phi}_{P}$  are
$\kappa_n$--rate optimal in the sense of the above definition
with $\kappa_n\asymp \sqrt{\ln n}$ in the Gaussian case and $\kappa_n\asymp \ln n$ in the Poisson case.
In particular, this implies that
the detection rates of these tests are within
a $\sqrt{\ln n}$ (resp., $\ln n$)--factor of the  rate profile $\rho^*$ of the ``oracle detector'' -- (the best) detection procedure which ``knows'' the node $\gamma\in \Gamma$ at which an event may occur.
This property of the proposed tests allows also for the following interpretation: consider the Gaussian problem setting in which the standard deviation $\sigma$ of noise is inflated by the factor $\kappa_n$. Then for every $i\in \{1,...,2n\}$ there is no test of hypothesis $H_0$ vs. $H^i(\rho_i)$ with risk $\le \epsilon$, provided that $\rho_i<\rho^G_i(\epsilon)$.
\mypar
Note that it can be proved that the price --
the $\sqrt{\ln n}$--factor -- for testing {multiple} hypotheses cannot
be eliminated at least in some specific settings \cite{GJTZ2008}.
\mypar
An important property of the proposed procedures is that they can be efficiently implemented -- when the nuisance set $\cV$ is computationally tractable (e.g., is a polyhedral convex set, an ellipsoid, etc.), the optimization problems $(G^{i,\chi}_\epsilon)$, $(P^i_\epsilon)$ are well structured {and} convex and thus can be efficiently solved using modern optimization tools even in relatively large dimensions.
\subsubsection{Numerical illustration: signal detection in the convolution model}\label{sec:numerics}
\begin{figure}
$$
\begin{array}{cc}
\resizebox{\twohundred}{\onefifty}{
 \includegraphics{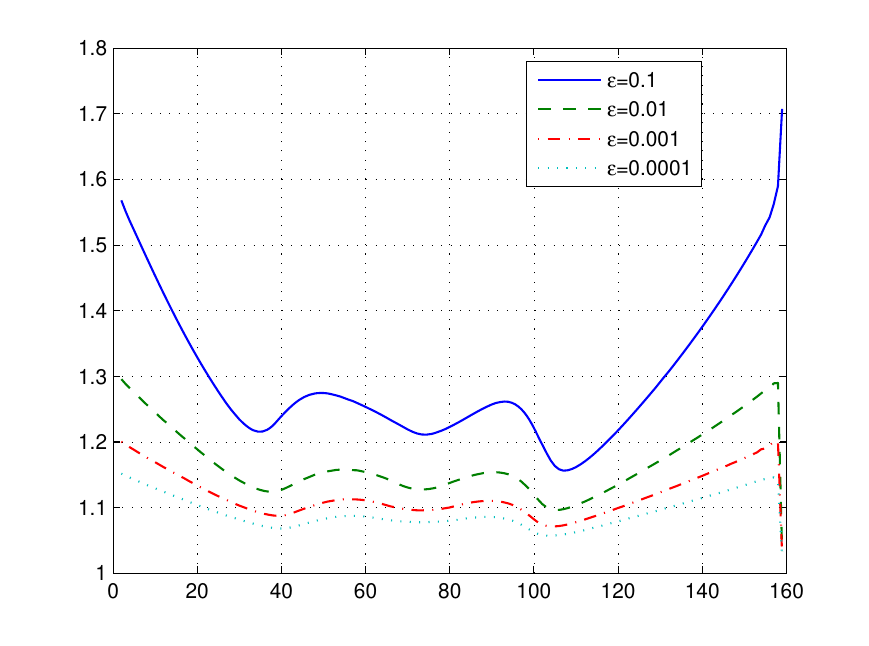}}
&
\resizebox{\twohundred}{\onefifty}{
 \includegraphics{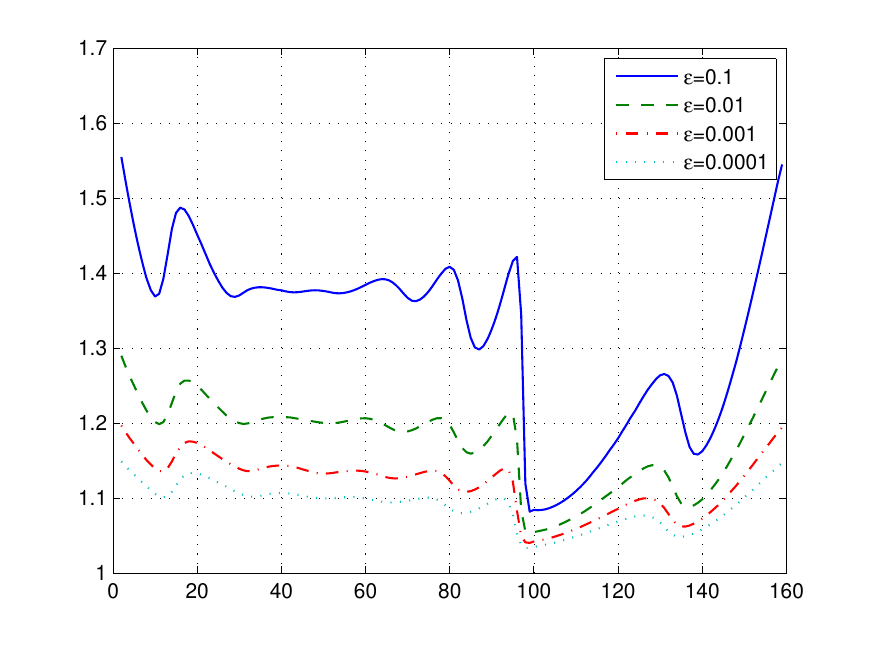}}
\\
{\rm (a)} & {\rm (b)}
\\
\resizebox{\twohundred}{\onefifty}{
 \includegraphics{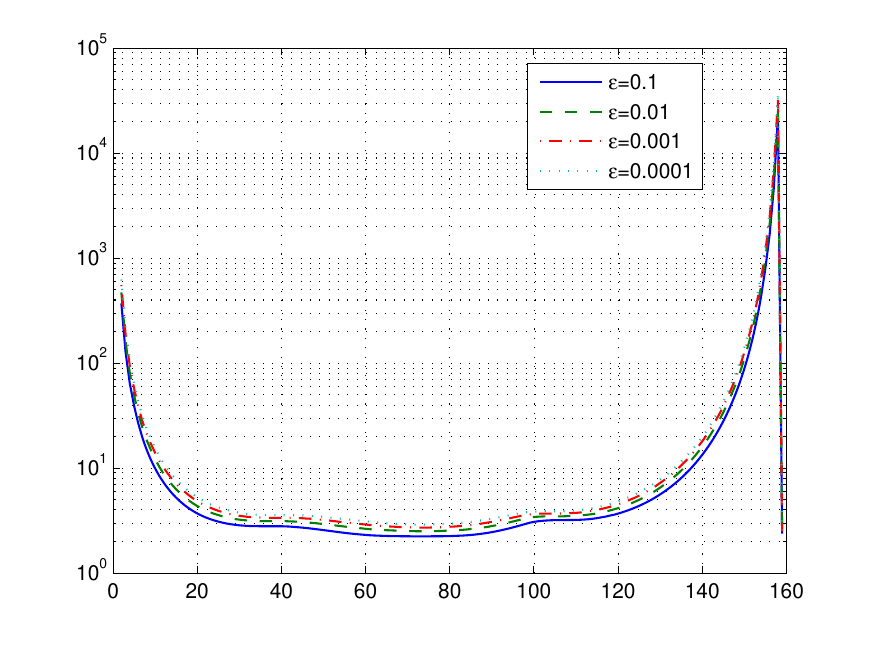}}
&
\resizebox{\twohundred}{\onefifty}{
 \includegraphics{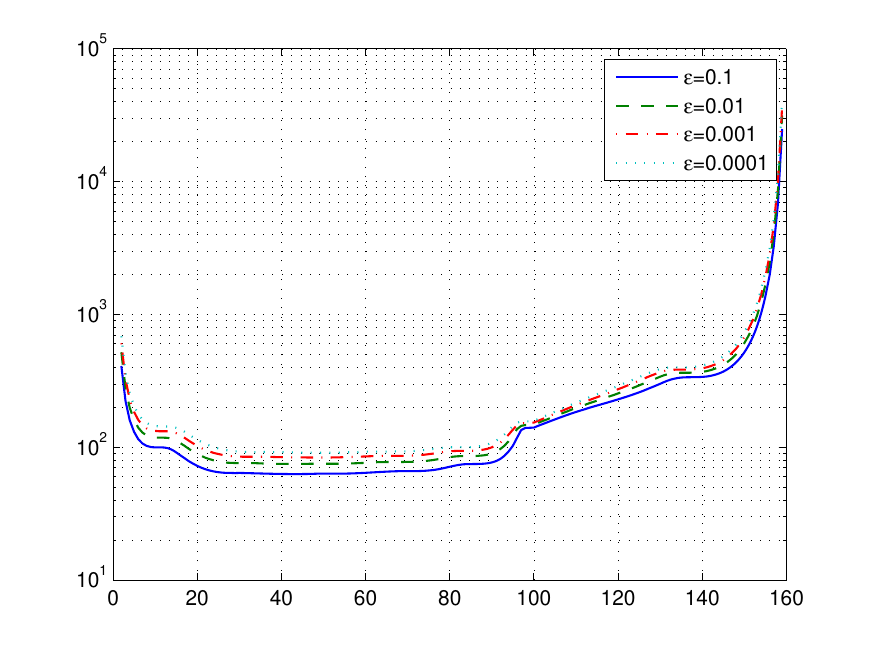}}
\\
{\rm (c)} & {\rm (d)}
\end{array}
$$
\caption{ \label{fig:2}  {\small The left pane (plots (a) and (c)) represents the experiment with ``step'' signals, the right pane (plots (c) and (d)) corresponds to the experiment with the signals which are proportional to basis orths. Nuisance parameter is set to $L=0.1$ and  $\sigma=1$ in both experiments. Plots (a) and (b): the value of $\rho^G[\epsilon]/\rho^*[\epsilon]$ for different values of $\epsilon$; plots (c) and (d): corresponding rate profiles (logarithmic scale).}}
\end{figure}

We consider here the ``convolution model'' with observation $\omega=A(s+v )+\xi$, where $s,v \in \bbr^n$, and $\xi\sim\cN(0,\sigma^2I_m)$ with known $\sigma>0$, and $A$ is as follows. Imagine that we observe at  $m$ consecutive moments the output of a discrete time
     linear dynamical system with a given impulse response (``kernel'') $\{g_k\}$ supported on a finite time horizon $k=1,...,T$. In this case,
     our observation $y\in\bbr^m$ is the linear image of $n$-dimensional ``signal'' $x$ which is system's input on the observation horizon, augmented by the input at $T-1$ time instants preceding this horizon (that is, $n=m+T-1$). $A$ is exactly the $m\times n$ matrix (readily given by $m$ and the kernel) of the just described linear mapping $x\mapsto y$.
     \mypar
We want to detect the presence of the signal $s=re[i]$, where $e[i],\;i=1,...,n,$ are some given vectors in $\bbr^n$. In other words, we are to decide between the hypotheses $H_0:\;\mu\in A\cV$ and $H_1(\rho)=\cup_{1\leq i\leq n,\chi=\pm1} H^{\chi,i}(\rho_i)$, with the hypotheses $H^{\chi,i}(\rho_i)$ defined in \rf{gausgep}. The setup for our experiment is as follow: we use $g_k=(k+1)^2(T-k)/T^3$, $k=0,...,T-1$, with $T=60$, and $m=100$, which results in $n=159$.
 The signatures $e[i]$, $1\leq i\leq n$ are the standard basic orths in $\bbr^n$ or unit step functions: $e_k[i]=1_{\{k\leq i\}}$, $k=1,...,n$, and the nuisance set $\cV$ is defined as $\cV_{L}=\{u\in \bbr^n:\;, |u_i-2u_{i-1}-u_{i-2}|\le L,\;i=3,...,n\}$, where $L$ is experiment's parameter.
 \mypar
  The
goal of the  experiment was to illustrate how large in the outlined problem is the  (theoretically, logarithmic in $n$) ``nonoptimality factor'' $\kappa_n(\epsilon)$ of the detector $\widehat{\phi}_G$, specifically, how it scales with the risk $\epsilon$.  To this end, we have computed, for different values of $\epsilon$, first, the ``baseline profile'' ---  the vector  with the entries\\
\be
\rho_i^*(\epsilon)=\max_{\rho,r,u,v}\left\{\rho:
\|A(u-re[i]-v)\|_2\leq 2\sigma \,\ErfInv(\epsilon/2),
r\geq\rho,\,u,v\in\cV\right\}
\ee{lower_G}
(cf. $(G^{i,1}_\epsilon)$); $\rho_i^*(\epsilon)$ is just the smallest $\rho$ for which the hypotheses $H_0$ and $H^{1,i}(\rho)$ can be distinguished with error probabilities $\leq\epsilon$ (recall that we are in the Gaussian case). Second, we computed the profile $\rho^G[\epsilon]$ of the test with detector $\widehat{\phi}_G$ underlying Proposition \ref{pro:normal}.
The results are presented on figure \ref{fig:2}.
Note that for $\epsilon\leq 0.01$ we have $\rho^G(\epsilon)/\rho^*(\epsilon)\le 1.3$ in the reported experiments.
\begin{figure}
$$
\begin{array}{cc}
\resizebox{\twohundred}{\onefifty}{
 \includegraphics{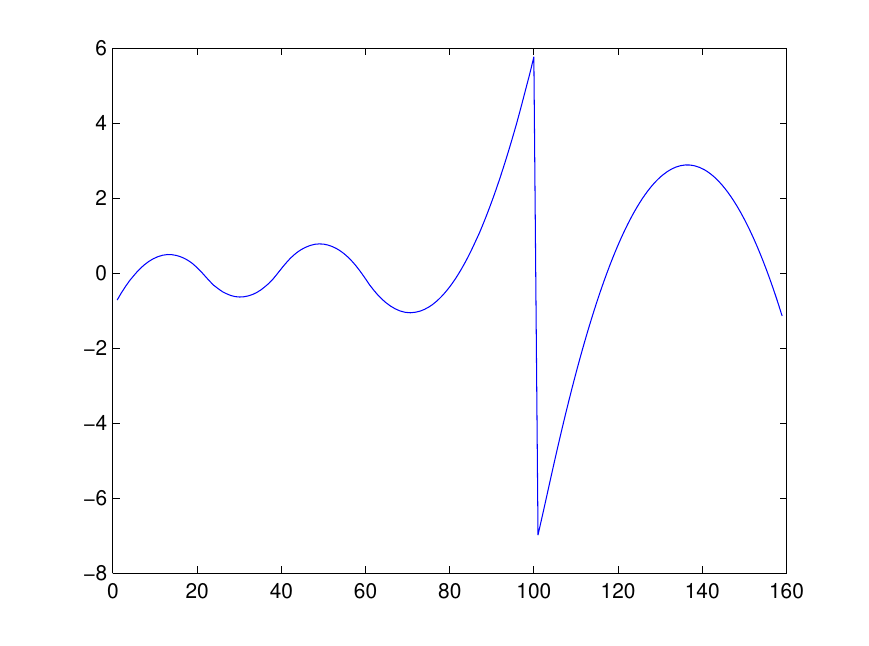}}
&
\resizebox{\twohundred}{\onefifty}{
 \includegraphics{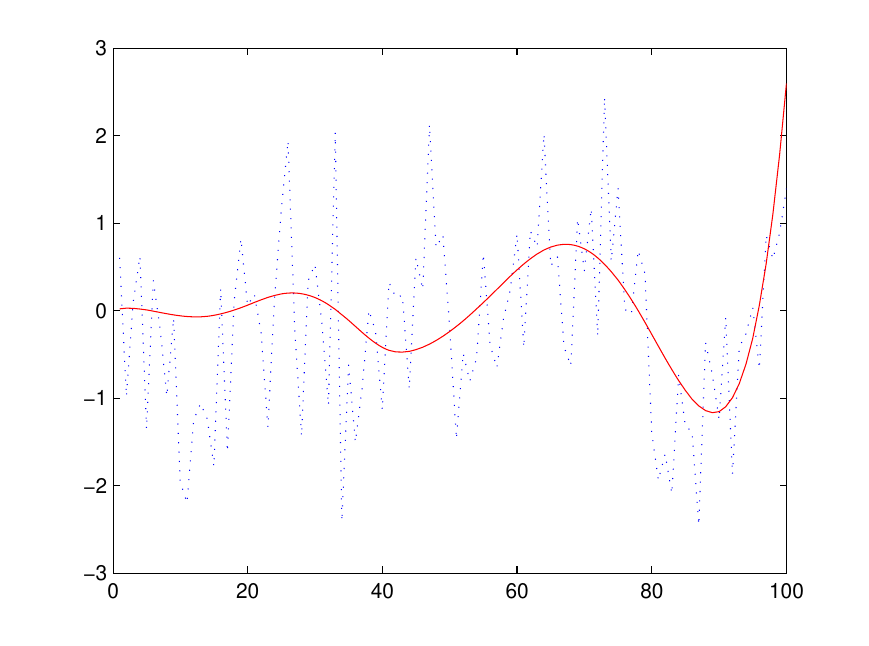}}
\end{array}
$$
\caption{ \label{fig:3}  {\small ``Hard to detect'' signal $\rho^G_{i}(\epsilon)e[i]+v^{i,1}-u^{i,1}$, where $\rho^G_{i}(\epsilon)$, $v^{i,1}$ and $u^{i,1}$ are components of an optimal solution to $(G^{i,\chi}_\epsilon)$ with $\epsilon=0.05$ and $i=100$ (left plot), and its image $Ax$ with a noisy observation (right plot). Experiment with ``step'' useful signals, nuisance parameter  $L=0.1$ and  $\sigma=1$.}}
\end{figure}
\paragraph{Quantifying conservatism.}While the baseline profile $\rho^*$ establishes an obvious lower bound for the $\rho$-profile of any test in our detection problem, better lower bounds can be computed by simulations. Indeed, let
\[
x_0^{i,\chi}=\chi u^{i},\;x_1^{i,\chi}=\chi \rho_i e[i]+v^{i},\; i=1,...,n, \;\chi\in \{-1,1\},
\] where $v^{i}$ and $u^{i}$ are some vectors in $\cV$. It is clear that the optimal risk
in the 
problem of distinguishing $H_0$ and $H_1(\rho)=\bigcup_{i=1}^n H^{\chi,i}(\rho_i)$
(cf. \rf{norm_null} and \rf{gausgep}) is lower bounded by the risk of distinguishing
\[
\bar{H}_0:\;\mu\in \{Ax_0^{i,\chi},i=1,...,n, \;\chi\in \{-1,1\}\},\;\mbox{and}\;
\bar{H}_1(\rho):\;\mu\in \{Ax_1^{i,\chi},i=1,...,n, \;\chi\in \{-1,1\}\},
\]
which, in its turn, is lower bounded by the risk of distinguishing of the hypothesis
$\tilde{H}_0:\;\mu=0$ from the alternative \[
\tilde{H}_1(\rho):\; \mu\in \{Az^{i,\chi},\;z^{i,\chi}=x_1^{i,\chi}-x_0^{i,\chi}=\chi (\rho_i e[i]+v^{i}-u^i),\; i=1,...,n, \;\chi\in \{-1,1\}\}.
\]
On the other hand, the latter risk is clearly bounded from below by the risk of the Bayesian test problem
as follows:
\[({\cal D}^\nu_\rho)\;\;\hbox{\begin{tabular}{l}Given
$\rho=[\rho_1;...;\rho_n]>0$, test the hypothesis $H_0:\;\mu=0$ against the\\
alternative $H^\nu_1(\rho):\;\mu= \chi A(\rho_i e[i]+v^{i}-u^i)$ with probability $\nu_{\chi i}$\\
where   $v^{i},u^i\in  \cV$, and  $\nu$ is a probability on $\{\chi i\}$,\;$i=1,...,n, \;\chi\in \{-1,1\}$.
 \end{tabular}}
 \]
 %
We conclude that the risk of deciding  between $H_0$ and $H_1(\rho)$ may be lower bounded by the risk of the optimal (Bayesian) test in the Bayesian testing problem $({\cal D}^\nu_\rho)$.
Note that we are completely free to choose the distribution $\nu$ and the points $u^i, v^i\in \cV$, $i=1,...,n$.
One can choose, for instance, $v^{\cdot,\chi}$ and $v^{\cdot,\chi}$ as components of an optimal solution to
\rf{lower_G} and a uniform on $\{\pm 1,...,\pm n\}$ prior probability $\nu$. Let us consider the situation where the matrix $A$ is an $n\times n$ Toeplitz matrix of periodic convolution on $\{1,...,n\}$ with kernel $g$, $g_k=({\koverT})^2(1-{\koverT})$, $k=1,..,T$, signatures $e[i]=e_{\cdot-i}$ are the shifts of the same signal $e_k=k/n$, $k=1,...,n$, and the nuisance set
\[
\cV_L=\{u\in \bbr^n:\;, |u_{i}-2u_{i-1{\,\rm mod\,} n}-u_{i-2{\,\rm mod\,} n}|\le L,\;i=1,...,n\}
\]
is symmetric and shift-invariant. Let us fix $\epsilon>0$ and choose $v^{i}=-u^i$ as components of an optimal solution to the corresponding optimization problem $(G^{i,\chi}_\epsilon)$. Because of the shift-invariance of the problem setup  the optimal values $\rho^*_i(\epsilon)$ and $\rho^G_i(\epsilon)$ do not depend on $i$ and are equal to the same $\rho^*(\epsilon)$ and, respectively, $\rho^G(\epsilon)$, and all $v^i$ are the shifts of the same $v\in \bbr^n$. In this case the risk of the Bayesian test corresponding to the uniform on $\{\pm 1,..., \pm n\}$ prior distribution $\nu$ is a lower bound of the optimal risk for the corresponding
detection problem $({\cal D}_\rho)$.
\par
\begin{figure}
$$
\begin{array}{cc}
\resizebox{\twohundred}{\onefifty}{
 \includegraphics{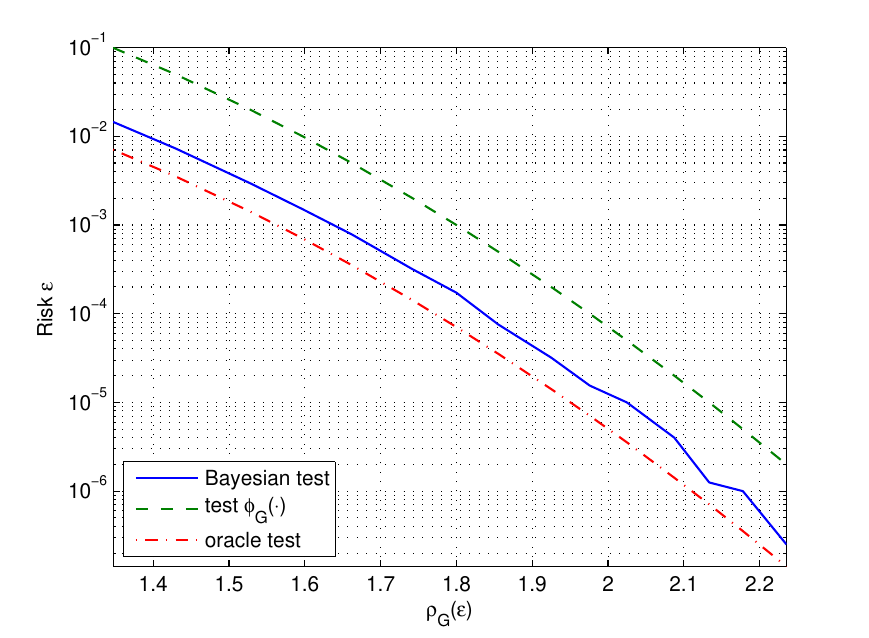}}
&
\resizebox{\twohundred}{\onefifty}{
 \includegraphics{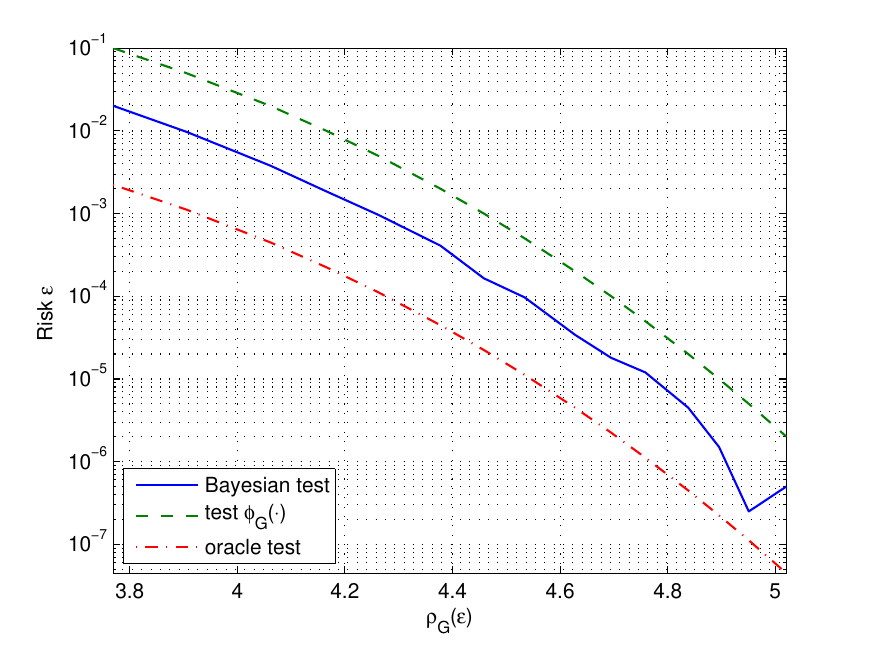}}
\end{array}
$$
\caption{ \label{fig:4}  {\small Estimated risk of the Bayes test as a function of test rate $\rho_G(\epsilon)$, compared to the risk of the baseline test and that of the simple test with data \rf{phiG} ($L=0.01$ and  $\sigma=1$). Simulation for $n=100$ (left plot) and $n=1000$ (right plot).}}
\end{figure}

On figure \ref{fig:4} we present the results of two simulation for $n=100$ and $n=1000$, the value $L=0.01$ of the parameter of the nuisance class, and $\sigma=1$. For different values of $\epsilon$ we have first computed corresponding rates  $\rho^*(\epsilon)$ and $\rho^G(\epsilon)$, as well as components $v^i=-u^i$ of the optimal solution (recall that due to the shift-invariance of the problem, $v^i_k=v^1_{k-i+1 \,{\rm mod} n}$). Then an estimation of the risk of the Bayesian test with the uniform prior is computed over $N=10^7$ random draws. Note that already for $\epsilon=0.01$ rate $\rho^G(\epsilon)$ of the simple test is only $7\%$ higher than the corresponding Bayesian lower bound for $n=1000$ ($15\%$ for $n=100$).

\subsubsection{Numerical illustration: signal identification in the convolution model}\label{signalident}
The purpose of the experiment we report on in this section is to illustrate an application of the approach to multiple hypotheses testing presented in section \ref{multiple}. The experiment in question is a modification of that described in section \ref{sec:numerics}, the  setup is as follows.
On time horizon $t=1,...,m$, we observe the output, contaminated by noise, of a discrete-time linear dynamic system  with ``memory'' $T$ (that is, the impulse response $g$ is zero before time 0 and after time $T-1$). The input $x$ to the system is an impulse of amplitude $\geq\rho>0$ ($\rho$ is known) at unknown time $\tau$ known to satisfy $-T+2\leq \tau\leq m$. Setting $n=m+T-1$, our observation is
$$
\omega=[\omega_1;...;\omega_{m}]=Ax+\xi,\,\,\xi\sim\cN(0,I_m),
$$
with $m\times n$ matrix $A$ readily given by the impulse response $g$. We have $n$ hypotheses about $x$, the $i$-th of them stating that $x\in X_i=\{x=re_i,r\geq\rho\}$, where $e_i$, $i=1,...,n$, are the standard basic orths in $\bbr^n$. Given an observation, we want to decide to which of the sets $X_1,...,X_n$  the actual input belongs, that is, we need to distinguish between $n$ hypotheses $H_1,...,H_n$ on the distribution of $\omega$, with $H_i$ stating that this distribution is $\cN(Ax,I_m)$ for some $x\in X_i$.
\par
The problem can be processed as follows. Let us choose two nonnegative integers $\mu$ (``margin'') and $\nu$ (``resolution''), and imagine that we do not care much about distinguishing between the ``boundary hypotheses'' $H_i$ (those with $i\leq \mu$ and with $i\geq n-\mu+1$) and all other hypotheses, same as we do not care much about distinguishing between ``close to each other'' hypotheses $H_i$ and $H_j$, those with $|i-j|\leq\nu$. What we do care about is not to miss the true hypothesis and to reject any non-boundary hypothesis which is not close to the true one. Note that when $\mu=\nu=0$, we ``care about everything;'' this, however, could require large amplitude $\rho$ in order to get a reliable test, since the impulses at times $t$ close to the endpoints of the time segment $-T+2\leq t\leq m$ could be poorly observed, and impulses at close to each other time instants could be difficult to distinguish.  Operating with positive margins and/or resolutions, we, roughly speaking, sacrifice the ``level of details'' in our conclusions in order to make these conclusions reliable for smaller values of the amplitude $\rho$.
\par
With the approach developed in section \ref{multiple}, our informally described intentions can be formalized as follows. In the terminology and notation of section \ref{multiple}, let us define the set $\C$ of pairs $(i,j)$, $1\leq i,j\leq n$, $i\neq j$, i.e., the pairs with ``$H_j$  close to $H_i$,'' as follows:
\begin{itemize}
\item for a ``boundary hypothesis'' $H_i$ (one with $1\leq i\leq \mu$ or $n-\mu+1\leq i\leq n$), every other hypothesis $H_j$ is close to $H_i$;
\item for a ``non-boundary hypothesis'' $H_i$ (one with $1+\mu\leq i\leq n-\mu$), close to $H_i$ hypotheses $H_j$ are those with $1\leq|i-j|\leq\nu$.
\end{itemize}
Detectors $\phi_{ij}(\omega)$   we intend to use are the Gaussian log-likelihood detectors
\begin{equation}\label{phiijs}
\begin{array}{c}
\phi_{ij}(\omega)={1\over 2}[\xi_{ij}-\eta_{ij}]^T\omega+{1\over 4}[\eta_{ij}^T\eta_{ij}-\xi_{ij}^T\xi_{ij}],\\
\xi_{ij}=Ax_{ij},\,\eta_{ij}=Ay_{ij}, [x_{ij}=y_{ij}]=\argmin_{r,s}\{\|rAe_i-sAe_j\|_2:r\geq\rho,s\geq\rho\},\\
\end{array}
\end{equation}
which allows to specify the quantities $\epsilon_{ij}$ in (\ref{detectorsasusual1}) as
\begin{equation}\label{epsilonsij}
\epsilon_{ij}=\exp\{-(\xi_{ij}-\eta_{ij})^T(\xi_{ij}-\eta_{ij})/8\},
\end{equation}
see section \ref{sect:Gauss}.
\par
Applying the construction from section \ref{multiple}, we arrive at a risk bound $\epsilon$ and a test which, given an observation $\omega$, accepts some of the hypotheses $H_i$, ensuring the following. Let the true hypothesis be $H_{i_*}$. Then (all probabilities are taken according to the distribution specified by $H_{i_*}$)
\begin{itemize}
\item[A.] The probability for $H_{i_*}$ to be rejected by the test is at most $\epsilon$;
\item[B.] The probability of the event that the list of accepted hypotheses contains a hypothesis $H_j$ such that both $H_j$ is {\sl not} close to $H_{i_*}$ and $H_{i_*}$ is not close to $H_j$ is at most $\epsilon$.\par
Note that with our definition of closeness, the latter claim implies that when $H_{i_*}$ is not a boundary hypotheses, the probability for the list of accepted hypotheses to contain a {\sl non-boundary} hypothesis $H_j$ with $|i-j|>\nu$ is at most $\epsilon$.
\end{itemize}
The outlined model demonstrates the potential of {\sl asymmetric} closeness: when a boundary hypothesis is difficult to distinguish from other hypotheses, it is natural to declare all these hypotheses to be close to the boundary one. On the other hand, there are no reasons to declare a boundary hypothesis to be close to a well identifiable ``inner'' hypothesis.
\par
As we have seen in section \ref{multiple}, given $\rho$, the risk $\epsilon$ can be efficiently computed via convex optimization, and we can use this efficient computation to
find the smallest amplitude $\rho$ for which $\epsilon$ takes a given target value $\varepsilon$. This is what was done in the numerical experiment we are about to report. In this experiment, we used $T=m=16$ (i.e., the number of hypotheses $n$ was 31), and the impulse response was similar to the one reported earlier in this section, namely the nonzero entries in $g$ were
$$
g_t=\alpha (t+1)^2(T-t),\,0\leq t\leq T-1,
$$
while $\alpha$ was selected to ensure $\max_tg_t=1$. For various values of margins $\mu$ and resolutions $\nu$, we computed the minimal amplitude $\rho=\rho(\mu,\nu)$ which still allowed for our test to guarantee  risk $\epsilon\leq0.01$. The results are presented in table \ref{Pulstable}. A simple lower bound $\underline{\rho}(\mu,\nu)$
on the smallest $\rho$ such that there exists ``in the nature''  a test capable to ensure A and B with $\epsilon=0.01$, amplitudes of impulses being $\rho$, may be constructed by lower bounding the probability of a union of events by the largest among the probabilities of these events.
 In the table we present, along with the values of $
\rho(\cdot,\cdot)$, the ``excess value'' $\rho(\mu,\nu)/\underline{\rho}(\mu,\nu)-1$. Observe that while $\rho(\mu,\nu)$ itself strongly depends on the margin $\mu$, the excess is nearly independent of $\mu$ and $\nu$. Of course, $40\%$ excess is unpleasantly large; note, however, that the lower bound $\underline{\rho}$ definitely is optimistic. In addition, this ``overly pessimistic'' excess decreases as the target value of $\epsilon$ decreases; what was 40\% for $\varepsilon=0.01$, becomes 26\% for $\varepsilon=0.001$  and $19\%$ for $\varepsilon=1$.e-4. \par
In the reported experiment, along with identifying $\rho(\cdot,\cdot)$, we were interested also in the effect of optimal shifts $\phi_{ij}(\cdot)\mapsto\phi_{ij}(\cdot)-\bar{\alpha}_{ij}$, see section \ref{multiple}. To this end we compute the smallest $\rho=\widetilde{\rho}(\mu,\nu)$ such that the version of our test utilizing $\alpha_{ij}\equiv 0$ is capable to attain the risk $\varepsilon=0.01$. Table \ref{Pulstable} presents, along with other data, the ratios $\widetilde{\rho}(\mu,\nu)/\rho(\mu,\nu)$ which could be considered as quantifying the effect of shifting the tests. We see that the effect of the shift is significant when the margin $\mu$ is positive.
\begin{table}
$$
\begin{array}{|c|c|c|c|c|}\cline{2-5}
\multicolumn{1}{c|}{}&\nu=0&\nu=1&\nu=2&\nu=3\\
\hline
\mu=0&\begin{array}{c}276.0 (+40.1\%)\\
1.00\\
\end{array}&\begin{array}{c}71.0 (+40.0\%)\\
1.00\\
\end{array}&\begin{array}{c}31.5 (+40.4\%)\\
1.00\\
\end{array}&\begin{array}{c}18.1 (+44.4\%)\\
1.03\\
\end{array}\\
\hline
\mu=1&\begin{array}{c}133.2 (+40.5\%)\\
1.88\\
\end{array}&\begin{array}{c}48.0 (+40.5\%)\\
1.48\\
\end{array}&\begin{array}{c}23.6 (+40.3\%)\\
1.33\\
\end{array}&\begin{array}{c} 14.1(+40.5\%)\\
1.25\\
\end{array}\\
\hline

\mu=2&\begin{array}{c}102.0(+40.2\%)\\
1.44\\
\end{array}&\begin{array}{c}36.8 (+40.0\%)\\
1.93\\
\end{array}&\begin{array}{c}19.4 (+40.3\%)\\
1.64\\
\end{array}&\begin{array}{c}11.9 (+40.1\%)\\
1.48\\
\end{array}\\
\hline
\mu=3&\begin{array}{c}77.5 (+40.1\%)\\
1.33\\
\end{array}&\begin{array}{c}29.8 (+40.0\%)\\
1.61\\
\end{array}&\begin{array}{c}16.3 (+40.3\%)\\
1.94\\
\end{array}&\begin{array}{c}10.4 (+40.1\%)\\
1.70\\
\end{array}\\
\hline
\end{array}
$$
\caption{\label{Pulstable} Identifying signals in the convolution model. In a cell, top: $\rho(\mu,\nu)$ and excess $\rho(\mu,\nu)/\underline{\rho}(\mu,\nu)-1$ (in brackets, percents); bottom: $\widetilde{\rho}(\mu,\nu)/\rho(\mu,\nu)$.}
\end{table}
\corr{}{
\begin{figure}
$$
\begin{array}{ccc}
\epsfxsize=140pt\epsfysize=140pt\epsffile{LnEps_mu3_nu3}&\epsfxsize=140pt\epsfysize=140pt\epsffile{LnEpsShifted_mu3_nu3}&
\epsfxsize=140pt\epsfysize=140pt\epsffile{Shifts_mu3_nu3}\\
\end{array}
$$
\caption{\label{figscale} Optimal shift of detectors (\ref{phiijs}), $\mu=\nu=3,\rho=\rho(3,3)$. Left: $\ln(\epsilon_{ij})$, see  (\ref{epsilonsij}); middle: $\ln(\epsilon_{ij})+
\alpha_{ij}$, where $\alpha_{ij}$ are optimal for (\ref{selectingalphas}); right: matrix $\alpha_{ij}$ of shifts as a function of $i,j$.}
\end{figure}}
\subsection{Testing from indirect observations}
  \subsubsection{Problem description}
  \label{sec:lf-test}
  Let $\F$ be a class of cumulative distributions on $\bbr$. Suppose that for $\ell=1,...,L$, we are given $K_\ell$ independent realizations of random variable $\zeta^\ell$.
We assume that the c.d.f. $F_{\zeta^\ell}$ of $\zeta^\ell$ is a linear transformation of unknown c.d.f. $F_\xi$ of ``latent'' random variable $\xi$,  $F_\xi\in \F$. In this section we consider two cases of the sort; in both of them, $\eta^\ell$ is an independent of $\xi$ random variable (``nuisance'') with known c.d.f. $F_{\eta^\ell}$. In the first case (``deconvolution model''), $\zeta^\ell= \xi+\eta^\ell$,  so that the distribution of $\zeta^\ell$ is $F_{\zeta^\ell}(t)=\int_{\bbr} F_\xi(t-s)dF_{\eta^\ell}(s)$. In the second case (``trimmed observations''), observations  are trimmed:
$\zeta^\ell= \max\{\xi,\eta^\ell\}$, so that $F_{\zeta^\ell}(t)=F_\xi(t)F_{\eta^\ell}(t)$.
\par
We consider here the testing problem where our objective is to test, for given $t\in \bbr,\;\alpha\in (0,1)$ and $\rho>0$, the hypotheses\footnote{A related problem of {\em estimation} of the c.d.f. $F_\xi$  in the deconvolution model, a special case of linear functional estimation \cite{Donoho1987,Donoho1991,JN2009}, have received much attention in the statistical literature (see, e.g., \cite{gaffey1959,zhang1990,fan1991,dattner2011} and \cite[Section~2.7.2]{meister2009}
for a recent review of corresponding contributions).}
\[
H_{{1}}:\;F_\xi(t)<\alpha-\rho\;\;\mbox{and}\;\;H_{{2}}: \;F_\xi(t)>\alpha+\rho \eqno{(C_{\alpha,t}[\rho])}
\]
given observations $\zeta_k^\ell,\;k=1,...,K_\ell,\;\ell=1,...,L.$
\par
Under minor regularity conditions on $F_{\eta^\ell}$ and $F_\xi$, $(C_{\alpha,t}[\rho])$ may be approximated by the discrete decision problem as follows.
Let $\xi$ be a  discrete random variable with unknown distribution $x$ known to belong to a given closed convex subset $\X$ of $n$-dimensional probabilistic simplex.
We want to infer about $x$ given indirect observations of $\xi$ obtained by $L$ different ``observers'': the observations $\omega^\ell_i$, $i=1,...,K_\ell$ of  $\ell$-th observer are independent realizations of random variable $\omega^\ell$ taking values $1,...,m_\ell$ with distribution $\mu^\ell=A^\ell x$, where $A^\ell$ is a known stochastic matrix.
For instance, when $\xi$ takes values $1,...,n$ and $\omega^\ell=\xi+\eta^\ell$ with nuisance $\eta^\ell$ taking values $1,...,n_\ell$ and distribution $u^\ell$, $A^\ell$ is $(n_\ell+n-1)\times n$ matrix,  and the nonzero entries of the matrix are given by
$A^\ell_{ij}=u^\ell_{i-j+1},\;1\leq j\leq i\leq j+n_\ell-n$. We assume in the sequel that $A^\ell x>0$ whenever $x\in\X$, $1\leq \ell \leq L$.
\par
Let  $g(x)=g^Tx$, $g\in \bbr^n$, be a given linear functional of the distribution $x$.
Given
$\alpha$ and $\rho>0$, our goal is to decide on
the hypotheses about the distribution $x$ of $\xi$
\[
H_{{1}}[\rho]:\; x\in \X,\,g(x)\leq \alpha-{\rho},\;\; H_{{2}}[\rho]:\; x\in\X,\,g(x)\geq \alpha+{\rho}.\eqno{(\D_{g,\alpha}[\rho])}\]
given observations $\omega^1,...,\omega^\ell$. We denote by $\rho_{\max}$ the largest $\rho$ for which both these hypotheses are nonempty, and assume from now on that $\rho_{\max}>0$ (as far as our goal is concerned, this is the only nontrivial case).
Now let us fix $0<\epsilon<1$ and, given a decision rule $T(\cdot)$, let us denote $\rho_T[\epsilon]$ the smallest $\rho\geq0$ such that the risk of the rule $T(\cdot)$ in the problem $(\D_{g,\alpha}[\rho])$ does not exceed $\epsilon$. We refer to $\rho_T[\epsilon]$ as the $\epsilon$-resolution of $T(\cdot)$ and denote by $\rho^*[\epsilon]=\inf_{T(\cdot)}\rho_T[\epsilon]$ (``$\epsilon$-rate'') the best $\epsilon$-resolution achievable in our problem.
 Our goal is given $\epsilon$, to design a test with $\epsilon$-resolution close to $\rho^*[\epsilon]$.
\par
The resulting observation scheme fits the definition of the direct product of Discrete observation schemes of section \ref{secrepeatedI} -- we have $K=\sum_{\ell=1}^LK_\ell$ ``simple'' (or $L$ $K_\ell$-repeated) Discrete observation schemes, the $k$-th scheme yielding the observation $\omega_k$, $k=1,...,K$, of one of $L$ types.
\par
Given an $\epsilon\in(0,1)$, we put
\be
\rho[\epsilon]=\max_{x,y,r}
\left\{
r:\;\begin{array}{l}
\sum_{\ell=1}^LK_\ell\ln\left(\sum_{i=1}^{m_\ell}\sqrt{[A^\ell x]_i[A^\ell y]_i }\right)
\geq \ln \epsilon,\\
x,y\in \X,\;g(x)\leq \alpha-r,\;g(y)\geq \alpha+r.
\end{array}\right\}
\ee{rho_A}
Clearly, $0\leq \rho[\epsilon]\leq\rho_{\max}$ due to $\rho_{\max}>0$. {\sl We assume from now on that $\rho[\epsilon]<\rho_{\max}$}.
Let now $\rho\in[\rho[\epsilon],\rho_{\max}]$. Consider the optimization problem
\[
\Opt[\rho]=\max_{x,y}
\left\{\Psi(x,y):\;\begin{array}{l}
\Psi(x,y)=\sum_{\ell=1}^LK_\ell\ln\left(\sum_{i=1}^{m_\ell}\sqrt{[A^\ell x]_i[A^\ell y]_i }\right),\\
x,y\in \X,\;g(x)\leq \alpha-\rho,\;g(y)\geq \alpha+\rho.
\end{array}\right\}.
\eqno{(F_{g,\alpha}[\rho])}
\]
This problem is feasible (since $\rho\leq\rho_{\max})$ and thus solvable, and from $\rho\geq\rho[\epsilon]$ and $\rho[\epsilon]<\rho_{\max}$ it easily follows (see item 1$^0$ in the proof of Proposition \ref{pro:cdf}) that $\Opt[\rho]\leq\epsilon$. Let $(x_{\rho},y_{\rho})$ be an optimal solution.
Consider a simple test $\widehat{T}_\rho$ given by the detector $\widehat{\phi}(\cdot)$,
\begin{equation}\label{testwidehatTrho}
\widehat{\phi}(\omega)=\widehat{\phi}_{\rho}(\omega):=\sum_{k=1}^K\phi_k(\omega_k),\,\,\phi_k(\omega_k)=\half \ln\left([A^{\ell(k)} x_{\rho}]_{\omega_k}/[A^{\ell(k)} y_{\rho}]_{\omega_k}\right),
\end{equation}
 with $\ell(k)$ uniquely defined by the relations
 \[\sum_{\ell<\ell(k)}K_\ell<k\leq\sum_{\ell\leq\ell(k)}K_\ell.
\]
We have the following simple corollary of Proposition \ref{propnonstI}:
\begin{proposition}\label{pro:cdf} Assume that $\rho_{\max}>0$ and $\rho[\epsilon]<\rho_{\max}$, and let $\epsilon\in(0,1/4)$. Then
\begin{equation}\label{wehaveq7}
\rho[\epsilon]\leq \vartheta(\epsilon)\rho^*[\epsilon],\;\;
\vartheta(\epsilon)={2\ln(1/\epsilon)\over
\ln [1/(4\epsilon)]}.
\end{equation}
In other words, there is no decision rule in the problem $(\D_{g,\alpha}[\rho])$  with the risk $\leq \epsilon$ if $\rho<\rho[\epsilon]/\vartheta(\epsilon)$.
\par
On the other hand, when $\rho\in[\rho[\epsilon],\rho_{\max}],$ the risk of the simple test $\widehat{\phi}_\rho$ in the problem $({\cal D}_{g,\alpha}[\rho])$ is $\leq$ $\exp\left(\Opt[\rho]\right)\leq\epsilon$.

\end{proposition}
Note that $\vartheta(\epsilon)\to 2$ as $\epsilon\to 0$.
 Under the premise of Proposition \ref{pro:cdf}, the test associated with detector $\widehat{\phi}_{\rho[\epsilon]}(\cdot)$ is well defined and distinguishes between the hypotheses $H_{{1}}[\rho[\epsilon]]$, $H_{{2}}[\rho[\epsilon]]$ with risk $\leq\epsilon$. We refer to the quantity $\rho[\epsilon]$ as to {\sl resolution} of this test.

\subsubsection{Numerical illustration}
We present here some results on numerical experimentation with the testing problem $(C_{\alpha,t}[\rho])$. For the sake of simplicity, we suppose that the distributions with c.d.f.'s from $\F$ are supported on $[-1,1]$. We start with an appropriate discretization of the continuous problem.
\paragraph{Discretizing continuous model.}
\begin{enumerate}
\item  Let $n\in \bbz_+$, and let $-1=a_0<a_1<a_2<...<a_{n}=1$ be a partition of $(-1,1]$ into $n$ intervals
$I_i=({a}_{i-1},a_i]$, $i=1,...,n$.
We associate with a c.d.f. $F\in F$  the $n$-dimensional probabilistic vector $x=x[F]$ with the entries $x_k=\Prob_{\xi\sim F}\{\xi\in I_k\}$ and $\bar{a}_k=(a_{k-1}+a_k)/ 2$, the central point of $I_k$, $k=1,...,n$, and denote by $\F_n$ the image of $\F$ under the mapping $F\mapsto x[F]$.
\item We build somehow a convex compact subset $\X\supset \F_n$ of the $n$-dimensional probabilistic simplex.
\item Depending on the observation scenario, we act as follows.
\begin{enumerate}\item
{\em Deconvolution problem:} $\zeta^\ell$ satisfy $\zeta^\ell=\xi+\eta^\ell$.
Let $0<\delta<1$ (e.g., $\delta=K_\ell^{-1}$), $m_\ell\in \bbz_+$, and let
\[
b^\ell_1=a_0+q_{\eta^\ell}(\delta),\;
b^\ell_{m_\ell-1}=a_n+q_{\eta^\ell}(1-\delta),
\]
where $q_{\eta^\ell}(p)$ is the $p$-quantile of $\eta_\ell$. Note that
$\prob\{\zeta^\ell\notin 
[b^\ell_1,b^\ell_{m_\ell-1}]
\}\le 2\delta.$
Let now $-\infty=b^\ell_0<b^\ell_1<b^\ell_2<...<b^\ell_{m_\ell-1}<b^\ell_m=\infty$ be a partition of $\bbr$ into $m_\ell$ intervals $J^\ell_i=(b^\ell_{i-1},b^\ell_i]$, $i=1,...,m_\ell-1$, $J_{m_\ell}=(b^{\ell}_{m_\ell-1}, \infty)$. We put
$\mu^\ell_i=\Prob\{\zeta\in J_i\}$, $i=1,...,m_\ell$ and define the $m_\ell\times n$ matrix stochastic matrix $A^\ell=(A^{\ell}_{jk})$ with elements
 \bse
 A^\ell_{ij}=\Prob\left\{{\bar{a}_j}+\eta^\ell\in  J_i\right\},%
\ese
the approximations of conditional probabilities
$\Prob\{\zeta^\ell \in J_i|\xi\in I_j\}$.
\item {\em Trimmed observations:} $\zeta^\ell=\max\{\xi,\eta^\ell\}$. We partition $\bbr$ into $m_\ell=n+1$ intervals, $I_i,\;i=1,...,n$ as above and an ``infinite bin'' $I_{n+1}=(a_n,a_{n+1}=\infty)$. We put
$\mu^\ell_i=\Prob\{\zeta\in J_i\}$, $i=1,...,m_\ell$ and define the $m_\ell\times n$ matrix $A^\ell$ with elements
\[A^\ell_{ij}=\delta_{ij}\Prob\{\eta^\ell\leq a_j\}+1_{\{i>j\}}\Prob\{\eta^\ell\in I_i\},
\]
where $\delta_{ij}=1$ if $i=j$ and zero otherwise,
which are the estimates of the probability of $\zeta^\ell$ to belong to $I_i$, given that $\xi\in I_j$.
\end{enumerate}
\item We denote $g=g(t)\in \bbr^n$, with entries $g_i=1_{\{\bar{a}_i\le t\}}$,
$i=1,\ldots, n$, so that
$g^Tx$ is an approximation of $F(t)$.
\item Finally, we consider discrete observations $\omega^\ell_k\in \{1,...,m_\ell\}$,
\[
\omega_k^\ell=i\,1_{\{\zeta^\ell_k\in J^\ell_i\}}
\;\; k=1,...,K^\ell,\; \ell=1,...,L.
\]
\end{enumerate}
We have specified the data of a testing problem of the form $(\D_{g,\alpha}[\rho])$. Note that the discrete observations we end up with are deterministic functions of the ``true'' observations $\zeta^\ell$, so that a test for the latter problem induces a test for the problem of interest $(C_{\alpha,t}[\rho])$.
When distributions from $\F$, same as distributions of the nuisances $\eta^\ell$,
possess some regularity, and  the partitions $(I_i)$ and $(J_i)$ are
``fine enough'',
the problem $(\D_{g,\alpha}[\rho])$ can be considered as a good proxy of the problem of actual interest.
\par
\paragraph{Simulation study.}\label{sec:numeric}
We present results for three distributions of the nuisance:
\begin{description}
\item{(i)} Laplace distribution ${\cal L}(\mu, a)$ (i.e., the density ${(2a)^{-1}}e^{-{|x-\mu|/a}}$) with parameter $a=\half$ and $\mu=0$;
\item{(ii)} distribution $\Gamma(0,2,1/(2\sqrt{2}))$
with the location 0, shape parameter $2$ and the scale $\frac{1}{2\sqrt{2}}$
(the standard
deviation of the error is equal to $0.5$).\footnote{Recall that $\Gamma$-distribution with parameters $\mu$,
$\alpha$, $\theta$ has the density $
[\Gamma(\alpha)\theta^\alpha]^{-1} (x-\mu)^{\alpha-1}
\exp\{-(x-\mu)/\theta\}1_{\{x\ge \mu\}}$.}
\item{(iii)} mixture of Laplace distributions
$\half {\cal L}(-1,\half)+\half{\cal L}(1,\half)$.
\end{description}
The interval $[-1,1]$ was split into $n=100$ bins of equal lengths. The discretized distributions $x=x[F]$, $F\in\F$, are assumed to have bounded second differences, specifically, when denoting $h$ the length of the bin,
\[
|x_{i+1}-2x_i+x_{i-1}|\leq h^2\LL, \;i=2,...,n-1;
\]
in the presented experiments, $\X$ is comprised of all probabilistic vectors satisfying the latter relation with $\LL=0.4$.

\begin{figure}
$$
\begin{array}{cc}
\epsfxsize=220pt\epsfysize=160pt\epsffile{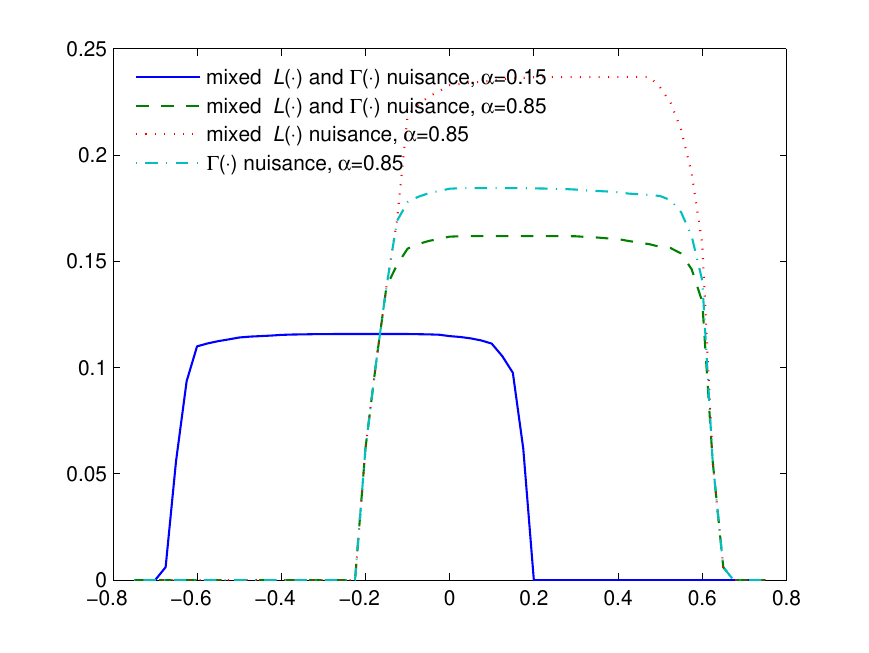}
&\epsfxsize=220pt\epsfysize=160pt\epsffile{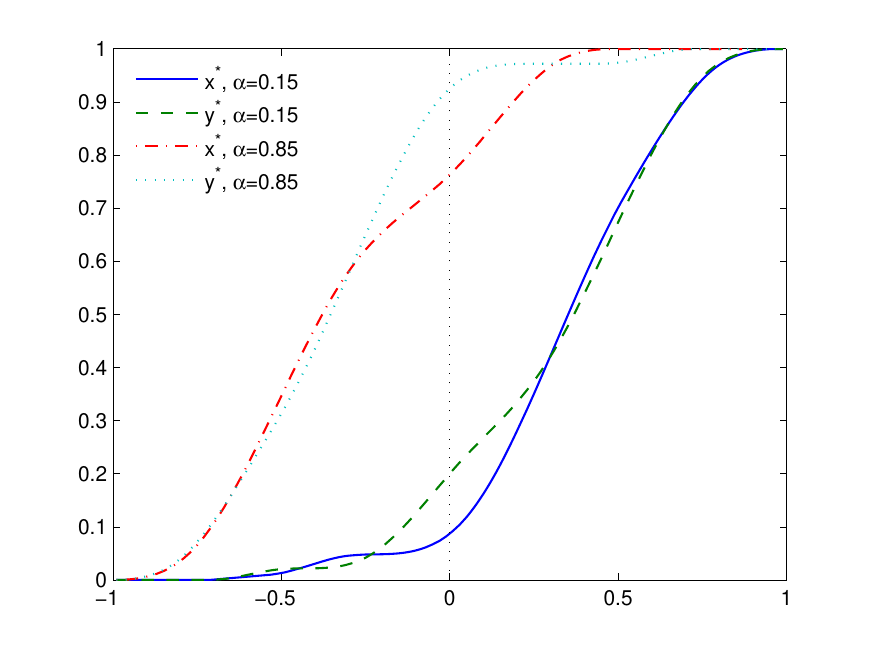}\\
(a)&(b)\\
\epsfxsize=220pt\epsfysize=160pt\epsffile{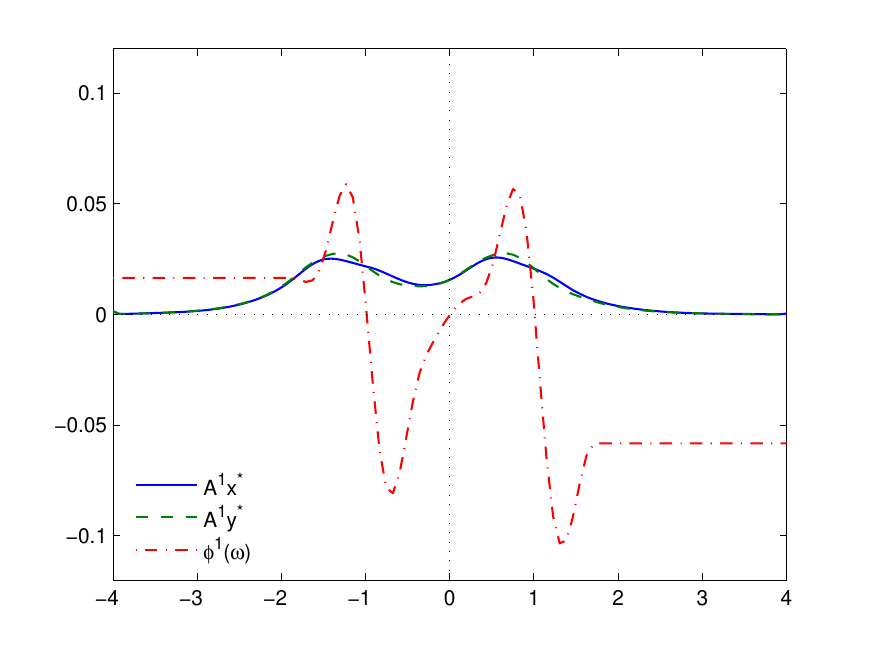}
&\epsfxsize=220pt\epsfysize=160pt\epsffile{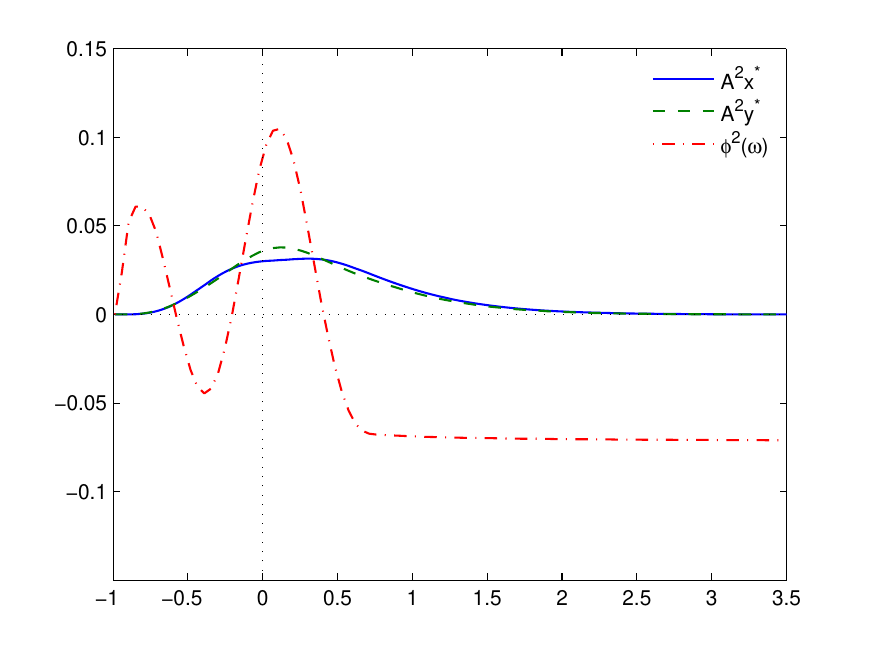}\\
(c)&(d)
\end{array}
$$
\caption{\label{fig:dom1}  Deconvolution experiment, $K_\ell=1000$, $k=1,2$, $\epsilon=0.05$. In the upper row: (a) resolution of the simple test as a function of $t\in [-1,1]$; (b) c.d.f. of the ``difficult to test'' distributions $x^*$ and $y^*$, corresponding optimal solutions to $(F_{g,\alpha}[\rho])$ for $g=g(0)$ (testing hypotheses about $F(0)$). Bottom row: convolution images of optimal solutions to $(F_{g,\alpha}[\rho])$, $\alpha=.85$ and $g=g(0)$, and corresponding detector $\phi$: (c) convolution with mixed Laplace distribution, (d) convolution with $\Gamma(\cdot)$ distribution.}
\end{figure}

\begin{figure}
$$
\begin{array}{cc}
&\epsfxsize=220pt\epsfysize=160pt\epsffile{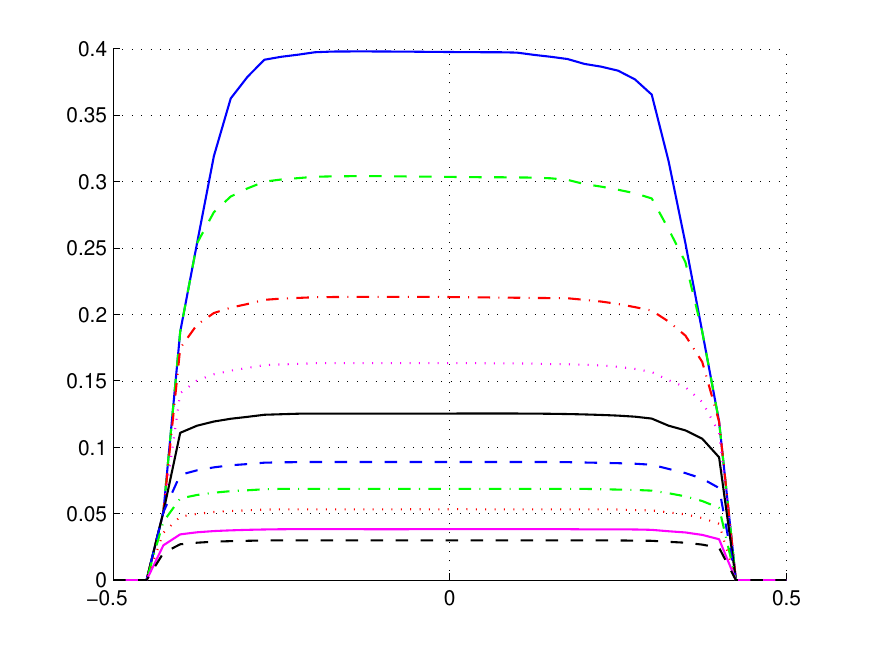}
\epsfxsize=220pt\epsfysize=160pt\epsffile{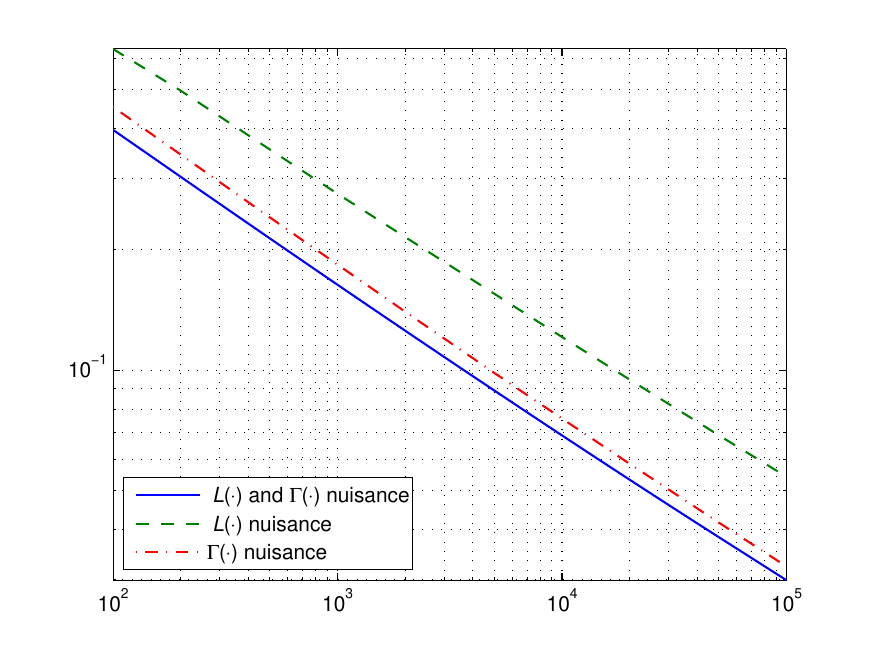}
\end{array}
$$
\caption{\label{fig:dom2} Deconvolution experiment, $\epsilon=0.05$, $\alpha=0.5$; $K_\ell=[100,200,500,1000,...,100\,000]$, $\ell=1,2$. On the left:  resolution of the simple test as a function of $t\in [-1,1]$ for different $K^\ell$, mixed Laplace and $\Gamma(\cdot)$ distributions of the observation noise;
on the right:  resolution at $t=0$ as a function of $K^\ell$; the test resolution clearly exhibits $C\,K^{-1/3}$ behavior.}
\end{figure}

On figures \ref{fig:dom1} and \ref{fig:dom2} we present details of the test in the deconvolution model with $L=2$ observers. Each observer acquires $K_\ell$ noisy observations $\zeta^\ell_k$, $k=1,...,K_\ell$. The distribution of the nuisance is mixed Laplace for the first observer and $\Gamma(0,2,1/2/\sqrt(2))$  for the second observer. The discretized model has the following parameters:
the observation spaces $\Omega_\ell=\bbr$, $\ell=1,2$ of each of 2 $K_\ell$--repeated observation schemes were split into $m_\ell=102$ ``bins'': we put $b_1^\ell=-1+q_{\eta^\ell}([K^\ell]^{-1})$ and $b_{100}^\ell=1+q_{\eta^\ell}(1-[K^\ell]^{-1})$, and split the interval $(b_1^\ell,b^\ell_{100}]$ into 100 equal length bins; then we add two  bins $(-\infty,b_1^\ell]$ and $(b_{100}^\ell,\infty)$.

\begin{figure}
$$
\begin{array}{cc}
\epsfxsize=220pt\epsfysize=160pt\epsffile{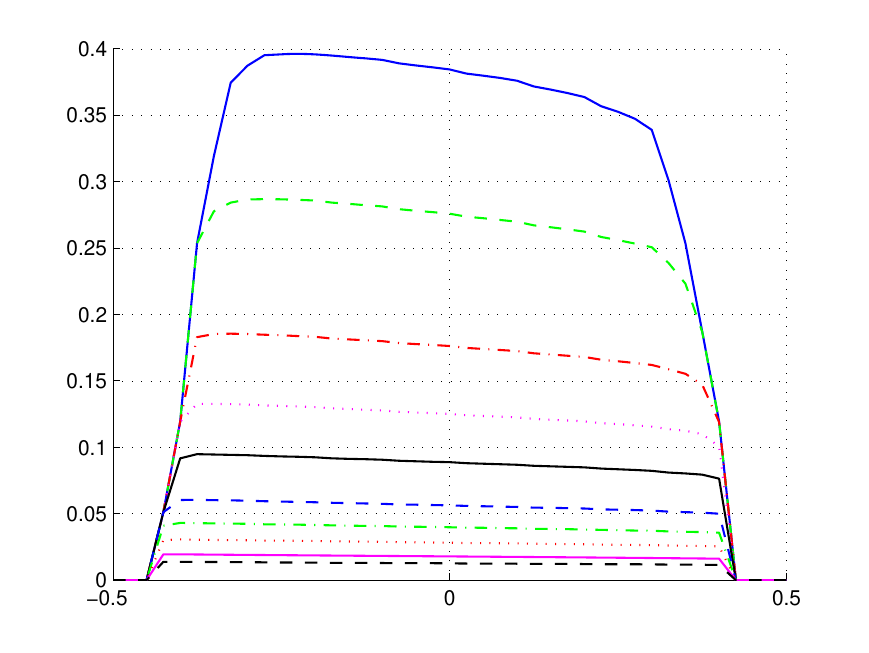}
&\epsfxsize=220pt\epsfysize=160pt\epsffile{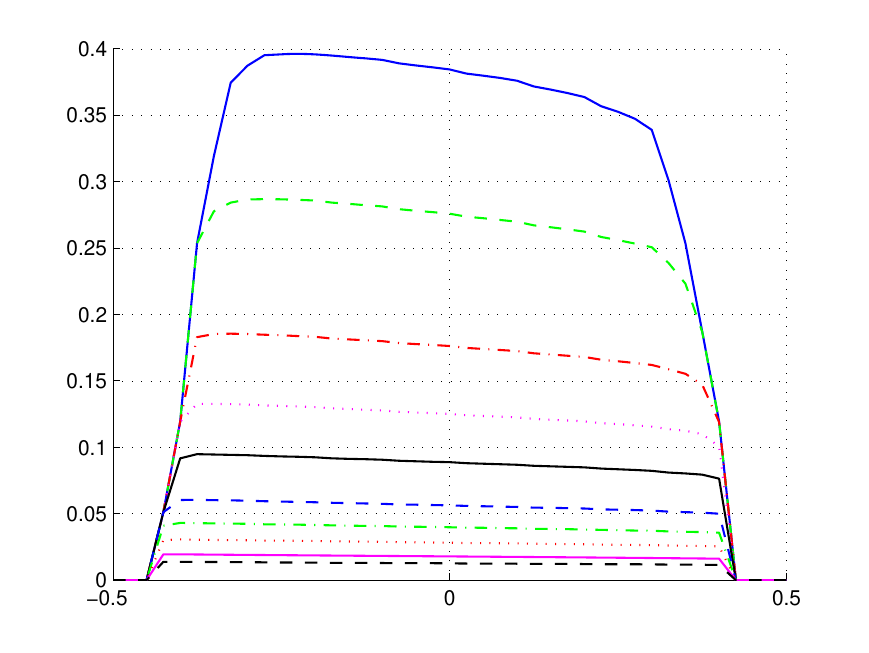}\\
(a)&(b)\\
\epsfxsize=220pt\epsfysize=160pt\epsffile{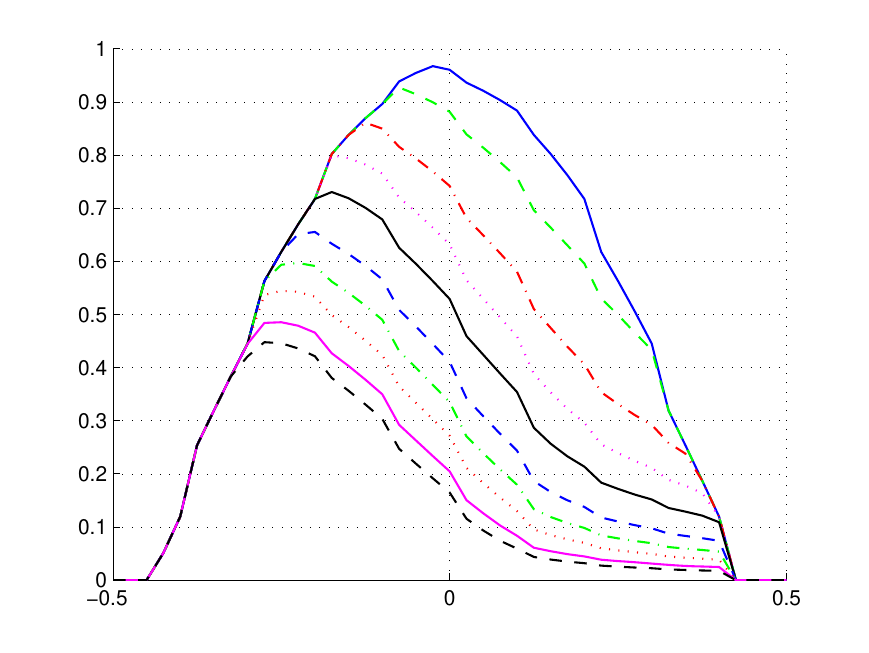}
&\epsfxsize=220pt\epsfysize=160pt\epsffile{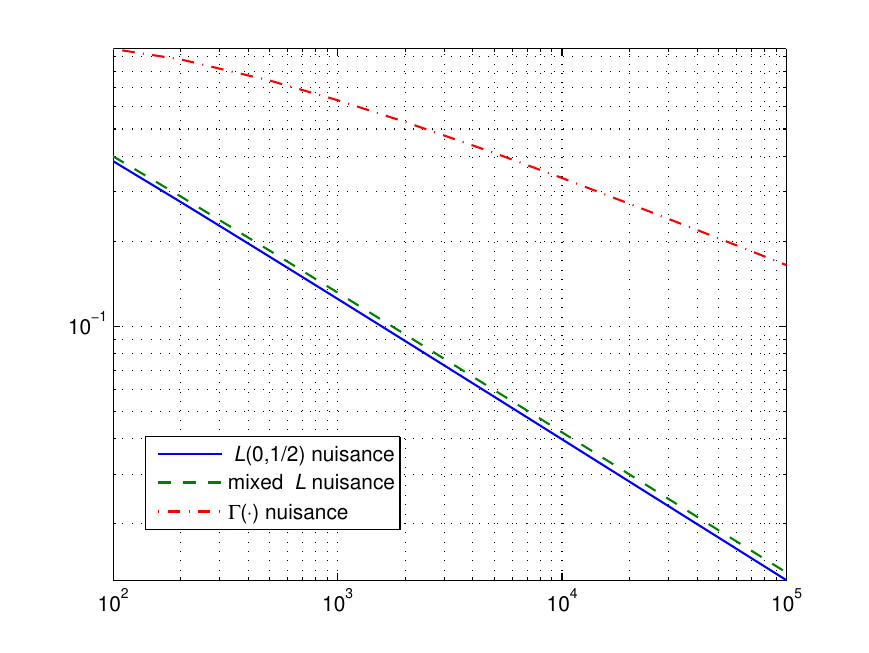}\\
(c)&(d)
\end{array}
$$
\caption{\label{fig:dom3}  Trimmed observation experiment,  resolution of the simple test for different $K$, $\epsilon=0.05$, $\alpha=0.5$; $K=[100,200,500,1000,...,100\,000]$. Plot (a): resolution of the test as a function of $t\in [-1,1]$, $L(0,\half)$ nuisance; plot (b) same for mixed Laplace nuisance; plot (c): resolution of the test with $\Gamma(\cdot)$ nuisance distribution. On plot (d): resolution at $t=0$ as a function of sample size $K$. While the test resolution exhibits $C\,K^{-1/3}$ behavior in the case of Laplace an mixed Laplace nuisance, convergence is slow (if any) in the case of $\Gamma(\cdot)$ nuisance distribution.}
\end{figure}

On figure \ref{fig:dom3} we present simulation results for the experiments with trimmed observations. Here $L=1$, the observations are $\omega_k=\max[\xi_k,\eta_k]$, $1\leq k\leq K$, with the $\LL(0,\half)$ nuisances $\eta_k$. The partition of the support $[-1,1]$ of $\xi$ is the same as in the deconvolution experiments, and the observation domain was split into $m=101$ bins -- 100 equal length bins over the segment $[-1,1]$ and the bin $(1,\infty)$.

\paragraph{Quantifying conservatism.} When building the test $\widehat{T}_\rho$ deciding on the hypotheses $H_\imath[\rho]$, $\imath=1,2$ (see $(\D_{g,\alpha}[\rho])$) via $K$ observations $\omega^K=(\omega_1,...,\omega_K)$, we get, as a byproduct, two probability distributions $x_\rho\in\F$, $y_\rho\in\F$, of the latent random variable $\xi$, see (\ref{testwidehatTrho}). These distributions give rise to two simple hypotheses, $\overline{H}_1$, $\overline{H}_2$, on the distribution  of observation $\omega^K$, stating that these observations come from the distribution $x_\rho$, resp., $y_\rho$, of the latent variable. The risk of any test deciding on the two simple hypotheses $\overline{H}_1$, $\overline{H}_2$, the observation being $\omega^K$, is lower-bounded by the quantity $\widehat{\epsilon}[K]=\sum_{\omega^K} \min[p_1^K(\omega^K),p_2^K(\omega^K)]$, where $p_i^K(\omega^K)$ is the probability to get an observation $\omega^K$ under hypothesis $\overline{H}_i$, $i=1,2$. The quantity $\widehat{\epsilon}[K]$, which can be estimated by Monte-Carlo simulation, by its origin is a lower bound on the risk of a whatever test deciding, via $\omega^K$,  on the composite ``hypotheses of interest'' $H_\imath[\rho]$, $\imath=1,2$. We can compare this lower risk bound with the upper bound $\epsilon[K]=\exp\{\Opt[\rho]\}$ on the risk of the test $\widehat{T}_\rho$, see $(F_{g,\alpha}[\rho])$, and thus quantify the conservatism of the latter test.
The setup of the related experiments was completely similar to the one in the just reported experiments, with the Laplace distribution $\LL(0,1/2)$ of the nuisance and with $n=500$ and $m=1002$ bins in the supports of $\xi$ and of $\omega$, respectively. We used $t=0$, $\alpha=0.5$, and $2\times 10^6$ Monte-Carlo simulations to estimate $\widehat{\epsilon}[K]$. In our experiments, given a number of observations $K$ and a prescribed risk level $\epsilon\in\{0.1,0.01,0.001,0.0001\}$, the parameter $\rho$ of the test $\widehat{T}_\rho$  was adjusted to ensure $\epsilon[K]=\epsilon$; specifically, we set $\rho=\rho[\epsilon]$, see (\ref{rho_A}).  The results are presented in table \ref{tab:lower31}. \par
Recall that by Proposition \ref{propnonstI} we have $\epsilon[K']\leq(\epsilon[K])^{{K'/K}}$ when $K'\geq K$, so that the ratios  $r[k]=\ln(\widehat{\epsilon}[K])/\ln(\epsilon[K])$ presented in the table upper-bound the nonoptimality of $\widehat{T}_\rho$ in terms of the number of observations required to achieve the risk $\widehat{\epsilon}[K]$: for the ``ideal'' test, at least $K$ observations are required to attain this risk, and for the test $\widehat{T}_\rho$ -- at most $\rfloor r[k]K\lfloor$ observations are enough. The data in table \ref{tab:lower31} show that the ratios $r[K]$ in our experiments newer exceeds 1.82 and steadily decrease when $\epsilon[K]$ decreases.
\begin{table}
\begin{center}
{\small \begin{tabular}{|c|c|c|c|c|c|c|c|}
  \hline
\backslashbox {$\epsilon$}{$K$}&200&500&1000&2000&5000&10000&20000
\\\hline
\multirow{2}{10mm}{1.0e-1}&1.5e-2&1.5e-2&1.7e-2&1.6e-2&1.6e-2&1.6e-2&1.5e-2\\
&  1.82&1.82&1.78&1.80&1.80&1.80&1.82\\
\hline
\multirow{2}{10mm}{1.0e-2}&1.3e-3&1.2e-3&1.2e-3&1.2e-3&1.2e-3&1.2e-3&1.2e-3\\
&1.45&1.46&1.46&1.46&1.46&1.45&1.46\\
\hline
\multirow{2}{10mm}{1.0e-3}&1.0e-4&0.9e-4&1.1e-4&1.1e-4&1.1e-4&0.9e-4&1.1e-4\\
& 1.33&1.35&1.32&1.32&1.32&1.34&1.32\\
\hline
\multirow{2}{10mm}{1.0e-4}&1.1e-5&0.9e-5&1.0e-5&0.9e-5&1.1e-5&0.7e-5&0.9e-5\\
&1.24&1.26&1.25&1.26&1.24&1.29&1.26\\
\hline
\end{tabular}
}
\caption{\label{tab:lower31} Quantifying conservatism of $\widehat{T}_\rho$ in Deconvolution experiment; in a cell: top -- $\widehat{\epsilon}[K]$, bottom -- the ratio ${\ln \widehat{\varepsilon}[K]\over \ln{\epsilon}[K]}$.}
\end{center}
\end{table}

\subsection{Testing hypotheses on Markov chains}
In this section, we present some applications of our approach to Markov chain related hypotheses testing.
For a positive integer $n$, let $\Delta_n=\{x\in\bbr^n_+:\sum_ix_i=1\}$, and  $\cS_n$ be the set of all $n\times n$ stochastic matrices..
\subsubsection{Deciding on two simple hypotheses}\label{simpleMarkov}
\paragraph{Situation.} The simplest setting of the Markov chain related hypotheses testing is as follows. We are given two $n\times n$ stochastic matrices $S^1$ and $S^2$ with positive entries, specifying two hypotheses on an $n$-state Markov chain. Both hypotheses state that the probability distribution of the initial (at time 0) state $\iota_0$ of the chain is a vector from some convex compact set $X\subset\rint\Delta_n$; in addition hypothesis $H_1$ ($H_2$)  states that the transition matrix of the chain is $S^1$($S^2$). We observe on a given time horizon $K$ a realization $\iota_0,\iota_1,...,\iota_K$ of the trajectory of the chain and want to decide on the hypotheses.
\paragraph{Construction and result.} With transition matrix fixed, the distribution of chain's trajectory on a fixed time horizon depends linearly on the distribution of the initial state. Consequently, our decision problem is to distinguish between two convex sets  of probability distributions on the finite set of all possible  chain trajectories from time 0 to time $K$ inclusively.
 According to the Discrete case version of our results, a nearly optimal test  is as follows: we solve the optimization problem
\begin{equation}\label{asgivenby}
\varepsilon_\star=\max_{p,q\in X} \sum_{1\leq \iota_0,\iota_1,...,\iota_K\leq n} \sqrt{\left[p_{\iota_0}S^1_{\iota_1\iota_0}S^1_{\iota_2\iota_1}...S^1_{\iota_K\iota_{K-1}}\right]
\left[q_{\iota_0}S^2_{\iota_1\iota_0}S^2_{\iota_2\iota_1}...S^2_{\iota_K\iota_{K-1}}\right]};
\end{equation}
denoting the optimal solution $(p_*,q_*)$ and setting
$$
\phi(\iota_0,...,\iota_K)={1\over 2}\ln\left({p_{\iota_0}S^1_{\iota_1\iota_0}S^1_{\iota_2\iota_1}...S^1_{\iota_K\iota_{K-1}}\over
q_{\iota_0}S^2_{\iota_1\iota_0}S^2_{\iota_2\iota_1}...S^2_{\iota_K\iota_{K-1}}}\right),
$$
the near-optimal test, the observed trajectory being $\iota^K=(\iota_0,...,\iota_K)$, accepts $H_1$ when $\phi(\iota^K) \geq 0$, and accepts $H_2$ otherwise. The risk of this test is upper-bounded by $\varepsilon_\star$ given by (\ref{asgivenby}). \par
Optimization problem (\ref{asgivenby}) clearly is convex and solvable, and whenever $(p,q)$ is feasible for the problem, so is $(q,p)$, the values of the objective at these two solutions being the same. As a result, there exists an optimal solution $(p_*,q_*)$ with $p_*=q_*$. The test $\phi$ associated with such a solution is completely independent of $p_*$ and is just the plain likelihood ratio test:
$$
\phi(\iota^K=(\iota_0,...,\iota_K))={1\over 2}\sum_{\tau=1}^K\ln\left({S^1_{\iota_\tau\iota_{\tau-1}}\over S^2_{\iota_\tau\iota_{\tau-1}}}\right).
$$
The (upper bound on the) risk of this test is immediately given by (\ref{asgivenby}):
$$
\varepsilon_\star=\max_{p\in X}\sum_{j=1}^m\left(\sum_{i=1}^m
(S^1_{ij}S^2_{ij})^{t/2}\right)p_j.
$$
\paragraph{Numerical illustration.} Consider a queuing system ($M/M/s/s+b$) with $s$ identical servers, with services times following exponential distribution $\E(\mu)$ with parameter $\mu$, and a common buffer of capacity $b$. The input stream of customers is Poisson process with rate $\lambda$. Upon arrival, a customer either starts to be served, if there is a free server, or joins the buffer, if all servers are busy and there are less than $b$ customers in the buffer, or leaves the system, if all servers are busy and there are $b$ waiting customers in the buffer. The system is observed at time instances $0,1,...,K$, and we want to distinguish between two systems differing only in the value of $\mu$, which is $\mu_1$ for the first, and $\mu_2$ for the second system. The observations form a Markov chain with $n=s+b+1$ states,
 a state $j\in\{1,...,n\}$ at time $t=1,2,...$ meaning that at this time there are $s(j):=\min[j-1,s]$ busy servers and {$j-s(j)-1$} customers in the buffer. Under hypothesis $H_\chi$, $\chi=1,2$, the  transition matrix of the chain is $S^\chi=\exp\{L^\chi\}$, where $L^\chi=L(\lambda,\mu_\chi)$ is a 3-diagonal {\em transition rate matrix} with zero column sums and $[L^\chi]_{j-1,j}=s(j)\mu_\chi$, $[L^\chi]_{j+1,j}=\lambda$. In table \ref{tabletwochains}, we present a sample of (the smallest) observation times $K$ ensuring that the upper bound $\varepsilon_\star$ on the risk of the simple test developed in this section is $\leq0.01$. We restrict ourselves to the case when distribution of the initial state is not subject to any restrictions, that is, $X=\Delta_{s+b+1}$.
\begin{table}
{\footnotesize
$$
\begin{array}{||c|c|c||c|c|c||c|c|c||c|c|c||c|c|c||c|c|c||}
\hline\hline
\multicolumn{6}{||c||}{\lambda=50}&\multicolumn{6}{c||}{\lambda=100}&\multicolumn{6}{c||}{\lambda=200}\\
\hline
\mu_1&\mu_2&K&\mu_1&\mu_2&K&\mu_1&\mu_2&K&\mu_1&\mu_2&K&\mu_1&\mu_2&K&\mu_1&\mu_2&K\\
\hline\hline
1.00&0.90&144&1.00&1.11&146&
1.00&0.90&91&1.00&1.11&74&
1.00&0.90&1929&1.00&1.11&1404\\
\hline
1.00&0.75&21&1.00&1.33&21&
1.00&0.75&19&1.00&1.33&11&
1.00&0.75&326&1.00&1.33&133\\
\hline
1.00&0.50&6&1.00&2.00&5&
1.00&0.50&8&1.00&2.00&3&
1.00&0.50&86&1.00&2.00&7\\
\hline\hline
\end{array}
$$}
\caption{\label{tabletwochains} Deciding with risk $\varepsilon_*=0.01$ on two simple hypotheses on the parameter $\mu$ of a queuing system with $s=100$, $b=20$.}
\end{table}

\subsubsection{Deciding on two composite hypotheses}\label{twocompositeMarkov}
In the previous example, we dealt with two simple hypotheses on a Markov chain with fully observable trajectory. Now consider the case of two composite hypotheses and indirect observations of state transitions.\footnote{One problem of testing specific composite hypotheses about Markov chains has been studied in \cite{Birge1983M} using a closely related approach. The techniques we discuss here are different and clearly aimed at numerical treatment of the problem.}
 More specifically, we intend to consider the case when a ``composite hypothesis'' specifies a set in $\cS_n$ containing the transition matrix of the chain we are observing, and ``indirectness of observations'' means that instead of observing consecutive states of the chain trajectory, we are observing some encodings of these states (e.g., in the simplest case, the state space of the chain is split into non-overlapping subsets -- {\sl bins}, and our observations are the bins to which the consecutive states of the chain belong).
\paragraph{Preliminaries.} Probability distribution $P_t$ of the trajectories, on time horizon $t$, of a Markov chain depends nonlinearly on the transition matrix of the chain. As a result, to utilize our convexity-based approach, we need to work with composite hypotheses of ``favorable structure,'' meaning that the family $\P_t$ of distributions $P_t$ associated with transition matrices allowed by the hypothesis admits a reasonable convex approximation. We start with specifying the main ingredient of such ``favorable structure.''  \par
Let $K_1,...,K_n$ be closed cones, all different from $\{0\}$, contained in $\bbr^n_+$. The collection $K^n=\{K_1,...,K_n\}$ gives rise to the following two entities:
\begin{itemize}
\item The set of stochastic matrices \[
\cS=\{S=[S_{ij}]_{i,j=1}^n\in\bbr^{n\times n}: \Col_j[S]\in K_j,\sum_iS_{ij}=1,\,j=1,...,n\}
\]
(from now on, $\Col_j[S]$ is the $j$-th column of $S$);
\item The convex set
\[
\cP=\{P=[P_{ij}]_{i,j=1}^n\in\bbr^{n\times n}: \Col_j[P]\in K_j,\,1\leq j\leq n, \sum_{i,j}P_{ij}=1\}.
\]
\end{itemize}
One has\footnote{Indeed, for $S\in\cS$, $x\in\Delta_n$ the matrix $P$ given by $\Col_j[P]=x_j\Col_j[S]$, $1\leq j\leq n$,  clearly belongs to $\cP$. Vice versa, if $P\in\cP$, then, setting $x_j=\sum_iP_{ij}$ and specifying the $j$-th column of $S$ as $\Col_j[P]/x_j$ when $x_j\neq0$ and as a whatever vector from $K_j\cap\Delta_n$ when $x_j=0$, we get $S\in\cS$, $x\in\Delta_n$ and $\Col_j[P]=x_j\Col_j[S]$ for all $j$. 
}
\begin{equation}\label{eq1957}
\cP=\{P=[P_{ij}]_{i,j=1}^n: \exists (S\in\cS,x\in\Delta_n):\Col_j[P]=x_j\Col_j[S],\,j=1,...,n\}.
\end{equation}
As a result, in a pair $(S,x)$ associated with $P\in\cP$ according to {\rm (\ref{eq1957})}, $x$ is uniquely defined by $P$:
$$
x_j=\sum_{i}P_{ij},\,1\leq j\leq n;
$$
besides this, for every $j$ such that $\sum_iP_{ij}>0$, $\Col_j[S]$ is the probabilistic normalization of $\Col_j[P]$.
\paragraph{Remark.}
{\rm The role played by the just defined entities in our context stems from the following immediate observation: consider a Markov chain with transition matrix $S$ from $\cS$, and let $x\in\Delta_n$ be the distribution of the state $\iota_{\tau-1}$ of this chain at time $\tau-1$. Denoting by $\iota_\tau$ the state of the chain at time $\tau$, the distribution of the state transition $(\iota_{\tau-1},\iota_\tau)$ clearly is
$$
p_{ij}=S_{ij}x_j,\,1\leq i,j\leq n.
$$
According to \rf{eq1957}, $\cP$ is nothing but the convex hull of all distributions of this type stemming from different $x\in\Delta_n$ and $S\in\cS$.}
\paragraph{Situation.} Assume that for $\chi=1,2$ we are given
\begin{itemize}
\item collection of cones $K^{n_\chi}_\chi=\{K_1^\chi,...,K_{n_\chi}^\chi\}$ of the type described in the preliminaries. This collection, as explained above, specifies a set $\cS_\chi$ of stochastic $n_\chi\times n_\chi$ matrices and a set $\cP_\chi$ of $n_\chi\times n_\chi$ matrices with nonnegative entries summing up to 1.
\item $m\times n_\chi^2$ ``observation matrix'' $A_\chi$ with {\em positive entries} and  unit column sums. We think of the $n_\chi^2$ columns of $A_\chi$ as being indexed by the pairs $(i,j)$, $1\leq i,j\leq n_\chi$.
\end{itemize}
The outlined data specify, for $\chi=1,2,$
\begin{itemize}
\item the family $\cM_\chi$ of Markov chains. Chains from $\cM_\chi$ have $n_\chi$ states, and their transition matrices belong to $\cS_\chi$;
\item observation scheme for transitions of a chain from $\cM_\chi$. Specifically, observation $\omega_\tau$ of the transition $\iota_{\tau-1}\to\iota_\tau$ takes values in $\{1,2,...,m\}$, and its conditional, the past of chain's state trajectory being given, distribution is the column $\Col_{(\iota_{\tau-1},\iota_\tau)}[A_\chi]$ of $A_\chi$.
\end{itemize}
Now assume that ``in the nature'' there exist two Markov chains, indexed by $\chi=1,2$, with $n_\chi$ states and transition matrices $S_\chi$, such that chain $\chi$ belongs to $\cM_\chi$, and we observe one of these two chains as explained above, so that, independently of $\chi$, our observation $\omega_t$ at time $t$ takes values in $\{1,...,m\}$. Given observation $\omega^K=(\omega_1,...,\omega_K)$, we want to decide on the hypotheses $H_\chi$, $\chi=1,2$, where $H_\chi$ states that the chain we are observing is chain $\chi$.
\paragraph{Construction and result.}
We can approach our goal as follows. Every $P\in\cP_\chi$ is a nonnegative $n_\chi\times n_\chi$ matrix with unit sum of entries and as such can be thought of as a probability distribution on $\cI_\chi=\{(i,j):1\leq i,j\leq n_\chi\}$. Matrix $A_\chi$ naturally associates with such a distribution a probability distribution $\cA_\chi(P)$ on $\{1,...,m\}$:
$$
\cA_\chi(P)=\sum_{i,j=1}^{n_\chi} P_{ij}\Col_{(i,j)}(A_\chi).
$$
Note that the mapping $P\mapsto\cA_\chi(P)$ is linear.
\par
Let us define the convex compact subsets $X^\chi$ of  the probabilistic simplex $\Delta_m$  by the relation
$$
X^\chi=\{p\in\Delta_m: \exists P\in\cP_\chi: p=\cA_\chi(P)\},\,\,\chi=1,2.
$$
By the above remark,
 \begin{quote}
 (!) {\sl For a chain from $\cM_\chi$ and every time instant $\tau\geq1$, the conditional, given chain's trajectory prior to instant $\tau-1$, distribution of the state transition $(\iota_{\tau-1},\iota_\tau)$ belongs to $\cP_\chi$, and, consequently, the conditional, by the same condition, distribution of the observation $\omega_\tau$ belongs to $X^\chi$.}
 \end{quote}
Note that $X^\chi\subset\rint\Delta_m$ due to entrywise positivity of $A_\chi$.
\par
For $t=1,2,...$, let $\jmath_{t,1}$, $\jmath_{t,2}$  be the states of chain 1 and chain 2 at time $t$, let $\zeta_{t,\chi}=(\jmath_{t,\chi},\jmath_{t-1,\chi})$, $\chi=1,2$,  and let $X_t=X^1$, $Y_t=X^2$. With this setup,
 we arrive at the situation considered in Proposition \ref{propnonstationarynew}: for $\chi=1,2$, under hypothesis $H_\chi$ $\omega_t$ is a deterministic function of $\zeta^t_\chi=(\zeta_{1,\chi},...,\zeta_{t,\chi})$, the conditional, given $\zeta^{t-1}_\chi$, distribution of $\omega_t$ depends deterministically on $\zeta^{t-1}_\chi$ and, by (!), belongs to $X^\chi$. Hence, Proposition \ref{propnonstationarynew} implies
\begin{proposition}\label{proptransitions} In the situation and under assumptions of this section, let the sets $X^1$, $X^2$ do not intersect. Let $p^*_1$, $p^*_2$, form the optimal solution to the problem
\begin{equation}\label{Markovvarepsilon}
\varepsilon_\star=\max_{p_1,p_2}\left\{\sum_{\omega=1}^m \sqrt{[p_1]_\omega [p_2]_\omega}:\;p_1\in X^1,\;p_2\in X^2,\right\},
\end{equation}
and let
$$
\phi(\omega)={1\over 2}\ln\left({[p^*_1]_{\omega}\over [p^*_2]_\omega}\right).
$$
Then the risk of the test which, given observations $\omega_1,...,\omega_K$, accepts $H_2$ when $\sum_{\tau=1}^K\phi(\omega_\tau) \geq0$ and accepts $H_2$ otherwise, is at most $\varepsilon_\star^K$.
\end{proposition}
\paragraph{Remark.} By inspecting the proof, Proposition \ref{proptransitions} remains valid in the situation where $\cM_\chi$ are families of {\sl non-stationary} Markov chains
 with $n_\chi$ states $1,...,n_\chi$. In such a chain, for every $\tau>0$, the conditional, given the trajectory
$\iota_0,...,\iota_{\tau-1}$ of the chain from time 0 to time $\tau-1$, distribution of state $\iota_\tau$ at time $\tau$ is selected, in a non-anticipative fashion, from the set $K_{\iota_{\tau-1}}^\chi\cap\Delta_n$.
\paragraph{Numerical illustration: random walk.} Consider a toy example where the Markov chains $\cM_\chi$, $\chi=1,2$, represent a random walk along $n=16$-element grid on the unit circle; thus, each chain has 16 states. The ``nominal'' transition matrices $S^n_\chi$ correspond to the walk where one stays in the current position with probability $1-2p_\chi$ and jumps to a neighbouring position with probability $2p_\chi$, with equal probabilities to move clock- and counter-clockwise; in our experiment, $p_1=0.2$ and $p_2=0.4$. The actual transition matrix $S_\chi$ of chain $\cM_\chi$ is allowed to belong to the ``uncertainty set''
\[
\cU_\chi=\{S_\chi\in\cS_n: (1-\rho)S_\chi^n\leq S_\chi\leq(1+\rho)S_\chi^n\},
 \]
 where the inequalities are entrywise. In other words, the cones $K_j^\chi$, $j=1,2,...,n$,
 are the conic hulls of the sets
 \[
 \{q\in \Delta_n:(1-\rho)\Col_j[S^n_\chi]\leq q\leq (1+\rho)\Col_j[S^n_\chi]\}.
  \]
  In our experiments, we used $\rho=0.1$.
\par
We  have considered two observation schemes: ``direct observations'', where we observe the positions of the walker  at times 0,1,..., and ``indirect observations,'' where the 16 potential positions are split into 8 ``bins,'' two states per bin, and what we see at time instant $t$ is the bin to which $t$-th position of the walker belongs. In the latter case we used a random partition of the states into the bins which was common for the chains $\cM_1$ and $\cM_2$ (i.e., in our experiments the ``observation matrices''  $A_1$ and $A_2$ always coincided with each other).\par
The results of a typical experiment are presented in  table \ref{tabwalk}. For each of our two observation schemes, we start with observation time which, according to Proposition \ref{proptransitions}, guarantees the risk $\epsilon=0.01$, and then decrease the observation time to see how the performance of the test deteriorates. In different simulations, we used different transition matrices allowed by the corresponding hypotheses, including the ``critical'' ones -- those associated with the optimal solution to (\ref{Markovvarepsilon}).
Evaluating the results of the experiment
is not easy -- in the first place, it is unclear what could be a natural ``benchmark'' to be compared to, especially when the observations are indirect.
In the case of direct observations we have considered as a contender the likelihood ratio test (see section \ref{simpleMarkov}) straightforwardly adjusted to the  uncertainty in the transition matrix.\footnote{Specifically, given the chain trajectory $\iota_0,...,\iota_t$, we can easily compute the maximal and the minimal values, $\psi_{\max}$ and $\psi_{\min}$, of the logarithm of likelihood ratio as allowed by our uncertainties in the transition matrices. Namely,   $\psi_{\max}=\max_{\{S_{\tau,1},S_{\tau,2}\}_{\tau=1}^t}\sum_{\tau=1}^t\ln([S_{\tau,1}]_{\jmath_\tau,\jmath_{\tau-1}}/[S_{\tau,2}]_{\jmath_\tau,\jmath_{\tau-1}})$, where $S_{\tau,\chi}$ run through the uncertainty sets associated with hypotheses $H_\chi$, $\chi=1,2$; $\psi_{\min}$ is defined similarly, with $\max_{\{S_{\tau,1},S_{\tau,2}\}_{\tau=1}^t}$ replaced with $\min_{\{S_{\tau,1},S_{\tau,2}\}_{\tau=1}^t}$.  We accept $H_1$ when a randomly selected point in $[\psi_{\min},\psi_{\max}]$ turns out to be nonnegative, and accept $H_2$ otherwise.} Such test turns out to be essentially less precise than the test presented in Proposition \ref{proptransitions}; e.g., in the experiment reported in column A of table \ref{tabwalk}, with observation time 71 the risks of the adjusted likelihood  test were as large as $0.01/0.06$.
\begin{table}[h]
\centering
\subfloat[$\varepsilon_\star=0.9368$]{
\begin{tabular}{||c|c|c|c||}
\hline
$t$&$\varepsilon_\star^t$&Risk(T)&Risk(ML)\\
\hline
\hline
71&0.0097&0.0004/0.0008&0.0094/0.0551\\
\hline
48&0.0436&0.0038/0.0018&0.0192/0.0798\\
\hline
32&0.1239&0.0226/0.0118&0.0390/0.1426\\
\hline
21&0.2540&0.0230/0.0610&0.0620/0.1903\\
\hline
14&0.4011&0.0870/0.0508&0.1008/0.2470\\
\hline
10&0.5207&0.0780/0.1412&0.1268/0.2649\\
\hline
7&0.6333&0.1184/0.1688&0.1824/0.3368\\
\hline
5&0.7216&0.1040/0.2682&0.2190/0.2792\\
\hline
3&0.8222&0.3780/0.1166&0.3000/0.4027\\
\hline
2&0.8777&0.1814/0.3780&0.1814/0.3780\\
\hline
1&0.9368&0.4230/0.2064&0.4230/0.2064\\
\hline
\end{tabular}
}
\qquad%
\subfloat[$\varepsilon_\star=0.9880$]{
%
\begin{tabular}{||c|c|c||}
\hline
$t$&$\varepsilon_\star^t$&Risk(T)\\
\hline
\hline
381&0.0099&0.0000/0.0000\\
\hline
254&0.0462&0.0000/0.0000\\
\hline
170&0.1277&0.0000/0.0002\\
\hline
113&0.2546&0.0002/0.0008\\
\hline
76&0.3982&0.0002/0.0054\\
\hline
51&0.5393&0.0022/0.0168\\
\hline
34&0.6626&0.0086/0.0412\\
\hline
23&0.7569&0.0210/0.0758\\
\hline
15&0.8339&0.0540/0.1018\\
\hline
10&0.8860&0.0872/0.1530\\
\hline
7&0.9187&0.1420/0.1790\\
\hline
5&0.9413&0.1386/0.2878\\
\hline
3&0.9643&0.2812/0.2638\\
\hline
2&0.9761&0.2078/0.3824\\
\hline
1&0.9880&0.3816/0.2546\\
\hline
\end{tabular}
}
\caption{\label{tabwalk} Random walk. (a) - direct observations; (b) - indirect observations. In the table:\newline
$t$: observation time; $\varepsilon_\star^t$ and Risk(T): theoretical upper bound on the risk of the test from Proposition \ref{proptransitions}, and empirical risk of the test;
Risk(ML): empirical risk of the likelihood ratio test adjusted for uncertainty in transition probabilities.
$\epsilon_1/\epsilon_2$ in ``risk'' columns: empirical, over 5000 simulations, probabilities to reject hypothesis $H_1$  ($\epsilon_1$)  and $H_2$ $(\epsilon_2)$ when the hypothesis is true.
Partition of 16 states of the walk into 8 bins in the reported experiment is $\{1,8\}$, $\{4,6\}$, $\{5,7\}$, $\{9,11\}$, $\{3,19\}$, $\{2,15\}$, $\{12,16\}$, $\{13,14\}$.
}
\end{table}
\subsubsection{Two composite hypotheses revisited}\label{twocompositeMarkovrevisited}
In the situation of section \ref{twocompositeMarkov} (perhaps, indirect) observations of {\sl transitions} of a Markov chain were available. We are about to consider the model in which we are only allowed  to observe how frequently  the chain visited different (groups of) states on a given time horizon, but do not use information in which order these states were visited.
\paragraph{Preliminaries.}
  For $Q\in \cS_n$ and $\rho\geq0$, let
 $$
 \cS_n(Q,\rho)=\{S\in\cS_n: \|S-Q\|_{1,1}\leq\rho\},
  $$
 where for a $p\times q$ matrix $C$
$$
\|C\|_{1,1}=\max_{1\leq j\leq q}\|\Col_j[C]\|_1
$$
is the norm of the mapping $u\mapsto Cu:\bbr^q\times \bbr^p$ induced by the norms $\|\cdot\|_1$ on the argument and the image spaces.
\paragraph{Situation} we consider here is as follows. ``In the nature'' there exist two Markov chains, indexed by $\chi=1,2$. Chain $\chi$  has $n_\chi$ states and transition matrix $S_\chi$. Same as in section \ref{twocompositeMarkov}, we do not observe the states exactly, and our observation scheme is as follows.
For $\chi=1,2$, we are given  $m\times n_\chi$ matrices $A_\chi$ with positive entries and all column sums equal to 1. When observing chain $\chi$, our observation
  $\eta_\tau$ at time $\tau$ takes values $1,...,m$, and the conditional, given the trajectory of the chain since time $0$ to time $\tau$ inclusively, distribution of $\eta_\tau$ is the $\iota_\tau$-th column $\Col_{\iota_\tau}[A_\chi]$ of $A_\chi$.
\par
Now assume that all we know about $S_\chi$, $\chi=1,2$, is that $S_\chi\in \cS_{n_\chi}(Q_\chi,\rho_\chi)$ with known $Q_\chi$ and $\rho_\chi$. We observe the sequence $\eta^t=(\eta_1,...,\eta_t)$ coming from one of two chains, and want to decide on the hypotheses $H_\chi$, $\chi=1,2$, stating that $S_\chi\in\cS_{n_\chi}(Q_\chi,\rho_\chi)$.
\paragraph{Construction and result.} Our approach is as follows. Given a positive integer $\kappa$, for $\chi=1,2$ let
$$
Z_\chi=\Conv\{A_\chi v:\,v\in \Delta_{n_\chi},\; \mbox{and}\; \exists j: \, \|v-\Col_j[Q_\chi^\kappa]\|_1\leq\kappa\rho_\chi\}
\subset\Delta_m.
$$
Note that $Z_\chi\subset\rint \Delta_m$ (since the column sums in $A_\chi$ are equal to one, and all entries of $A_\chi$ are positive).\par
It is immediately seen that
\begin{itemize}
\item Under hypothesis $H_\chi$, $\chi=1,2$, for every positive integer $t$, the conditional,  given the state $\jmath_{\kappa(t-1),\chi}$ of the Markov chain $\chi$ at time $\kappa(t-1)$, distribution of observation $\eta_{\kappa t}$ belongs to $Z_\chi$.
    {\small \begin{quotation}
    Indeed,
     $S_\chi$ and $Q_\chi$ are stochastic matrices with $\|S_\chi-Q_\chi\|_{1,1} \leq\rho_\chi$ (we are under hypothesis $H_\chi$),  and for stochastic matrices  $A,B,\bar{A}$ and $\bar{B}$ one has
    \[
    \|\bar{A}\bar{B}-AB\|_{1,1}\leq \|\bar{A}-A\|_{1,1}+\|\bar{B}-B\|_{1,1}
    \] due to
    \[\begin{array}{l}
    \|\bar{A}\bar{B}-AB\|_{1,1}\leq \|\bar{A}(\bar{B}-B)\|_{1,1}+\|(\bar{A}-A)B\|_{1,1}\\
    \leq \|\bar{A}\|_{1,1}\|\bar{B}-B\|_{1,1}+
    \|\bar{A}-A\|_{1,1}\|B\|_{1,1}=\|\bar{B}-B\|_{1,1}+
    \|\bar{A}-A\|_{1,1}.
    \end{array}
    \]
    Whence $\|S_\chi^\kappa-Q_\chi^\kappa\|_{1,1}\leq\kappa\rho_\chi$, so that the probabilistic vector
    $v=\Col_{\jmath_{\kappa(t-1),\chi}}[S_\chi^\kappa]$ satisfy $\|v-\Col_{\jmath_{\kappa(t-1),\chi}}[Q_\chi^\kappa]\|_1\leq\kappa\rho_\chi$. We conclude that the distribution of $A_\chi v$ of  $\eta_{\kappa t}$  belongs to $Z_\chi$.
    \end{quotation}
    } 
\item $Z_\chi$ is a polyhedral convex set with an explicit  representation:
 \[Z_\chi=\left\{z:\exists \alpha,v^1,...,v^{n_\chi}\in\bbr^{n_\chi}:\begin{array}{l}
 z=A_\chi\sum_{j=1}^{n_\chi}v^j,\;v^j\geq0,\;\sum_{i=1}^{n_\chi}v^j_i=\alpha_j,\;\alpha\in\Delta_{n_\chi},\\
\|v^j-\alpha_j\Col_j[Q_\chi^\kappa]\|_1\leq\alpha_j\kappa\rho_\chi,\;\;1\leq j\leq n_\chi.
\end{array}\right\}
\]
\end{itemize}
Setting $\omega_t=\eta_{\kappa t}$, $\zeta_{t,\chi}=\jmath_{t\kappa,\chi}$, $\chi=1,2$, and $X_t=Z_1$, $Y_t=Z_2$, $t=1,2,...$, we arrive at the situation considered in Proposition \ref{propnonstationarynew}: under hypothesis $H_\chi$, $\chi=1,2$, $\omega_t$ is a deterministic function of $\zeta^t_\chi=(\zeta_{0,\chi},...,\zeta_{t,\chi})$, and the conditional, given $\zeta^{t-1}_\chi$, distribution of $\omega_t$ is $\mu_t=A_\chi\Col_{\jmath_{(t-1)\kappa,\chi}}[S_\chi^\kappa]$, which is a deterministic function of $\zeta^{t-1}_\chi$. Besides this, $\mu_t\in X_t\equiv Z_1$ under hypothesis $H_1$, and $\mu_t\in Y_t\equiv Z_2$ under hypothesis $H_2$. For these reasons, Proposition \ref{propnonstationarynew} implies
\begin{proposition}\label{Markov} Let $\kappa$ be such that $Z_1$ does not intersect $Z_2$. Let, further, $(x_*,y_*)$ be an optimal solution to the convex optimization problem
$$
\varepsilon_\star=\max_{x\in Z_1,y\in Z_2}\sum_{i=1}^m \sqrt{x_iy_i},
$$
and let
$$
\phi_*(i)={1\over 2}\ln([x_*]_i/[y_*]_i),\,1\leq i\leq m.
$$
Then for every positive integer $K$, the risk of the test $\phi_*^K$ which, given observation $\omega^K$, accepts $H_1$ whenever
\begin{equation}\label{eqMarkov}
\sum_{t=1}^K\phi_*(\omega_t)=\sum_{i=1}^m\phi_*(i)\Card\{t\leq K:\omega_t=i\}
\end{equation}
is nonnegative and accepts $H_2$ otherwise, does not exceed $\varepsilon_\star^K$.
\end{proposition}
\paragraph{Remarks.} Note that $\kappa$ meeting the premise of Proposition \ref{Markov} does exist, provided that $\rho_\chi$ are small enough and that $A_1 e\neq A_2 f$ for every pair of steady-state distributions $e=Q_1e$, $f=Q_2f$ of the chains with transition matrices $Q_1$ and $Q_2$.
\par
Note that in order to compute the test statistics \rf{eqMarkov} {\em we do not need to observe the trajectory $\omega_1,\omega_2,...,\omega_K$; all what matters is the ``histogram'' $\{p_i=\Card\{t\leq K:\omega_t=i\}\}_{i=1}^m$ of $\omega_1,...,\omega_K$}.
Furthermore, we lose nothing if instead of observing {\sl a single and long} $\omega$-trajectory,  we observe a {\em population of independent ``short'' trajectories}. Indeed, assume that $N$ independent
trajectories are observed on time horizon $L\kappa\leq K\kappa$;
all the trajectories start at time $\tau=0$
in a once for ever fixed state and then move from state to state independently of each other and utilizing the same transition matrix $S$. Our observations now are the total, over $N$ trajectories,  numbers $p_i$, $i=1,...,m$, of time instants of the form $\kappa t$, $t\geq1$, spent by the trajectories in state
$i$. If our goal is to decide which of the chains $\chi=1,2$ we are observing, it is immediately seen that Proposition \ref{propnonstationarynew} implies that under the premise and in the notation of Proposition \ref{Markov}, the test which accepts $H_1$ when $\sum_{i=1}^m \phi_*(i)p_i\geq0$ and accepts $H_2$ otherwise (cf. (\ref{eqMarkov})) obeys the upper risk bound $\varepsilon_\star^{LN}$. In other words, the risk of the test would be exactly the same as if instead of (aggregated partial) information on $N$ trajectories of length $L\kappa$ each we were collecting similar information on a single  trajectory of length $K=LN\kappa$.
\paragraph{Numerical illustration.} Consider a queuing system ($M/M/s/s+b$) with several identical servers and  a single buffer of capacity $b$. The service times of each server and inter-arrival times are exponentially distributed, with distributions $\E(\mu)$ and $\E(\lambda)$ respectively. Upon arrival, a customer either starts being served, when there are free servers, or joins the buffer queue, if all servers are busy and there are $<b$ customers in the buffer queue, or leaves the system immediately when all servers are busy and there are $b$ customers in the buffer. We assume that the
parameters $\lambda$, $\mu$ are not known exactly; all we know is that
\[
|\lambda-\bar{\lambda}|\leq\delta_\lambda\;\mbox{and}\;|\mu-\bar{\mu}|\leq\delta_\mu,
 \]
 with given $\bar{\lambda}>0$, $\bar{\mu}>0$ and
$\delta_\lambda<\bar{\lambda}$, $\delta_\mu<\bar{\mu}$.\par
We observe the number of customers in the buffer at times $t=1,2,...$, and want to decide on the hypotheses $H_1$ stating that the number of servers in the system is $s_1$, and $H_2$, stating that this number is $s_2$.
\par
In terms of the hidden Markov chain framework presented above, the situation is as follows.
Under hypothesis $H_\chi$ the queuing system can be modeled by Markov chain with $n_\chi=s_\chi+b+1$ states with the transition matrix of the chain $S_\chi=\exp\{L_\chi\}$, where the transition rate matrix $L_\chi=L_\chi(\lambda,\mu)$ satisfies
\[ [L_\chi]_{j-1,j}=s(j)\mu, \;\;[L_\chi]_{j,j}= -(s(j)\mu+\lambda),\;\;[L_\chi]_{j+1,j}=\lambda,\;\;s(j):=\min[j-1,s_\chi],\;1\leq j\leq n_\chi.
\]
It is immediately seen that if $Q_\chi=\exp\{L_\chi(\bar{\lambda},\bar{\mu})\}$, it holds\footnote{Indeed, we have $S_\chi=\lim_{k\to\infty}(I+{1\over k}L_\chi(\lambda,\mu))^k$; for large $k$, the matrix $N_k(\lambda,\chi)=I+{1\over k}L_\chi(\lambda,\mu)$ is stochastic, and we clearly have $\|N_k(\lambda,\mu)-N_k(\bar{\lambda},\bar{\mu})\|_{1,1}\leq k^{-1}\rho_\chi$. Whence, as we have already seen, \[\|N_k^k(\lambda,\mu)-N_k^k(\bar{\lambda},\bar{\mu})\|_{1,1}\leq \rho_\chi.
 \]When passing to the limit as $k\to\infty$, we get the desired bound on $\|S_\chi-Q_\chi\|_{1,1}$.}
  \[
   \|S_\chi-Q_\chi\|_{1,1}\leq \rho_\chi:=2\delta_\lambda+2s_\chi\delta_\mu.
   \]
We can now apply the outlined scheme to decide between the hypotheses $H_1$ and $H_2$. A numerical illustration is presented in table \ref{tabMarkov}; in this illustration, we use $\kappa=1$, that is, observations used in the test are the numbers of customers in the buffer at times $t=1,2,...,K$.
\begin{table}[h]
\begin{center}{\small
\begin{tabular}{||c||c||c||c|c||c|c||c|c||}
\cline{4-9}
\multicolumn{3}{c||}{}&\multicolumn{2}{|c||}{$K=K_*$}&\multicolumn{2}{c||}{$K=\rfloor K_*/2\lfloor$}
&\multicolumn{2}{c||}{$K=\rfloor K_*/3\lfloor$}\\
\hline
$s_1,s_2,b$&$\varepsilon_\star$&$K_*$&${\epsilon}_1$&${\epsilon}_2$&${\epsilon}_1$&${\epsilon}_2$&${\epsilon}_1$&${\epsilon}_2$\\
\hline
$s_1=10$, $s_2=9$, $b=5$&0.993240&679&0.0000&0.0000&0.0035&0.0015&0.0119&0.0104\\
\hline
$s_1=10$, $s_2=7$, $b=5$&0.894036&42&0.0002&0.0002&0.0093&0.0100&0.0260&0.0273\\
\hline
\end{tabular}}
\end{center}
\caption{\label{tabMarkov} Experiments with toy queuing systems. $\bar{\lambda}=40,\bar{\mu}=5,\rho_1=\rho_2=0$. ${\epsilon}_\chi$: empirical, over sample of $10^4$ experiments with observation time $K$ each, probability to reject $H_\chi$ when the hypothesis is true. $\varepsilon_\star$ is defined in Proposition \ref{Markov}, $K_*=\rfloor\ln(1/0.01)/\ln(1/\varepsilon_\star)\lfloor$ is the observation time, as defined by Proposition \ref{Markov}, resulting in risk $\leq 0.01$.}
\end{table}


\appendix
\section{Proofs}
\subsection{Proof of Theorem \ref{the1}}
\paragraph{1$^0$.} The fact that the function (\ref{Phi})
is continuous on its domain, convex in $\phi(\cdot)\in\F$ and concave in  $[x;y]\in X\times Y$ is readily given by our basic assumptions. Let us set
\begin{equation}\label{PsiPsi}
\Psi([x;y])=\inf_{\phi\in\F} \Phi(\phi,[x;y]).
\end{equation}
We claim that the function
\[
\phi_{x,y}(\omega)={1\over 2}\ln(p_{x}(\omega)/p_{y}(\omega))
\]
(which, by our assumptions, belongs to $\F$) is an optimal solution  to the right hand side minimization problem in (\ref{PsiPsi}), so that
\begin{equation}\label{sothat}
\forall (x\in X,y\in Y): \Psi([x;y]):=\inf_{\phi\in\F} \Phi(\phi,[x;y])=2\ln\left(\int_\Omega \sqrt{p_{x}(\omega)p_{y}(\omega)}
P(d\omega)\right).
\end{equation}
Note that $\Psi$, being the infinum of a family of concave functions of $[x;y]\in \M\times \M$, is concave on $\M\times\M$.
Indeed,
we have
$$
\exp\{-\phi_{x,y}(\omega)\}p_{x}(\omega)=\exp\{\phi_{x,y}(\omega)\}p_{y}(\omega)=g(\omega):=\sqrt{p_{x}(\omega)p_{y}(\omega)},
$$
whence $\Phi(\phi_{x,y},[x;y]) =2\ln\left(\int_\Omega g(\omega) P(d\omega)\right)$. On the other hand,
for $\phi(\cdot)=\phi_{x,y}(\cdot)+\delta(\cdot)\in\F$ we have
$$
\begin{array}{ll}
&\int_\Omega g(\omega)P(d\omega)=\int_\Omega \left[\sqrt{g(\omega)}\exp\{-\delta(\omega)/2\}\right]\left[\sqrt{g(\omega)}\exp\{\delta(\omega)/2\}\right]P(d\omega)\\
(a)&\leq
\left(\int_\Omega g(\omega)\exp\{-\delta(\omega)\}P(d\omega)\right)^{1/2} \left(\int_\Omega g(\omega)\exp\{\delta(\omega)\}P(d\omega)\right)^{1/2}\\
&=\left(\int_\Omega\exp\{-\phi(\omega)\}p_{x}(\omega) P(d\omega) \right)^{1/2}\left(\int_\Omega\exp\{\phi(\omega)\}p_{y}(\omega) P(d\omega) \right)^{1/2}\\
(b)&\Rightarrow 2\ln\left(\int_\Omega g(\omega)P(d\omega)\right)\leq \Phi(\phi,[x;y]),\\
\end{array}
$$
and thus $\Phi(\phi_{x,y},[x,y])\leq \Phi(\phi,[x;y])$ for every $\phi\in\F$.
\begin{remark}\label{rem166} {\rm Note that the inequality in $(b)$ can be equality only when the inequality in $(a)$ is so. In other words, if $\bar{\phi}$ is a minimizer of
$\Phi(\phi,[x;y])$ over $\phi\in\F$, setting $\delta(\cdot)=\bar{\phi}(\cdot)-\phi_{x,y}(\cdot)$, the functions $\sqrt{g(\omega)}\exp\{-\delta(\omega)/2\}$ and $\sqrt{g(\omega)}\exp\{\delta(\omega)/2\}$, considered as elements of $L_2[\Omega,P]$, are proportional to each other. Since $g$ is positive and $g,\delta$ are continuous, while the support of $P$ is the entire $\Omega$, this ``$L_2$-proportionality'' means that the functions in question differ by a constant factor, or, which is the same, that $\delta(\cdot)$ is constant. Thus, {\sl the minimizers of $\Phi(\phi,[x;y])$ over $\phi\in \F$ are exactly the functions of the form $\phi(\omega)=\phi_{x,y}(\omega)+{\rm const}$.}}
\end{remark}
\paragraph{2$^0$.}
We are about to verify that $\Phi(\phi,[x;y])$ has a saddle point   ($\min$ in $\phi$, $\max$ in $[x;y]$) on $\F\times (X\times Y)$.
Indeed, observe, first, that on the domain of $\Phi$ it holds
\begin{equation}\label{shift}
\Phi(\phi(\cdot)+a,[x;y])=\Phi(\phi(\cdot),[x;y])\,\,\forall (a\in\bbr,\phi\in\F).
 \end{equation}
 Let $\bar{x}\in\M$ and let $\bar{P}$ be the probability measure with density $p_{\bar{x}}$ w.r.t. $P$. Since the observation scheme in quesiton is good, for $\phi\in\F$ we have $\int_\Omega\exp\{\pm\phi(\omega)\}\bar{P}(d\omega)<\infty$, implying that $\phi\in L_1[\Omega,\bar{P}]$. Let $\F_0=\{f\in\F: \int_\Omega \phi(\omega)\bar{P}(d\omega)=0\}$, so that $\F_0$ is a linear subspace in $\F$, and all functions from $\F$ are obtained from functions from $\F_0$ by adding constants. Invoking (\ref{shift}), to prove existence of a saddle point of $\Phi$ on $\F\times(X\times Y)$ is the same as to prove that $\Phi$ has a saddle point on $\F_0\times(X\times Y)$.
Since $X\times Y$ is a convex compact set,  $\Phi$ is continuous on $\F_0\times(X\times Y)$ and convex-concave, all we need in order to verify the existence of a saddle point is to show that $\Phi$ is coercive in the first argument, that is, for every fixed $[x;y]\in X\times Y$ one has $\Phi(\phi,[x;y])\to+\infty$ as $\phi\in\F_0$ and  $\|\phi\|\to\infty$ (whatever be the norm $\|\cdot\|$ on $\F_0$; recall that $\F_0$ is a finite-dimensional linear space). Setting $\Theta(\phi)=\Phi(\phi,[x;y])$ and taking into account that $\Theta$ is convex and finite on $\F_0$, in order to prove that $\Theta$ is coercive, it suffices to verify that $\Theta(t\phi)\to\infty$, $t\to\infty$, for every nonzero $\phi\in \F_0$, which is evident: since $\int_\Omega\phi(\omega)\bar{P}(d\omega)=0$ and  $\phi$ is nonzero, we have $\int_\Omega\max[\phi(\omega),0]\bar{P}(d\omega)= \int_\Omega\max[-\phi(\omega),0]\bar{P}(d\omega)>0$, whence $\Theta(t\phi)\to\infty$ as $t\to\infty$ due to the fact that both $p_{x}(\cdot)$ and $p_{y}(\cdot)$ are positive everywhere and the support of $\bar{P}$ is the entire $\Omega$.
\paragraph{3$^0$.}
Now let $(\phi_*(\cdot);[x_*;y_*])$ be a saddle point of $\Phi$ on $\F\times(X\times Y)$. Shifting, if necessary, $\phi_*(\cdot)$ by a constant (by (\ref{shift}), this does not affect the fact that $(\phi_*,[x_*;y_*])$ is a saddle point of $\Phi$), we can assume that
\begin{equation}\label{balance1}
\varepsilon_\star:=\int_\Omega\exp\{-\phi_*(\omega)\}p_{x_*}(\omega)P(d\omega)= \int_\Omega\exp\{\phi_*(\omega)\}p_{y_*}(\omega)P(d\omega),
\end{equation}
so that the saddle point value of $\Phi$ is
\begin{equation}\label{saddlevalue}
\Phi_*:=\max_{[x;y]\in X\times Y}\min_{\phi\in\F}\Phi(\phi,[x;y])=\Phi(\phi_*,[x_*;y_*])=2\ln(\varepsilon_\star).
\end{equation}
The following lemma completes the proof of Theorem \ref{the1}.i:
\begin{lemma}\label{prop1} Under the premise of Theorem \ref{the1}, let $(\phi_*,[x_*;y_*])$ be a saddle point of $\Phi$ satisfying
{\rm (\ref{balance1})}, and let $\phi_*^a(\cdot)=\phi_*(\cdot)-a$, $a\in\bbr$. Then
\begin{equation}\label{then1}
\begin{array}{lrcl}
(a)&\int_\Omega \exp\{-\phi_*^a(\omega)\}p_{x}(\omega)P(d\omega)&\leq&\exp\{a\}\varepsilon_\star\,\,\forall x\in X,\\
(b)&\int_\Omega \exp\{\phi_*^a(\omega)\}p_{y}(\omega)P(d\omega)&\leq&\exp\{-a\}\varepsilon_\star\,\,\forall y\in Y.\\
\end{array}
\end{equation}
As a result, for the simple test associated with the detector $\phi_*^a$, the probabilities $\epsilon_X$ to reject $H_X$ when the hypothesis is true and $\epsilon_Y$ to reject $H_Y$ when the hypothesis is true can be upper-bounded according to {\rm (\ref{upperbounded}).}
\end{lemma}
{\bf Proof.} For $x\in X$, we have
$$
\begin{array}{rcl}
2\ln(\varepsilon_\star)&=&\Phi_*\geq \Phi(\phi_*,[x;y_*])\\
&=&\ln\left(\int_\Omega\exp\{-\phi_*(\omega)\}p_{x}(\omega)P(d\omega)\right)+
\ln\left(\int_\Omega\exp\{\phi_*(\omega)\}p_{y_*}(\omega)P(d\omega)\right)\\
&=&\ln\left(\int_\Omega\exp\{-\phi_*(\omega)\}p_{x}(\omega)P(d\omega)\right)+\ln(\varepsilon_\star),\\
\end{array}
$$
whence $\ln\left(\int_\Omega\exp\{-\phi_*^a(\omega)\}p_{x}(\omega)P(d\omega)\right)=
\ln\left(\int_\Omega\exp\{-\phi_*(\omega)\}p_{x}(\omega)P(d\omega)\right)+a\leq\ln(\varepsilon_\star)+a$, and (\ref{then1}.$a$) follows. Similarly,
 when $y\in Y$, we have
$$
\begin{array}{rcl}
2\ln(\varepsilon_\star)&=&\Phi_*\geq \Phi(\phi_*,[x_*;y])\\
&=&\ln\left(\int_\Omega\exp\{-\phi_*(\omega)\}p_{x_*}(\omega)P(d\omega)\right)+
\ln\left(\int_\Omega\exp\{\phi_*(\omega)\}p_{y}(\omega)P(d\omega)\right)\\
&=&\ln(\varepsilon_\star)+\ln\left(\int_\Omega\exp\{\phi_*(\omega)\}p_{y}(\omega)P(d\omega)\right),\\
\end{array}
$$
so that $\ln\left(\int_\Omega\exp\{\phi_*^a(\omega)\}p_{y}(\omega)P(d\omega)\right)=
\ln\left(\int_\Omega\exp\{\phi_*(\omega)\}p_{y}(\omega)P(d\omega)\right)-a\leq\ln(\varepsilon_\star)-a$, and (\ref{then1}.$b$) follows.
\mypar
Now let $x\in X$, and let $\epsilon(x)$ be the probability for the test, the detector being $\phi_*^a$, to reject $H_X$; this is at most the probability for $\phi_*^a(\omega)$ to be nonpositive when $\omega\sim p_{x}(\cdot)$, and therefore
$$
\epsilon(x)\leq \int_\Omega\exp\{-\phi_*^a(\omega)\}p_{x}(\omega)P(d\omega),
$$
so that $\epsilon(x)\leq\exp\{a\}\varepsilon_\star$ by (\ref{then1}.$a$). Thus, the probability for our test to reject the hypothesis $H_X$ when it is true is $\leq\exp\{a\}\varepsilon_\star$. Relation (\ref{then1}.$b$) implies in the same fashion that the probability for our test to reject $H_Y$ when this hypothesis is true is $\leq\exp\{-a\}\varepsilon_\star$.
\paragraph{4$^0$.}
Theorem \ref{the1}.ii is readily given by the following
\begin{lemma}\label{prop2} Under the premise of Theorem \ref{the1}, let $(\phi_*,[x_*;y_*])$ be a saddle point of $\Phi$, and let $\epsilon\geq0$ be such that there exists a (whatever) test for deciding between two simple hypotheses
\begin{equation}\label{AB1}
(A): \omega\sim p(\cdot):=p_{x_*}(\cdot),\quad
(B): \omega\sim q(\cdot):=p_{y_*}(\cdot)\\
\end{equation}
with the sum of error probabilities $\leq2\epsilon$. Then
\begin{equation}\label{then2}
\varepsilon_\star\leq 2\sqrt{(1-\epsilon)\epsilon}.
\end{equation}
\end{lemma}
{\bf Proof.} Under the premise of the lemma, $(A)$ and  $(B)$ can be decided with the sum of error probabilities $\leq2\epsilon$, and therefore the test affinity of $(A)$ and $(B)$ is bounded by $2\epsilon$:
$$
\int_\Omega\min[p(\omega),q(\omega)]P(d\omega)\leq 2\epsilon.
$$
On the other hand, we have seen that the saddle point value of $\Phi$ is $2\ln(\varepsilon_\star)$; since $[x_*;y_*]$ is a component of a saddle point of $\Phi$, it follows that
$\min_{\phi\in \F}\Phi(\phi,[x_*;y_*])=2\ln(\varepsilon_\star)$. The left hand side in this equality, as we know from item 1$^0$, is $\Phi(\phi_{x_*,y_*},[x_*;y_*])$, and we arrive at $2\ln(\varepsilon_\star)=\Phi({1\over 2}\ln(p_{x_*}(\cdot)/p_{y_*}(\cdot)),[x_*;y_*])=2\ln\left(\int_\Omega\sqrt{p_{x_*}(\omega)p_{y_*}(\omega)}P(d\omega)\right)$, so that $\varepsilon_\star=
\int_\Omega\sqrt{p_{x_*}(\omega)p_{y_*}(\omega)}P(d\omega)=\int_\Omega \sqrt{p(\omega)q(\omega)} P(d\omega)$. We now have {(cf. \cite[chapter 4]{Lecam1986})}
$$
\begin{array}{l}
\varepsilon_\star=\int_\Omega\sqrt{p(\omega)q(\omega)}P(d\omega)=\int_\Omega\sqrt{\min[p(\omega),q(\omega)]} \sqrt{\max[p(\omega),q(\omega)]} P(d\omega)\\
\leq \left(\int_\Omega \min[p(\omega),q(\omega)]P(d\omega)\right)^{1/2}\left(\int_\Omega \max[p(\omega),q(\omega)]P(d\omega)\right)^{1/2}
\leq \sqrt{2(2-2\epsilon)\epsilon}= 2\sqrt{(1-\epsilon)\epsilon}.
\end{array}
$$
\paragraph{5$^0$.} We have proved items (i) and (ii) of Theorem \ref{the1}. To complete the proof of the theorem, it remains to justify (\ref{phistar}). Thus, let
$(\phi_*,[x_*;y_*])$ be a saddle point of $\Phi$ satisfying (\ref{balance1}).  All we need to prove is that $\phi_*$ is nothing but $$\bar{\phi}(\cdot)={1\over 2}\ln\left(p_{x_*}(\cdot)/p_{y_*}(\cdot)\right).$$
Indeed, the function $\Phi(\cdot,[x_*;y_*])$ attains its minimum on $\F$ at the point $\phi_*$; by Remark \ref{rem166}, it follows that $\phi_*(\cdot)-\bar{\phi}(\cdot)$ is constant on $\Omega$; since both $\bar{\phi}$ and $\phi_*$ satisfy (\ref{balance1}), this constant is zero. \qed
\subsection{Proofs of Propositions \ref{propnonstation1pairnewnot} and \ref{propnonstationarynew}}\label{sect:proofspropnonstationary}
Proposition \ref{propnonstation1pairnewnot} is a simple particular case of  Proposition \ref{propnonstationarynew} which we prove here.\par
Observe that when $t\leq K$ and $p\in X_{t}$, so that $p\in X_{it}$ for some $i\in\I_t$, we have by definition of $\phi_t$, see (\ref{aijs}),
\be
\lefteqn{\int_{\Omega_t}\exp\{-\phi_t(\omega_t)\}p(\omega_t)P_t(d\omega_t)=
\int_{\Omega_t}\exp\{\min_{r\in\I_t}\max_{s\in\J_t}[a_{rst}-\phi_{rst}(\omega_t)]\}p(\omega_t)P_t(d\omega_t)}\nn
&\leq&\int_{\Omega_t}\exp\{\max_{s\in\J_t}[a_{ist}-\phi_{ist}(\omega_t)]\}p(\omega_t)P_t(d\omega_t)\leq \sum\limits_{s\in\J_t}\int_{\Omega_t}\exp\{a_{ist}-\phi_{ist}(\omega_t)\}p(\omega_t)P_t(d\omega_t)\nn
&\leq& \sum\limits_{s\in\J_t}\exp\{a_{ist}\}\epsilon_{ist}=\sum\limits_{s\in\J_t}h^t_s\epsilon_{ist}/g^t_i \;\;\hbox{[see (\ref{suchthatnew177}.$a$), (\ref{aijs})]}\nn
&=&[E_th^t]_i/g^t_i=\varepsilon_t
\;\;\;\hbox{[see  (\ref{Perron})]}.
\ee{ccase1}
Similarly, when $t\leq K$ and $p\in Y_{t}$, so that $p\in Y_{jt}$ for some $j\in \J_t$, we have
\be
\lefteqn{\int_{\Omega_t}\exp\{\phi_t(\omega_t)\}p(\omega_t)P_t(d\omega_t)=\int_{\Omega_t}\exp\{\max_{r\in\I_t}\min_{s\in\J_t}[\phi_{rst}(\omega_t)-a_{rst}]\}p(\omega_t)P_t(d\omega_t)
}\nn
&\leq&\int_{\Omega_t}\exp\{\max_{r\in\I_t}[\phi_{rjt}(\omega_t)-a_{rjt}]\}p(\omega_t)P_t(d\omega_t)\leq \sum\limits_{r\in\I_t}\int_{\Omega_t}\exp\{\phi_{rjt}(\omega_t)-a_{rjt}\}p(\omega_t)P_t(d\omega_t)\nn
&\leq& \sum\limits_{r\in\I_t}\exp\{-a_{rjt}\}\epsilon_{rjt}=\sum_{r\in \I_t}g^t_r\epsilon_{rjt}/h^t_j \hbox{\ [see (\ref{suchthatnew177}.$b$), (\ref{aijs})]}\nn
&=&[E_t^Tg^t]_j/h^t_j=\varepsilon_t\hbox{\ [see (\ref{Perron})]}.
\ee{ccase2}
Now let $H_1=H_X$ be true, let $\bE_{|\zeta^{t-1}_1}\{\cdot\}$ stand for the conditional expectation, $\zeta^{t-1}_1$ being fixed,
and let $p_{\zeta^{t-1}_1}(\cdot)$ be conditional, $\zeta^{t-1}_1$ being fixed, probability density of $\omega_t$  w.r.t. $P_t$, so that $p_{\zeta^{t-1}_1}(\cdot)\in X_t$ for all $\zeta^{t-1}_1$ and all $t\leq K$. We have
$$
\begin{array}{l}
\bE\left\{\exp\{-\phi_1(\omega_1)-...-\phi_{t}(\omega_t)\}\right\}=\bE\left\{\exp\{-\phi_1(\omega_1)-...-\phi_{t-1}(\omega_{t-1})\}\bE_{|\zeta^{t-1}_1}\{\exp\{-\phi_t(\omega_t)\}\}\right\}\\
=\bE\left\{\exp\{-\phi_1(\omega_1)-...-\phi_{t-1}(\omega_{t-1})\}\int_{\Omega_t}\exp\{-\phi_t(\omega_t)\}p_{\zeta^{t-1}_1}(\omega_t)P_t(d\omega_t)\}\right\}\\
\leq \varepsilon_t\bE\left\{\exp\{-\phi_1(\omega_1)-...-\phi_{t-1}(\omega_{t-1})\}\right\},\\
\end{array}
$$
where the concluding inequality is due to (\ref{ccase1}). From the resulting recurrence,
$$
\bE\{\exp\{-\phi^K(\omega^K)\}\}\leq {\prod}_{t=1}^K\varepsilon_t.
$$
This inequality combines with the description of our test to imply that the probability to reject $H_X$ when it is true is at most $\prod_{t=1}^K\varepsilon_t$.
 \par
Now assume that $H_2=H_Y$ holds true, so that the conditional, $\zeta^{t-1}_2$ being fixed, distribution $p_{\zeta^{t-1}_2}(\cdot)$ of $\omega_t$ belongs to $Y_t$ for all $\zeta^{t-1}_2$ and all $t\leq K$. Applying the previous reasoning to $-\phi^K$ in the role of $\phi^K$, $\zeta^t_2$ in the role of $\zeta^t_1$, and (\ref{ccase2}) in the role of (\ref{ccase1}), we conclude that the probability to reject $H_Y$ when it is true is at most $\prod_{t=1}^K\varepsilon_t$. \qed
\subsection{Proof of Proposition \ref{prop_alpha1}}
\paragraph{1$^0$.} The matrix  $\bar{E}=\left[p_i\epsilon_{ij}\right]_{1\leq i,j\leq m}$ has zero diagonal and positive off-diagonal entries. By the Perron-Frobenius theorem, the largest in magnitude eigenvalue of $\bar{E}$ is some positive real $\rho$, and the corresponding eigenvector $g$ can be selected to be nonnegative.
In addition, $g\geq0$ is in fact positive, since the relation
$$
\rho g_i=[\bar{E}g]_i
$$
along with the fact the all $p_i$ and all off-diagonal entries in $E$ are positive, allows for $g_i=0$ only if all the entries $g_j$ with $j\neq i$ are zeros, that is, only when $g=0$, which is impossible. Since $g>0$, we can set
$$
\alpha_{ij}=\bar{\alpha}_{ij}:=\ln(g_j)-\ln(g_i),
$$
thus ensuring $\alpha_{ij}=-\alpha_{ji}$ and
$$
p_i\varepsilon_i=\sum_{j=1}^mp_i\epsilon_{ij}\exp\{\alpha_{ij}\}=\sum_{j=1}^mp_i\epsilon_{ij}g_j/g_i=g_i^{-1}\sum_{j=1}^mp_i\epsilon_{ij}g_j=
g_i^{-1}[\bar{E}g]_i=\rho.
$$
Thus, with our selection of $\alpha_{ij}$ we get
$$
\varepsilon=\rho.
$$
\paragraph{2$^0$.} We claim that in fact  $\varepsilon_*=\rho$, that is, the feasible solution $[\bar{\alpha}_{ij}]$ is optimal for (\ref{varepsilonproblem}). Indeed, otherwise
there exists a feasible solution $[\alpha_{ij}=\bar{\alpha}_{ij}+\delta_{ij}]_{i,j}$ with $\delta_{ij}=-\delta_{ji}$ such that
$$
\bar{\rho}=\max_i \left[p_i\sum_j\epsilon_{ij}\exp\{\alpha_{ij}\}\right]<\rho.
$$
As we have shown, for every $i$ we have $\rho=\sum_jp_i\epsilon_{ij}\exp\{\bar{\alpha}_{ij}\}$. It follows that the convex functions
$$
f_i(t)=\sum_jp_i\epsilon_{ij}\exp\{\bar{\alpha}_{ij}+t\delta_{ij}\}
$$
all are equal to $\rho$ when $t=0$ and are $\leq\bar{\rho}<\rho$ when $t=1$, whence, due to convexity of  $f_i$,  for every $i$ one has
$$
0>{d\over dt}\big|_{t=0}f_i(t)=\sum_j p_i\epsilon_{ij}\exp\{\bar{\alpha}_{ij}\} \delta_{ij}=p_i\sum_jg_jg_i^{-1}\epsilon_{ij}\delta_{ij}.
$$
Multiplying the resulting inequalities by $g_i^2/p_i>0$ and summing up the results over $i$, we get
$$
0>\sum_{i,j}g_ig_j\epsilon_{ij}\delta_{ij}.
$$
This is impossible, since $\epsilon_{ij}=\epsilon_{ji}$ and $\delta_{ij}=-\delta_{ji}$, and the right hand side in the latter inequality is zero.\qed
\subsection{Proof of Proposition \ref{theverylatestopt}}
 In the notation and under the premise of the proposition, let $\widehat{\epsilon}_{ij}$ be the risks of detectors $\phi_{ij}$ as defined in Theorem \ref{the1}, so that  $\widehat{\epsilon}_{ij}^K$ are the risks of $\phi_{ij}^K$. Denote $\delta$ the maximum of the risks $\widehat{\epsilon}_{ij}$ taken over all ``far from each other'' pairs of indexes $(i,j)$, that is, pairs such that $i,j$ do not belong to the same group $\I_\ell,\,\ell=1,...,L$, and let $\bar{i},\bar{j}$ be two ``far from  each other'' indexes such that  $\delta=\widehat{\epsilon}_{\bar{i}\bar{j}}$. Test $\overline{T}$ clearly induces a test for deciding on the pair of hypotheses $H^1:=H_{\bar{i}}$, $H^2:=H_{\bar{j}}$ from observation $\omega^{\bar{K}}$ which does not accept $H^\chi$, $\chi=1,2$, when the hypothesis is true, with probability at most $\epsilon$, and never accepts both these hypotheses simultaneously. Same as in the proof of Proposition \ref{newcol1}, the latter implies that $\delta^{\bar{K}}=[\widehat{\epsilon}_{\bar{i}\bar{j}}]^{\bar{K}}\leq2\sqrt{\epsilon}$. Since the nonzero entries in the matrix $D=D_K$ participating in the description of the test $\widehat{\T}^K$ are of the form $\widehat{\epsilon}_{ij}^K$ with ``far from each other'' $i,j$, the entries in the entrywise nonnegative matrix $D_K$ do not exceed $\delta^K\leq [2\sqrt{\epsilon}]^{K/\bar{K}}$. Therefore the spectral norm of $D_K$ (which, as we know,  upper bounds the risk of $\widehat{\T}^K$) does not exceed $M[2\sqrt{\epsilon}]^{K/\bar{K}}$, and the conclusion of Proposition \ref{theverylatestopt} follows. \qed
\subsection{Proofs of Propositions \ref{propPoisson} and \ref{pro:normal}}
We prove here Proposition \ref{propPoisson}, the proof of Proposition \ref{pro:normal} can be conducted following same lines.
\paragraph{1$^0$.} Let us fix $i$. It is immediately seen that problem $(P^i_\epsilon)$ is solvable (recall that $Ae[i]\neq0$); let $\rho^i=\rho_i^P(\epsilon)$, $r^i$, $u^i$, $v^i$ be an optimal solution to this problem. We clearly have $r^i=\rho^i$. We claim that the optimal value in the optimization problem
$$
\min_{r,u,v} \left\{{1\over 2}\sum_\ell \left[\sqrt{[Au]_\ell}-\sqrt{[A(re[i]+v^i)]_\ell}\right]^2: u\in\cV,v\in\cV, \rho^i\leq r\leq R\right\}\eqno{(P)}
$$
is $\ln(\sqrt{n}/\epsilon)$, while $(r^i,u^i,v^i)$ is an optimal solution to the problem. Indeed, taking into account the origin of $u^i,v^i,\rho^i=r^i$ and the relation $R\geq \rho_i^P(\epsilon)$,
$(r^i,u^i,v^i)$ is a feasible solution to this problem with the value of the objective $\leq \ln(\sqrt{n}/\epsilon)$; thus, all we need in order to support our claim is to verify that the optimal value in $(P)$ is  $\geq\ln(\sqrt{n}/\epsilon)$. To this end assume for a moment that $(P)$ has a feasible solution
$(\bar{r},\bar{u},\bar{v})$ with the value of the objective $<\ln(\sqrt{n}/\epsilon)$. Then, setting $\rho^+=\rho^i+\delta$, $r^+=\bar{r}+\delta$, $u^+=\bar{u}$, $v^+=\bar{v}$ and choosing $\delta>0$ small enough, we clearly get a feasible solution to $(P^i_\epsilon)$ with the value of the objective $>\rho^i=\rho^P_i(\epsilon)$, which is impossible. Our claim is justified.
\paragraph{2$^0$.} Recalling the ``Poisson case'' discussion in section \ref{sec:mainres}, item {\bf 1$^0$} implies that the simple test associated with the detector $\phi_i(\cdot)$ given by (\ref{phiP})
decides between the hypotheses $H_0$ and $H^i(\rho^P_i(\epsilon))$ with probabilities of errors $\leq \epsilon/\sqrt{n}$. Since $H^i(r)$ ``shrinks'' as $r$ grows, we conclude that whenever $\rho_i\in[\rho^P_i(\epsilon),R]$, the same test decides between the hypotheses $H_0$ and $H^i(\rho_i)$ with probabilities of errors not exceeding $\epsilon/\sqrt{n}$. Now let $\rho=[\rho_1;...;\rho_n]$ satisfy the premise of Proposition \ref{propPoisson}, so that $\rho_i\geq\rho^P_i(\epsilon)$ for all $i$.  Note that the problem of testing $H_0:\,\mu\in X$ against $H_1(\rho):\mu\in \bigcup_{i=1}^n Y(\rho_i)$, along with the tests $\phi_{1i}(\cdot)=\phi_i(\cdot)$, $i=1,...,n$ satisfy the premise of Proposition \ref{propnonstation1pairnewnot} with $\epsilon_{1i}=\epsilon/\sqrt{n}$, $\varepsilon=\sqrt{\sum_{i=1}^n \epsilon_{1i}^2}(=\epsilon)$, and $a_{1i}=-\half\ln n$, $i=1,...,n$. As a result, by Proposition \ref{propnonstation1pairnewnot}, the risk of the test $\phi^P(\cdot)$ does not exceed $\epsilon$.

 \paragraph{3$^0$.} To justify the bound on rate optimality, let us set
 $$
 \Opt_i(\rho)=\min_{r,u,v} \left\{{1\over 2}\sum_\ell \left[\sqrt{[Au]_\ell}-\sqrt{[A(re[i]+v^i)]_\ell}\right]^2: u\in\cV,v\in\cV, \rho\leq r\leq R\right\}\eqno{[\rho\geq0]}
 $$
 The function $\Opt(\rho)$ by its origin is a nondecreasing convex function on the segment $0\leq\rho\leq R$, $\Opt_i(\rho)=+\infty$ when $\rho>R$, and $\Opt(0)=0$. It follows that
 \begin{equation}\label{ineqtriv}
 \forall (\rho\in[0,R],\theta\geq 1): \Opt_i(\theta\rho)\geq\theta\Opt_i(\rho)
 \end{equation}

 Now assume that for some $\rho=[\rho_1;...;\rho_n]$ and $\epsilon\in(0,{1/4})$ there exists a test which decides between $H_0$ and $H_1(\rho)$ with probability of error $\leq\epsilon$. Taking into account the union structure of $H_1(\rho)$, for every fixed $i$ this test decides with the same probabilities of errors between the hypotheses $H_0$ and $H^i(\rho_i)$. All we need in order to prove the bound on the rate of optimality of $\widehat{\phi}_P$ is to extract from the latter observation that $\rho^P_i(\epsilon)/\rho_i\leq\kappa_n:=\kappa_n(\epsilon)$ for every $i$. Let us fix $i$ and verify that  $\rho^P_i(\epsilon)/\rho_i\leq\kappa_n$. There is nothing to do when $\rho_i\geq\rho^P_i(\epsilon)$ (due to $\kappa_n\geq1$); thus, assume that $\rho_i<\rho^P_i(\epsilon)$. Note that $\rho_i>0$ (since otherwise the hypotheses $H_0$ and $H^i(\rho_i)$ have a nonempty intersection and thus cannot be decided with probabilities of errors $<1/2$, while we are in the case of $\epsilon<{1/4}$). Applying Theorem \ref{the1} to the pair of hypotheses $H_0$, $H^i(\rho_i)$, it is straightforward to see that in this case item (ii) of Theorem states exactly that $\exp\{-\Opt_i(\rho_i)\}\leq 2\sqrt{\epsilon}$, or, which is the same,
 $\Opt(\rho_i)\geq \delta:=\half\ln(1/\epsilon)-\ln(2)$; $\delta$ is positive due to $\epsilon\in(0,{1/4})$. Now let $\theta>\ln(\sqrt{n}/\epsilon)/\delta$, so that $\theta\geq 1$. By (\ref{ineqtriv}),
we either have $\theta\rho_i>R$, whence $\theta\rho_i\geq\rho^P_i(\epsilon)$ due to  $\rho^P_i(\epsilon)\leq R$, or $\theta\rho_i\leq R$ and $\Opt_i(\theta\rho_i)>\ln(\sqrt{n}/\epsilon)$. In the latter case, as we have seen in item {\bf 1}$^0$ of the proof, it holds $\Opt_i(\rho^P_i(\epsilon))=\ln(\sqrt{n}/\epsilon)$,
and thus $\rho^P_i(\epsilon)<\theta\rho_i$ since $\Opt_i$ is nondecreasing in $[0,R]$. Thus, in all cases $\theta\rho_i>\rho^P_i(\epsilon)$ whenever $\theta>\ln(\sqrt{n}/\epsilon)/\delta$. But the latter ratio is exactly $\kappa_n$, and we conclude that $\kappa_n\rho_i\geq\rho_i^P(\epsilon)$, as required. \qed
\corr{}{
\subsection{Proof of Proposition \ref{proploglog}}
Let $d>0$, $J\geq2$, $\epsilon\in(0,0.05)$ and $\kappa\geq1$ satisfy (\ref{letdsatisfy}). Relation (\ref{ifdsatisfies}) is trivially true when $S=1$ and thus $\bar{K}=\bar{k}(1)=1$. Assume from now on that $S>1$, so that (\ref{seq10}) is not satisfied when $S=1$ and therefore $d\leq \ln(\Theta)$.
\paragraph{1$^0$.} We claim that
\begin{equation}\label{eq12345}
S\leq \bar{S}:=\rfloor S_*:={(2\kappa+1)\ln(\Theta)/d}\lfloor.
\end{equation}
To justify this claim, it suffices to verify that
\begin{equation}\label{eq12346}
{\ln(\Theta \bar{S})/\bar{k}(\bar{S})}<d.
\end{equation}
Indeed, we have $S_*\geq 3$, whence $S_*\leq \bar{S}\leq {4\over 3}S_*$, and $\bar{k}(S)\geq \bar{S}\geq S_*$, and therefore the left hand side in (\ref{eq12346}) is at most $\C d$, with $\C$ given by
$$
\C={\ln\left({4\over 3}[2\kappa+1]\Theta\ln(\Theta)/d\right)\over[2\kappa+1]\ln(\Theta)}\leq
\bar{\C}:={\ln\left({4\over 3}[2\kappa+1]\Theta\ln(\Theta)\right)+\kappa \ln(\Theta)\over [2\kappa+1]\ln(\Theta)},
$$
where the concluding $\leq$ is due to (\ref{letdsatisfy}). From the relations $\Theta=2.1J^2/\epsilon\geq 168$ and $\kappa\geq1$ it immediately follows that
$\bar{\C}<1$, implying the validity of  (\ref{eq12346}) and thus -- the validity of (\ref{eq12345}).
\paragraph{2$^0$.} We are in the situation $S>1$, whence $d\leq {\ln([S-1]\Theta)\over \bar{k}(S-1)}\leq {\ln([\bar{S}-1]\Theta)\over \bar{k}(S-1)}\leq {\ln(S_*\Theta)\over \bar{k}(S-1)}$, so that
$$
\begin{array}{c}
\bar{K}=\bar{k}(S)\leq 2\bar{k}(S-1)\leq {2\ln(S_*\Theta)\over d}= {2\ln([2\kappa+1]\Theta\ln(\Theta)/d)\over d}\leq
{2\left[\ln([2\kappa+1]\Theta\ln(\Theta))+\kappa\ln(\Theta)\right]\over d}=\D{\kappa\ln(\Theta)\over d},\\
\D={2[\ln([2\kappa+1]\ln(\Theta))+\kappa+1]\over\kappa\ln(\Theta)}
\end{array}
$$
(we have used the already established relation $S\leq\bar{S}$ and (\ref{letdsatisfy})). Since $\Theta\geq168$ and $\kappa\geq1$, it is immediately seen that
$\D\leq 2$, and we arrive at (\ref{ifdsatisfies}). \qed}{}
\subsection{Proof of Proposition \ref{pro:cdf}}
\paragraph{1$^0$.} Let the premise in Proposition \ref{pro:cdf} hold true, and let us set $\varrho=\rho[\epsilon]$. Observe, first, that $\Opt[\varrho]=\ln \epsilon$. Indeed, problem (\ref{rho_A}) clearly is solvable, and $\bar{x},\bar{y},r=\varrho$ is an optimal solution to this problem. $(\bar{x},\bar{y})$ is a feasible solution to $(F_{g,\alpha}[\varrho])$, whence the optimal value in the latter problem is at least $\ln \epsilon$. Now let us lead to a contradiction the assumption that $\Opt[\varrho]>\ln \epsilon$. Under this assumption, let $x_0\in H_0[\rho_{\max}]$, $y_0\in H_1[\rho_{\max}]$, and let $(\hat{x},\hat{y})$ be an optimal solution to $(F_{g,\alpha}[\varrho])$, so that
\begin{equation}\label{sothat17}
\sum_{\ell=1}^LK_\ell\ln\left(\sum_{i=1}^{n_\ell}\sqrt{[A^\ell x]_i[A^\ell y]_i }\right)
> \ln \epsilon
\end{equation}
when $x=\hat{x}$, $y=\hat{y}$. Now let $x_t=\hat{x}+t(x_0-\hat{x})$, $y_t=\hat{y}+t(y_0-\hat{y})$. Since \rf{sothat17} hods true for $x=\hat{x}$, $y=\hat{y}$, for small enough positive $t$ we have
$$
g^Tx_t\leq \alpha-\varrho-t(\rho_{\max}-\varrho),\,\,g^Ty_t\geq \alpha+\varrho+t(\rho_{\max}-\varrho),\,\,\sum_{\ell=1}^LK_\ell\ln\left(\sum_{i=1}^{n_\ell}\sqrt{[A^\ell x_t]_i[A^\ell y_t]_i }\right)\geq {\ln} \epsilon.
$$
which, due to $\rho_{\max}>\varrho$, contradicts the fact that $\varrho$ is the optimal value in (\ref{rho_A}).
\paragraph{2$^0$.} Let us prove (\ref{wehaveq7}). This relation is trivially true when $\varrho=0$, thus assume that $\varrho>0$. Since $\rho_{\max}\geq0$, and $g^Tx$ takes on $X$ both values $\leq\alpha$ and values $\geq\alpha$, this implies, by convexity of $\X$, that $g^Tx$ takes value $\alpha$ somewhere on $X$. Therefore, the hypotheses $H_0[0]$ and $H_1[0]$ intersect, whence $\Opt[0]=0$. In addition to this, due to its origin, $\Opt[\rho]$ is a concave function of $\rho\in[0,\varrho]$. Thus, $\Opt[\theta\varrho]\geq \theta\Opt[\varrho]=\theta\ln\epsilon$ when $0\leq\theta\leq1$. Now, to prove (\ref{wehaveq7}) is exactly the same as to prove that when $0\leq\rho<\vartheta^{-1}(\epsilon)\varrho$, no test for problem $(\cD_{g,\alpha}[\rho])$ with risk $\leq\epsilon$ is possible. Assuming, on the contrary, that $0\leq\rho<\vartheta(\epsilon)\varrho$ and $(\cD_{g,\alpha}[\rho])$ admits a test with risk $\leq\epsilon$; same as in the proof of Theorem \ref{the1}.ii, this implies that for every $x\in H_0[\rho]$ and  $y\in H_1[\rho]$, the Hellinger affinity of the distributions of observations associated with $x$ and $y$ does not exceed $2\sqrt{\epsilon}$, whence $\Opt[\rho]\leq \ln(2\sqrt{\epsilon})$. On the other hand, as we have seen, $\Opt[\rho]\geq{\rho\over\varrho}\ln \epsilon$, and we arrive at ${\rho\over\varrho}\ln\epsilon\leq \ln(2\sqrt{\epsilon})$, whence $\vartheta^{-1}(\epsilon)>\rho/\varrho\geq {\ln(2\sqrt{\epsilon})\over\ln\epsilon}=\vartheta^{-1}(\epsilon)$, which is impossible.
\paragraph{3$^0$.} Let now $\rho\in[\varrho,\rho_{\max}]$, so that problem  $(F_{g,\alpha}[\rho])$ is solvable with optimal value $\Opt[\rho]$; clearly, $\Opt[\rho]$ is a nonincreasing function of $\rho$, whence $\Opt[\rho]\leq \Opt[\varrho]=\epsilon$, Applying Proposition \ref{propnonstI} (with no Gaussian and Poisson factors and $a=0$) and recalling the origin of $\Opt[\rho]$,  we conclude that the risk of the simple test with the detector $\widehat{\phi}_{\rho}$ does not exceed $\exp\{\Opt[\rho]\}\leq\epsilon$. \qed

\end{document}